\documentclass{amsart}
\usepackage[dvips,final]{graphics}
\usepackage{array}
\usepackage{arydshln}
\usepackage[makeroom]{cancel}
 \usepackage[all]{xy}
 \usepackage{url}
\usepackage{multirow, blkarray}
\usepackage{booktabs}
\usepackage{textcomp}
 \usepackage[final]{epsfig}
 \usepackage{color}
\usepackage[T1]{fontenc}      
\usepackage[english,french]{babel}
\usepackage[utf8]{inputenc}
\usepackage{blindtext}

\usepackage{amsfonts,amscd,array, mathdots, epigraph}
\usepackage{amsmath}
\usepackage{amssymb}
\usepackage{amsthm}
\usepackage{mathrsfs}
\usepackage{stmaryrd}
\usepackage{slashbox}
\usepackage{diagbox}
\usepackage{enumitem}

\usepackage{ulem}
\usepackage{tikz}
\usepackage{xcolor}
\usepackage{multicol}
\definecolor{ufogreen}{rgb}{0.24, 0.82, 0.44}

\usepackage{hyperref}

\begin{document}


\newtheorem{theorem}{Théorème}[section]
\newtheorem{theore}{Théorème}
\newtheorem{definition}[theorem]{Définition}
\newtheorem{proposition}[theorem]{Proposition}
\newtheorem{corollary}[theorem]{Corollaire}
\newtheorem{con}{Conjecture}
\newtheorem*{remark}{Remarque}
\newtheorem*{remarks}{Remarques}
\newtheorem*{pro}{Problème}
\newtheorem*{examples}{Exemples}
\newtheorem*{example}{Exemple}
\newtheorem{lemma}[theorem]{Lemme}


\title[$\lambda$-quiddités et minoration de la taille maximale des irréductibles sur un corps fini]{Étude de quelques familles de $\lambda$-quiddités et minoration de la taille maximale des $\lambda$-quiddités irréductibles sur un corps fini}

\author{Flavien Mabilat}

\date{}

\keywords{modular group, $\lambda$-quiddity, squares in finite fields, generators}

\address{
}
\def\emailaddrname{{\itshape Courriel}}
\email{flavien.mabilat@univ-reims.fr}
\subjclass[2020]{20H30, 05E16, 20F05, 20H05}

\maketitle

\selectlanguage{french}

\begin{abstract}
Les $\lambda$-quiddités de taille $n$ sont des $n$-uplets d'éléments d'un ensemble fixé qui sont solutions d'une équation matricielle qui est fondamentale dans l'étude de la combinatoire du groupe modulaire et des frises de Coxeter. Pour en savoir plus sur celles-ci, on utilise une notion d'irréductibilité qui permet de restreindre l'étude à un nombre limité d'éléments qui doivent être déterminés pour chaque ensemble. Notre objectif ici est de définir plusieurs familles de $\lambda$-quiddités sur les corps finis et d'étudier les propriétés d'irréductibilité de celles-ci, dans l'optique notamment de minorer la taille maximale des éléments irréductibles sur les ensembles $\mathbb{F}_{q}$.
\\
\end{abstract}

\selectlanguage{english}
\begin{abstract}
$\lambda$-quiddities of size $n$ are $n$-tuples of elements from a fixed set that are solutions to a matrix equation which is fundamental in the study of the combinatorics of the modular group and Coxeter's friezes. To gain further insight into these objects, we use a notion of irreducibility, which allows restricting the study to a limited number of elements that must be determined for each set. Our goal here is to define several families of $\lambda$-quiddities over finite fields and to study their irreducibility properties, with the specific aim of establishing lower bounds on the maximal size of irreducible elements over $\mathbb{F}_{q}$.
\\
\end{abstract}

\selectlanguage{french}

\thispagestyle{empty}

\noindent {\bf Mots clés:} groupe modulaire; $\lambda$-quiddité; carrés dans les corps finis; générateurs.   
\\
\begin{flushright}
\og \textit{Ce que l'on conçoit bien s'énonce clairement.} \fg
\\Nicolas Boileau, \textit{L'Art poétique, Chant 1.}
\end{flushright}

\section{Introduction}
\label{intro}

Depuis le début des années soixante-dix, les frises de Coxeter, qui sont des arrangements de nombres dans le plan appartenant à un ensemble fixé et vérifiant une règle arithmétique, appelée règle unimodulaire, n'ont eu de cesse de fasciner des générations de mathématiciens. Introduites en 1973 par H.\ S.\ M.\ Coxeter pour l'étude du Pentagramma mirificum (voir \cite{Cox}), celles-ci se sont rapidement affranchies de leur rôle d'intermédiaire et sont maintenant devenues un objet d'étude à part entière du fait de leurs fascinantes propriétés et de leurs nombreuses applications. En particulier, l'une des questions les plus étudiées les concernant est de savoir comment construire toutes les frises possibles sur l'ensemble choisi. Lorsque ce dernier est un anneau commutatif unitaire $A$, cette interrogation est étroitement liée à la résolution des équations matricielles suivantes sur $A$ (voir par exemple \cite{CH} proposition 2.4) :
\[M_{n}(a_1,\ldots,a_n):=\begin{pmatrix}
    a_{n} & -1_{A} \\
    1_{A}  & 0_{A} 
   \end{pmatrix} \ldots \begin{pmatrix}
    a_{1} & -1_{A} \\
    1_{A}  & 0_{A} 
   \end{pmatrix}=-Id.\]
	
Par ailleurs, les matrices $M_{n}(a_1,\ldots,a_n)$ jouent également un rôle fondamental dans l'étude de la combinatoire du groupe modulaire. En effet, en utilisant les générateurs classiques de ce groupe, à savoir
\[S:=\begin{pmatrix}
    0 & -1 \\
    1  & 0 
   \end{pmatrix}~~~~{\rm et}~~~~T:=\begin{pmatrix}
    1 & 1 \\
    0  & 1 
   \end{pmatrix},\]
\noindent on montre que, pour tout $B \in Sl_{2}(\mathbb{Z})$, il existe $n>0$ et $(a_{1},\ldots,a_{n}) \in (\mathbb{N}^{*})^{n}$ tels que :
\[B=T^{a_{n}}S \ldots T^{a_{1}}S=M_{n}(a_{1},\ldots,a_{n}).\]
\noindent Cette écriture n'étant pas unique (on a par exemple $-Id=M_{3}(1,1,1)=M_{4}(1,2,1,2)$), on est naturellement amené à chercher une manière de construire et de compter les différents uplets d'entiers strictement positifs associés à une matrice, ou à un ensemble de matrices. Le premier résultat dans cette direction est le théorème d'Ovsienko (voir \cite{O} Théorèmes 1 et 2) qui donne une construction récursive et une description combinatoire des différentes écritures de $\pm Id$, complété par le théorème de Conley-Ovsienko (voir \cite{CO2} Théorème 2.2) qui donne le nombre exact de ces dernières. À la lueur de ces résultats, plusieurs pistes de généralisations sont envisageables. La plus naturelle est d'effectuer les mêmes tâches pour d'autres éléments du groupe modulaire, notamment pour les générateurs $S$ et $T$ (voir \cite{M6} pour une présentation des résultats obtenus pour ce problème). Cela dit, la direction la plus prometteuse est d'interpréter $\{-Id,Id\}$ comme le centre de $SL_{2}(\mathbb{Z})$, ce qui invite à considérer d'autres sous-groupes intéressants, comme les sous-groupes de congruence :
\[\hat{\Gamma}(N):=\{B \in SL_{2}(\mathbb{Z}),~B \equiv \pm Id [N]\}.\]
\noindent L'objectif est d'alors de connaître les $n$-uplets $(a_1,\ldots,a_n)$ d'entiers strictement positifs pour lesquels $M_{n}(a_{1},\ldots,a_{n}) \in \hat{\Gamma}(N)$, ce qui est équivalent à la résolution de l'équation $M_{n}(a_1,\ldots,a_n)=\pm Id$ sur $\mathbb{Z}/N\mathbb{Z}$.
\\
\\ \indent On est ainsi naturellement amené à considérer, pour un anneau commutatif unitaire $A$, l'équation généralisée ci-dessous :
\begin{equation}
\label{p}
\tag{$E_{A}$}
M_{n}(a_1,\ldots,a_n)=\pm Id.
\end{equation}

\noindent Cette équation est parfois nommée équation de Conway-Coxeter et les solutions de celle-ci sont en général appelées $\lambda$-quiddités sur $A$. Il existe plusieurs façons d'étudier ces dernières. On peut par exemple chercher s'il existe un nombre fini de $\lambda$-quiddités de taille fixée sur certains anneaux infinis (voir notamment \cite{CHP}), on peut essayer de fournir des constructions géométriques de toutes les solutions sur certains anneaux simples (par exemple sur $\mathbb{Z}$, voir \cite{CH} Théorème 7.3) ou bien encore on peut tenter de dénombrer les $\lambda$-quiddités sur des anneaux finis (voir par exemple \cite{CM}). Cela dit, l'un des principaux outils dont on dispose pour étudier ces dernières est une notion d'irréductibilité, rappelée dans la section suivante, qui, suivant une logique classique en algèbre, permet de restreindre l'étude à un nombre limité d'objets. 
\\
\\ \indent Cette idée est renforcée par le fait qu'il n'y a qu'un nombre fini de $\lambda$-quiddités irréductibles sur les anneaux finis (voir \cite{M} Théorème 1.1), ce qui incite naturellement à approfondir ces situations. On dispose notamment d'un certain nombre de résultats de classification pour des anneaux finis de petits cardinaux (voir notamment \cite{M1} section 4). Toutefois, il est illusoire d'espérer obtenir des résultats généraux de classification et assez peu pertinent de développer une collection de résultats isolés centrés sur des anneaux disparates, choisis plus pour leurs propriétés accommodantes que pour leur intérêt intrinsèque. Cela conduit par conséquent à rechercher des propriétés vérifiées par toutes les $\lambda$-quiddités irréductibles sur des familles d'anneaux ou de définir des familles de $\lambda$-quiddités possédant de bonnes propriétés d'irréductibilité. L'objectif de ce texte est de réfléchir sur ces deux pistes de travail pour les corps finis et notamment d'obtenir des minorations intéressantes de la taille maximale des irréductibles. Notons par ailleurs que les résultats que nous allons obtenir contribuent à l'étude des différentes écritures des éléments des sous-groupes de congruence du groupe modulaire abordée dans le paragraphe précédent via les corps $\mathbb{Z}/p\mathbb{Z}$.
\\
\\ \indent D'un point de vue pratique, on va commencer par regrouper dans la section suivante les définitions précises des notions dont on aura besoin dans la suite et les principaux résultats obtenus dans ce texte. Puis, après avoir donné dans la section \ref{pre} quelques éléments préliminaires, on reprendra dans la section \ref{mondyn} les éléments développés sur les solutions monomiales minimales et dynomiales minimales dans le cas des anneaux $\mathbb{Z}/N\mathbb{Z}$ et on les considérera dans le contexte des corps finis. Ensuite, dans la section \ref{tri}, on définira une nouvelle famille de solutions : les solutions trinomiales minimales. Grâce à elles, on obtiendra l'excellente minoration suivante : la taille maximale des $\lambda$-quiddités irréductibles sur $\mathbb{F}_{q}$ est supérieure à $\frac{3(q-1)}{2}$ si $q$ n'est pas une puissance de 2 et $3(q-1)$ sinon. Pour finir, on abordera rapidement d'autres familles dont la définition semble naturelle mais dont les propriétés ne sont malheureusement pas à la hauteur des espérances qu'elles pouvaient légitimement susciter.

\section{Définitions et résultats principaux}
\label{prin}

L'objectif de cette section est de définir précisément l'ensemble des concepts évoqués dans l'introduction et de regrouper les résultats les plus significatifs démontrés dans les parties suivantes.
\\
\\\indent Dans tout ce qui suit, tous les anneaux considérés sont commutatifs et unitaires. Soit $(A,+,\times)$ l'un de ceux-ci, on note ${\rm car}(A)$ la caractéristique de $A$, $0_{A}$ le neutre de $A$ pour $+$ et $1_{A}$ le neutre pour $\times$. On supposera systématiquement que $0_{A} \neq 1_{A}$. Si $k \in \mathbb{N}^{*}$, on pose $k_{A}:=\sum_{l=1}^{k} 1_{A}$. $U(A)$ désigne le groupe des éléments inversibles de $A$. Si $N \geq 2$, on note, pour $a \in \mathbb{Z}$, $\overline{a}:=a+N\mathbb{Z}$. Lorsque l'on considérera un corps fini défini par un quotient de la forme $\mathbb{F}_{q}:=\frac{(\mathbb{Z}/p\mathbb{Z})[X]}{\left\langle P(X)\right\rangle}$, on notera de la même façon un élément de $(\mathbb{Z}/p\mathbb{Z})[X]$ et sa classe dans $\mathbb{F}_{q}$. La partie entière d'un réel $x$ est notée $E[x]$ et $\mathbb{P}$ désigne l'ensemble des nombres premiers. On note $\varphi$ la fonction indicatrice d'Euler.
\\
\\On commence par une définition précise des $\lambda$-quiddités.

\begin{definition}[\cite{C}, définition 2.2]
\label{01}

Soient $A$ un anneau commutatif unitaire et $n \in \mathbb{N}^{*}$. On dit que le $n$-uplet $(a_{1},\ldots,a_{n})$ d'éléments de $A$ est une $\lambda$-quiddité sur $A$ de taille $n$ si $(a_{1},\ldots,a_{n})$ est une solution de \eqref{p}, c'est-à-dire si $M_{n}(a_{1},\ldots,a_{n})=\pm Id.$ S'il n'y a pas d'ambiguïté on parlera simplement de $\lambda$-quiddité.

\end{definition}

Dans la suite de ce texte, on réservera le terme $\lambda$-quiddité aux résultats généraux. En revanche, lorsque l'on parlera des familles particulières que l'on va étudier, on parlera de solutions (suivies d'un adjectif qualifiant la famille en question), et cela afin que les termes restent cohérents avec les travaux déjà menés sur les anneaux $\mathbb{Z}/p\mathbb{Z}$. Par ailleurs, les $\lambda$-quiddités étant invariantes par permutations circulaires, on considérera, si $(a_{1},\ldots,a_{n})$ est une $\lambda$-quiddité de taille $n$ et si $k \in \mathbb{N}$, $a_{n+k}:=a_{k}$.

\begin{definition}[\cite{C}, lemme 2.7]
\label{02}

Soient $(a_{1},\ldots,a_{n}) \in A^{n}$ et $(b_{1},\ldots,b_{m}) \in A^{m}$. On définit l'opération ci-dessous: \[(a_{1},\ldots,a_{n}) \oplus (b_{1},\ldots,b_{m})= (a_{1}+b_{m},a_{2},\ldots,a_{n-1},a_{n}+b_{1},b_{2},\ldots,b_{m-1}).\] Le $(n+m-2)$-uplet obtenu est appelé la somme de $(a_{1},\ldots,a_{n})$ avec $(b_{1},\ldots,b_{m})$.

\end{definition}

\begin{examples}

{\rm On donne ci-dessous quelques exemples de sommes dans $\mathbb{Z}/11\mathbb{Z}$ :
\begin{itemize}
\item $(\overline{5},\overline{2},\overline{3}) \oplus (\overline{-1},\overline{-2},\overline{6})=(\overline{0},\overline{2},\overline{2},\overline{-2})$;
\item $(\overline{2},\overline{1},\overline{3},\overline{3}) \oplus (\overline{4},\overline{5},\overline{0},\overline{1})=(\overline{3},\overline{1},\overline{3},\overline{7},\overline{5},\overline{0})$;
\item $(\overline{1},\overline{4},\overline{3},\overline{1}) \oplus (\overline{3},\overline{0},\overline{0},\overline{4},\overline{4})=(\overline{5},\overline{4},\overline{3},\overline{4},\overline{0},\overline{0},\overline{4})$;
\item $n \geq 2$, $(\overline{a_{1}},\ldots,\overline{a_{n}}) \oplus (\overline{0},\overline{0}) = (\overline{0},\overline{0}) \oplus (\overline{a_{1}},\ldots,\overline{a_{n}})=(\overline{a_{1}},\ldots,\overline{a_{n}})$.
\end{itemize}
}
\end{examples}

Cette opération est particulièrement utile pour l'étude de \eqref{p} car elle possède la propriété fondamentale suivante : si $(b_{1},\ldots,b_{m})$ est une $\lambda$-quiddité alors $(a_{1},\ldots,a_{n}) \oplus (b_{1},\ldots,b_{m})$ est une solution de \eqref{p} si et seulement si $(a_{1},\ldots,a_{n})$ est une $\lambda$-quiddité (voir \cite{C,WZ} et \cite{M1} proposition 3.7). En revanche, il convient de noter que $\oplus$ n'est ni commutative ni associative (voir \cite{WZ} exemple 2.1) et que les $k$-uplets d'éléments de $A$ n'ont pas d'inverse pour $\oplus$ lorsque $k \geq 3$.

\begin{definition}[\cite{C}, définition 2.5]
\label{03}

Soient $(a_{1},\ldots,a_{n}) \in A^{n}$ et $(b_{1},\ldots,b_{n}) \in A^{n}$. On note $(a_{1},\ldots,a_{n}) \sim (b_{1},\ldots,b_{n})$ si $(b_{1},\ldots,b_{n})$ est obtenu par permutations circulaires de $(a_{1},\ldots,a_{n})$ ou de $(a_{n},\ldots,a_{1})$.

\end{definition}

On peut aisément vérifier que $\sim$ est une relation d'équivalence sur les $n$-uplets d'éléments de $A$ (voir \cite{WZ} lemme 1.7) De plus, si $(a_{1},\ldots,a_{n}) \sim (b_{1},\ldots,b_{n})$ alors $(a_{1},\ldots,a_{n})$ est une $\lambda$-quiddité si et seulement si $(b_{1},\ldots,b_{n})$ en est une aussi (voir \cite{C} proposition 2.6).
\\
\\ \noindent Les définitions de $\oplus$ et de $\sim$ étant fixées, on peut maintenant présenter la notion d'irréductibilité annoncée.

\begin{definition}[\cite{C}, définition 2.9]
\label{04}

Une solution $(c_{1},\ldots,c_{n})$ avec $n \geq 3$ de \eqref{p} est dite réductible s'il existe une solution de \eqref{p} $(b_{1},\ldots,b_{l})$ et un $m$-uplet $(a_{1},\ldots,a_{m})$ d'éléments de $A$ tels que : \begin{itemize}
\item $(c_{1},\ldots,c_{n}) \sim (a_{1},\ldots,a_{m}) \oplus (b_{1},\ldots,b_{l})$;
\item $m \geq 3$ et $l \geq 3$.
\end{itemize}
Une solution est dite irréductible si elle n'est pas réductible.

\end{definition}

\begin{remark} 

{\rm $(0_{A},0_{A})$ est toujours solution de \eqref{p}. Cependant, celle-ci n'est jamais considérée comme étant irréductible.}

\end{remark}

Cette notion d'irréductibilité nous permet donc de réduire l'étude de l'équation \eqref{p} à la bonne connaissance des $\lambda$-quiddités irréductibles. Dans le cas des anneaux finis, l'intérêt de cette réduction est accru par le théorème ci-dessous :

\begin{theorem}[\cite{M}, Théorème 1.1]
\label{05}

Soit $A$ un anneau commutatif unitaire fini. Il n'y a qu'un nombre fini de $\lambda$-quiddités irréductibles sur $A$. Notons $\ell_{A}$ la taille maximale des $\lambda$-quiddités irréductibles sur $A$. On a :
\begin{itemize}
\item Si $car(A)=2$ alors $4 \leq \ell_{A} \leq \frac{\left|SL_{2}(A)\right|}{\left|A\right|}+2$.
\item Si $car(A) \neq 2$ alors ${\rm max}(4,car(A)) \leq \ell_{A} \leq \frac{\left|SL_{2}(A)\right|}{2\left|A\right|}+2$.
\end{itemize}

\end{theorem}

Notons que les minorations générales fournies par ce théorème sont optimales. En effet, $\ell_{\mathbb{Z}/2\mathbb{Z}}=4$ et $\ell_{\mathbb{Z}/6\mathbb{Z}}=6$. En revanche, elles semblent assez éloignées des valeurs exactes de $\ell_{A}$ dans la plupart des cas.
\\
\\ \indent Lorsque l'on combine ces éléments à l'impossibilité concrète d'obtenir des listes complètes pour des anneaux de grand cardinaux, on est naturellement amené à chercher des propriétés vérifiées par l'ensemble des solutions irréductibles ou des familles de $\lambda$-quiddités possédant de bonnes propriétés d'irréductibilité. Ici, on va chercher à combiner ces deux idées en étudiant des familles particulières et en en déduisant des informations sur $\ell_{A}$. Pour approfondir cette piste, on va bien évidemment commencer par définir des familles particulières de solutions.
\\
\\Au cours des différents travaux menés sur $(E_{\mathbb{Z}/p\mathbb{Z}})$, on a défini les deux familles ci-dessous.
\begin{itemize}
\item Soit $\overline{k} \in \mathbb{Z}/p\mathbb{Z}$. La solution $\overline{k}$-monomiale minimale de $(E_{\mathbb{Z}/p\mathbb{Z}})$ est la solution dont toutes les composantes sont égales à $\overline{k}$ et de taille minimale pour cette propriété.
\item Soit $\overline{k} \in \mathbb{Z}/p\mathbb{Z}$. La solution $\overline{k}$-dynomiale minimale de $(E_{\mathbb{Z}/p\mathbb{Z}})$ est le $n$-uplet $(\overline{k},\overline{-k},\ldots,\overline{k},\overline{-k})$ avec $n$ le plus entier pour lequel il existe une solution de cette forme.
\end{itemize}

\noindent Ces deux familles se prolongent naturellement aux anneaux $\mathbb{Z}/N\mathbb{Z}$. Dans le cas des solutions monomiales minimales, on dispose de très nombreux résultats d'irréductibilité (voir par exemple \cite{M1,M2,M4}). En particulier, si $p$ est premier, on sait que toutes les solutions monomiales minimales non nulles de $(E_{\mathbb{Z}/p\mathbb{Z}})$ sont irréductibles (\cite{M1} Théorème 3.16). En ce qui concerne les solutions dynomiales minimales, on dispose du résultat suivant. Soient $p$ un nombre premier supérieur à 5 et $\overline{k} \in \mathbb{Z}/p\mathbb{Z}-\{\overline{0}\}$. Si $\overline{k}^{2}+\overline{8}$ n'est pas un carré dans $\mathbb{Z}/p\mathbb{Z}$ alors la solution $\overline{k}$-dynomiale minimale de $(E_{\mathbb{Z}/p\mathbb{Z}})$ est irréductible (\cite{M2} Théorème 2.7).
\\
\\ \indent On va commencer par généraliser ces résultats sur tous les corps finis. Pour commencer, on étendra, dans la section \ref{mondyn}, la définition des solutions monomiales minimales à tous les corps $\mathbb{F}_{q}$ et on obtiendra l'irréductibilité de toutes les solutions non nulles de cette famille. Ensuite, on généralisera le concept de solutions dynomiales minimales sur un anneau $A$ en donnant la définition suivante : la solution $(a,b)$-dynomiale minimale de \eqref{p} est la solution de la forme $(a,b,\ldots,a,b)$ de taille minimale. Après avoir démontré quelques résultats techniques généralisant ce qui avait été fait sur $\mathbb{Z}/p\mathbb{Z}$, on obtiendra le théorème d'irréductibilité énoncé ci-après. Soient $\mathbb{K}$ un corps fini de caractéristique différente de 2 et $(a,b) \in (\mathbb{K}^{*})^{2}$ avec $a \neq b$ et $a,b \neq \pm 1_{\mathbb{K}}$. Si $\Delta_{1}:=a^{2}+4_{\mathbb{K}}ab^{-1}(ab^{-1}-1_{\mathbb{K}})$ et $\Delta_{2}:=b^{2}+4_{\mathbb{K}}a^{-1}b(a^{-1}b-1_{\mathbb{K}})$ ne sont pas des carrés sur $\mathbb{K}$ alors la solution $(a,b)$-dynomiale minimale de $(E_{\mathbb{K}})$ est irréductible (Théorème \ref{35}).
\\
\\ \indent Enfin, dans la section \ref{tri}, on introduira une nouvelle classe de solutions : les solutions trinomiales minimales. Ce sont les solutions de la forme $(u,u^{-1},u^{-1},\ldots,u,u^{-1},u^{-1})$ et de taille minimale. Ces solutions sont celles qui nous seront le plus utile, grâce notamment au théorème suivant (Théorème \ref{42}) :

\begin{theorem}
\label{06}

Soient $\mathbb{K}$ un corps fini et $u \in U(\mathbb{K})$. Soient $m$ ta taille de la solution $u$-trinomiale minimale et $o(u)$ l'ordre de $u$ dans $\mathbb{K}^{*}$. On a trois cas :
\begin{itemize}
\item si ${\rm car}(\mathbb{K})=2$ alors $m=3o(u)$;
\item si ${\rm car}(\mathbb{K})\neq 2$ et $o(u)$ est pair alors $m=\frac{3o(u)}{2}$;
\item si ${\rm car}(\mathbb{K})\neq 2$ et $o(u)$ est impair alors $m=3o(u)$.
\end{itemize}

\noindent De plus, si $u \neq \{\pm 1_{\mathbb{K}}\}$, la solution $u$-trinomiale minimale de $(E_{\mathbb{K}})$ est irréductible si et seulement si $u$ n'est pas racine des polynômes $X^{2l}+X^{l+1}-1_{\mathbb{K}}$ et $X^{2l}-X^{l+1}-1_{\mathbb{K}}$ pour $1 \leq l \leq E\left[\frac{m}{6}\right]$.

\end{theorem}

Grâce à celui-ci, on obtiendra une minoration assez fine pour $\ell_{\mathbb{K}}$, donnée dans l'énoncé ci-dessous (Théorèmes \ref{410} et \ref{420}) et démontrée dans la section \ref{tri}.

\begin{theorem}
\label{07}

Soit $\mathbb{K}$ un corps fini de cardinal $q$.
\begin{itemize}
\item si ${\rm car}(\mathbb{K}) \neq 2$ alors $\ell_{\mathbb{K}} \geq \frac{3(q-1)}{2}$;
\item si ${\rm car}(\mathbb{K})=2$ alors $\ell_{\mathbb{K}} \geq 3(q-1)$.
\end{itemize}

\end{theorem}

\section{Résultats préliminaires}
\label{pre}

L'objectif de cette section est de rappeler quelques résultats classiques qui nous seront utiles dans toute la suite.

\subsection{Continuant}

On va donner dans cette sous-partie la définition et quelques propriétés du continuant qui nous seront d'une inappréciable utilité dans tout ce qui va suivre.

\begin{definition}
\label{11}

Soient $A$ un anneau commutatif unitaire, $n \in \mathbb{N}^{*}$ et $(a_{1},\ldots,a_{n})$ un $n$-uplet d'éléments de $A$. On pose $K_{-1}:=0_{A}$, $K_{0}:=1_{A}$ et on note \[K_{n}(a_{1},\ldots,a_{n}):=
\left|
\begin{array}{cccccc}
a_1&1_{A}&&&\\[4pt]
1_{A}&a_{2}&1_{A}&&\\[4pt]
&\ddots&\ddots&\!\!\ddots&\\[4pt]
&&1_{A}&a_{n-1}&\!\!\!\!\!1_{A}\\[4pt]
&&&\!\!\!\!\!1_{A}&\!\!\!\! a_{n}
\end{array}
\right|.\] 
\noindent $K_{n}(a_{1},\ldots,a_{n})$ est le continuant de $a_{1},\ldots,a_{n}$. 

\end{definition}

Notons qu'en développant ce déterminant suivant la première ou la dernière colonne, on obtient, pour $n \geq 1$, les formules ci-dessous :
\begin{itemize}
\item $K_{n}(a_{1},\ldots,a_{n})=a_{1}K_{n-1}(a_{2},\ldots,a_{n})-K_{n-2}(a_{3},\ldots,a_{n})$;
\item $K_{n}(a_{1},\ldots,a_{n})=a_{n}K_{n-1}(a_{1},\ldots,a_{n-1})-K_{n-2}(a_{1},\ldots,a_{n-2})$.
\\
\end{itemize} 

Pour calculer les continuants, on peut utiliser l'algorithme d'Euler. On peut décrire celui-ci de la façon suivante (voir par exemple \cite{CO}). $K_{n}(a_{1},\ldots,a_{n})$ est la somme de tous les produits possibles de $a_{1},\ldots,a_{n}$, dans lesquels un nombre quelconque de paires disjointes de termes consécutifs est supprimé. Chacun de ces produits étant multiplié par $-1_{A}$ puissance le nombre de paires supprimées. Par exemple, cela donne :
\begin{itemize}
\item $K_{2}(a_{1},a_{2})=a_{1}a_{2}-1_{A}$;
\item $K_{3}(a_{1},a_{2},a_{3})=a_{1}a_{2}a_{3}-a_{3}-a_{1}$;
\item $K_{4}(a_{1},a_{2},a_{3},a_{4})=a_{1}a_{2}a_{3}a_{4}-a_{3}a_{4}-a_{1}a_{4}-a_{1}a_{2}+1_{A}$.
\\
\end{itemize}

Ces polynômes sont très utiles pour les questions considérées dans ce texte car ils permettent d'exprimer les coordonnées des matrices $M_{n}(a_{1},\ldots,a_{n})$. En effet, on dispose de l'égalité suivante (voir \cite{CO}).

\begin{proposition}
\label{12}

Soient $A$ un anneau commutatif unitaire et $(a_{1},\ldots,a_{n}) \in A^{n}$, avec $n \in \mathbb{N}^{*}$.
\[M_{n}(a_{1},\ldots,a_{n})=\begin{pmatrix}
    K_{n}(a_{1},\ldots,a_{n}) & -K_{n-1}(a_{2},\ldots,a_{n}) \\
    K_{n-1}(a_{1},\ldots,a_{n-1})  & -K_{n-2}(a_{2},\ldots,a_{n-1}) 
   \end{pmatrix}.\]

\end{proposition}

\begin{proof}

Ce résultat est assez connu. Cependant, étant donné son importance dans la suite, on va le démontrer rapidement. 
\\
\\On raisonne par récurrence sur $n$. On a $K_{-1}=0_{A}$, $K_{0}=1_{A}$ et $K_{1}(a_{1})=a_{1}$, donc le résultat est vrai pour $n=1$.
\\
\\Supposons qu'il existe $n \in \mathbb{N}^{*}$ tel que la formule est vraie au rang $n$. On a :
\begin{eqnarray*}
M &=& M_{n+1}(a_{1},\ldots,a_{n+1}) \\ 
  &=& M_{1}(a_{n+1})M_{n}(a_{1},\ldots,a_{n}) \\
  &=& \begin{pmatrix}
    a_{n+1} & -1_{A} \\
    1_{A}  & 0_{A} 
   \end{pmatrix}\begin{pmatrix}
    K_{n}(a_{1},\ldots,a_{n}) & -K_{n-1}(a_{2},\ldots,a_{n}) \\
    K_{n-1}(a_{1},\ldots,a_{n-1})  & -K_{n-2}(a_{2},\ldots,a_{n-1}) 
   \end{pmatrix} \\
	&=& \begin{pmatrix}
    a_{n+1}K_{n}(a_{1},\ldots,a_{n})-K_{n-1}(a_{1},\ldots,a_{n-1}) & -a_{n+1}K_{n-1}(a_{2},\ldots,a_{n})+K_{n-2}(a_{2},\ldots,a_{n-1}) \\
    K_{n}(a_{1},\ldots,a_{n})  & -K_{n-1}(a_{2},\ldots,a_{n}) 
   \end{pmatrix} \\
	&=& \begin{pmatrix}
    K_{n+1}(a_{1},\ldots,a_{n+1}) & -K_{n}(a_{2},\ldots,a_{n+1}) \\
    K_{n}(a_{1},\ldots,a_{n})  & -K_{n-1}(a_{2},\ldots,a_{n}) 
   \end{pmatrix}. \\
\end{eqnarray*}		

\noindent Par récurrence, le résultat est vrai.

\end{proof}

\noindent On aura également besoin des formules regroupées dans les deux propositions suivantes.

\begin{proposition}
\label{13}

Soient $A$ un anneau commutatif unitaire, $u \in U(A)$, $n \in \mathbb{N}^{*}$ et $(a_{1},\ldots,a_{n}) \in A^{n}$.
\\
\\i) $K_{n}(a_{1},\ldots,a_{n})=K_{n}(a_{n},\ldots,a_{1})$.
\\ii) On suppose que $n=2l$ est pair. $K_{n}(u a_{1},u^{-1} a_{2},\ldots,u a_{2l-1}, u^{-1} a_{2l})=K_{n}(a_{1},\ldots,a_{2l})$.
\\iii) On suppose que $n=2l+1$ est impair. $K_{n}(u a_{1},u^{-1} a_{2},\ldots,u a_{2l-1}, u^{-1} a_{2l},u a_{2l+1})=u K_{n}(a_{1},\ldots,a_{n})$.

\end{proposition}

\begin{proof}

i) On procède par récurrence sur $n$. Pour $n=1$, l'égalité est immédiate et pour $n=2$ on a, pour tout $(a_{1},a_{2}) \in A^{2}$, $K_{2}(a_{1},a_{2})=a_{1}a_{2}-1_{A}=K_{2}(a_{2},a_{1})$. Supposons qu'il existe un $n \geq 2$ tel que, pour tout $(a_{1},\ldots,a_{n}) \in A^{n}$, $K_{n}(a_{1},\ldots,a_{n})=K_{n}(a_{n},\ldots,a_{1})$ et $K_{n-1}(a_{1},\ldots,a_{n-1})=K_{n-1}(a_{n-1},\ldots,a_{1})$. Soit $(a_{1},\ldots,a_{n+1}) \in A^{n+1}$. 
\begin{eqnarray*}
K &=& K_{n+1}(a_{1},\ldots,a_{n+1}) \\
  &=& a_{1}K_{n}(a_{2},\ldots,a_{n+1})-K_{n-1}(a_{3},\ldots,a_{n+1}) \\
	&=& a_{1}K_{n}(a_{n+1},\ldots,a_{2})-K_{n-1}(a_{n+1},\ldots,a_{3})~~({\rm hypoth\grave{e}se~de~r\acute{e}currence}) \\
	&=& K_{n+1}(a_{n+1},\ldots,a_{1}).
\end{eqnarray*}

\noindent Par récurrence, le résultat est démontré.
\\
\\ii) et iii). On propose deux preuves. La première est purement matricielle et repose sur l'égalité suivante :
\[\begin{pmatrix} 
u^{-1} & 0_{A} \\ 
0_{A} & 1_{A} 
\end{pmatrix}  
\begin{pmatrix} 
 u x & -1_{A} \\
1_{A} & 0_{A} 
\end{pmatrix}
\begin{pmatrix} 
u^{-1} y & -1_{A} \\ 
1_{A} & 0_{A} 
\end{pmatrix}
\begin{pmatrix} 
u & 0_{A} \\
0_{A} & 1_{A}
 \end{pmatrix}=\begin{pmatrix} 
                x & -1_{A} \\ 
								1_{A} & 0_{A} 
								\end{pmatrix}
               \begin{pmatrix} 
							 y & -1_{A} \\ 
							 1_{A} & 0_{A}
							\end{pmatrix}.\]

\noindent La seconde utilise l'algorithme d'Euler.
\\
\\ \uwave{Preuve matricielle :} Soit $n=2l$. Avec la formule ci-dessus, on a :
\[M_{n}(a_{1},\ldots,a_{n})=\begin{pmatrix} 
	u & 0_{A} \\
	 0_{A} & 1_{A} 
	\end{pmatrix}M_{n}(u a_{1},u^{-1} a_{2},\ldots,u a_{2l-1},u^{-1} a_{2l})\begin{pmatrix} 
	u^{-1} & 0_{A} \\ 
	0_{A} & 1_{A}
	\end{pmatrix}.\]

\noindent Donc, on a les égalités suivantes :
\begin{eqnarray*}
M &=& M_{n}(u a_{1},u^{-1} a_{2},\ldots,u a_{2l-1},u^{-1} a_{2l}) \\
  &=& \begin{pmatrix} 
	  K_{n}(u a_{1},u^{-1} a_{2},\ldots,u a_{2l-1}, u^{-1} a_{2l}) & -K_{n-1}(u^{-1} a_{2},\ldots,u a_{2l-1}, u^{-1} a_{2l}) \\ 
	  K_{n-1}(u a_{1},u^{-1} a_{2},\ldots,u a_{2l-1}) & -K_{n-2}(u^{-1} a_{2},\ldots,u a_{2l-1})
	    \end{pmatrix} \\
	&=& \begin{pmatrix} 
	  u^{-1} & 0_{A} \\ 
	  0_{A} & 1_{A}
	    \end{pmatrix}\begin{pmatrix} 
	  K_{n}(a_{1},\ldots,a_{n}) & -K_{n-1}(a_{2},\ldots,a_{n}) \\ 
	  K_{n-1}(a_{1},\ldots,a_{n-1}) & -K_{n-2}(a_{2},\ldots,a_{n-1})
	    \end{pmatrix}\begin{pmatrix} 
	  u & 0_{A} \\ 
	  0_{A} & 1_{A}
	    \end{pmatrix} \\
	&=& \begin{pmatrix} 
	  K_{n}(a_{1},\ldots,a_{n}) & -u^{-1}K_{n-1}(a_{2},\ldots,a_{n}) \\ 
	  u K_{n-1}(a_{1},\ldots,a_{n-1}) & -K_{n-2}(a_{2},\ldots,a_{n-1})
	    \end{pmatrix}.
\end{eqnarray*}
\noindent Ainsi, on a :
\begin{itemize}
\item $K_{n}(u a_{1},u^{-1} a_{2},\ldots,u a_{2l-1}, u^{-1} a_{2l})=K_{n}(a_{1},\ldots,a_{2l})$.
\item $K_{n-1}(u a_{1},u^{-1} a_{2},\ldots,u a_{2l-1})=u K_{n-1}(a_{1},\ldots,a_{n-1})$, ce qui en remplaçant $l$ par $l+1$ donne la formule souhaitée.
\\
\end{itemize}

\noindent \uwave{Preuve avec l'algorithme d'Euler :} Soient $n=2l$, $(a_{1},\ldots,a_{n}) \in A^{n}$ et $u \in U(A)$. Par l'algorithme d'Euler, on obtient que $K_{2l}(u a_{1},u^{-1} a_{2},\ldots,u a_{2l-1}, u^{-1} a_{2l})$ est une somme de produits de $a_{j}u$ ($j$ impair), de $a_{k}u^{-1}$ ($k$ pair) et éventuellement de $-1_{A}$. Dans chaque terme d'une de ces sommes, pour chaque $j$ impair $a_{j}u$ est suivi d'un $a_{k}u^{-1}$ avec $k$ pair. On fait alors la manipulation suivante :\[a_{j}ua_{k}u^{-1}=a_{j}a_{k}.\]
\noindent Avec cette manipulation, on a : \[K_{2l}(u a_{1},u^{-1} a_{2},\ldots,u a_{2l-1}, u^{-1} a_{2l})=K_{2l}(a_{1},\ldots,a_{2l}).\]

\noindent Soient $n=2l+1$, $u \in U(A)$, et $(a_{1},\ldots,a_{n}) \in A^{n}$. Par l'algorithme d'Euler, on sait que le continuant $K_{2l+1}(u a_{1},u^{-1} a_{2},\ldots,u a_{2l-1}, u^{-1} a_{2l},ua_{2l+1})$ est une somme de produits. Pour chaque terme de la somme, on a deux choix :
\begin{itemize}
\item le terme est de la forme $a_{j}u$ avec $j$ impair;
\item le terme est un produit de $2h+1$ éléments de la forme $a_{j}u$ ou $a_{k}u^{-1}$ (avec éventuellement $-1_{A}$) avec $h+1$ indices $j$ impair. Pour chaque $m$ pair $a_{m}$ est précédé d'un $a_{k}$ avec $k$ impair. On fait alors la manipulation suivante : $a_{k}u^{-1}a_{m}u=a_{k}a_{m}$.
\end{itemize}

\noindent Dans chaque terme du produit, il reste donc un seul $u$. Cela donne donc 
\[K_{2l+1}(u a_{1},u^{-1} a_{2},\ldots,u a_{2l-1}, u^{-1} a_{2l},ua_{2l+1})=u K_{n}(a_{1},\ldots,a_{n}).\]

\end{proof}

\begin{proposition}
\label{132}

Soient $A$ un anneau commutatif unitaire, $n \in \mathbb{N}^{*}$ et $x \in A$.
\[K_{n}(x,\ldots,x)=\sum_{i=0}^{E[\frac{n}{2}]} (-1)^{i} {n-i \choose i} x^{n-2i}.\]

\end{proposition}

\begin{proof}

On prouve cette formule par récurrence. Pour une preuve détaillée, voir par exemple \cite{M2} lemme 3.7.

\end{proof}

\noindent On termine cette section par le rappel du résultat très utile suivant :

\begin{lemma}
\label{131}

Soit $A$ un anneau commutatif unitaire. Soient $n \geq 1$ et $(a_{1},\ldots,a_{n}) \in A^{n}$ tels que $K_{n}(a_{1},\ldots,a_{n})=\epsilon \in \{\pm 1_{A}\}$. On pose $x=\epsilon K_{n-1}(a_{2},\ldots,a_{n})$ et $y=\epsilon K_{n-1}(a_{1},\ldots,a_{n-1})$. 
\[M_{n+2}(x,a_{1},\ldots,a_{n},y)=-\epsilon Id.\]

\end{lemma}

\begin{proof}

\begin{eqnarray*}
M &=& M_{n+2}(x,a_{1},\ldots,a_{n},y) \\
  &=& \begin{pmatrix}
    K_{n+2}(x,a_{1},\ldots,a_{n},y) & -K_{n+1}(a_{1},\ldots,a_{n},y) \\
    K_{n+1}(x,a_{1},\ldots,a_{n})  & -K_{n}(a_{1},\ldots,a_{n}) 
   \end{pmatrix} \\
	&=& \begin{pmatrix}
    K_{n+2}(x,a_{1},\ldots,a_{n},y) & -K_{n+1}(a_{1},\ldots,a_{n},y) \\
    K_{n+1}(x,a_{1},\ldots,a_{n})  & -\epsilon
   \end{pmatrix}. \\
\end{eqnarray*}

\noindent Comme $\epsilon^{2}=1_{A}$, $K_{n-1}(a_{2},\ldots,a_{n})=\epsilon x$ et $K_{n-1}(a_{1},\ldots,a_{n-1})=\epsilon y$.
\\
\\On a, $K_{n+1}(a_{1},\ldots,a_{n},y)=yK_{n}(a_{1},\ldots,a_{n})-K_{n-1}(a_{1},\ldots,a_{n-1})=\epsilon y-\epsilon y=0_{A}$. On procède de la même façon pour le calcul de $K_{n+1}(x,a_{1},\ldots,a_{n})$. Comme $M \in SL_{2}(A)$, $K_{n+2}(x,a_{1},\ldots,a_{n},y)=-\epsilon$. Ainsi, $M=-\epsilon Id$.

\end{proof}

\subsection{Carrés modulo un nombre premier $p$}
\label{carré}

On va s'intéresser aux carrés dans les corps finis, c'est-à-dire aux éléments qui peuvent s'écrire comme le carré d'un autre. On commence par deux résultats de comptage.

\begin{proposition}[\cite{G}, proposition VII.51]
\label{141}

Soit $q>1$ la puissance d'un nombre premier $p$. Si $p=2$ alors tous les éléments de $\mathbb{F}_{q}$ sont des carrés. Si $p \neq 2$, il y a $\frac{q+1}{2}$ carrés dans $\mathbb{F}_{q}$.

\end{proposition}

\begin{proposition}
\label{14}

Soient $q>1$ la puissance d'un nombre premier impair et $a \in \mathbb{F}_{q}^{*}$. 
\begin{itemize}
\item Si $-a$ est un carré dans $\mathbb{F}_{q}$ alors il y a $\frac{q+1}{2}$ éléments $x$ dans $\mathbb{F}_{q}$ tel que $x^{2}+a$ est un carré dans $\mathbb{F}_{q}$.
\item Si $-a$ n'est pas un carré dans $\mathbb{F}_{q}$ alors il y a $\frac{q-1}{2}$ éléments $x$ dans $\mathbb{F}_{q}$ tel que $x^{2}+a$ est un carré dans $\mathbb{F}_{q}$.
\end{itemize}

\end{proposition}

\begin{proof}

Soit $a \in \mathbb{F}_{q}^{*}$. On commence par compter le nombre de solutions sur $\mathbb{F}_{q}$ de l'équation $x^{2}-y^{2}=a$, c'est-à-dire $(x-y)(x+y)=a$.
\\
\\Soient $\mathcal{A}:=\{(x,y) \in (\mathbb{F}_{q})^{2},~xy=a\}$ et $\mathcal{B}:=\{(x,y) \in (\mathbb{F}_{q})^{2},~(x-y)(x+y)=a\}$. Comme $q$ est impair, $2_{\mathbb{F}_{q}}$ est inversible. On considère les deux applications suivantes :
\[\begin{array}{ccccc} 
\varphi_{a} : & \mathcal{A} & \longrightarrow & \mathcal{B} \\
  & (x,y) & \longmapsto & (2_{\mathbb{F}_{q}}^{-1}(x+y),2_{\mathbb{F}_{q}}^{-1}(x-y))  \\
\end{array}~~{\rm et}~~\begin{array}{ccccc} 
\psi_{a}  : & \mathcal{B} & \longrightarrow & \mathcal{A} \\
 & (x,y) & \longmapsto & (x+y,x-y)  \\
\end{array}.\]
\noindent $\varphi_{a}$ et $\psi_{a}$ sont des bijections réciproques. Donc, ${\rm card}(\mathcal{A})={\rm card}(\mathcal{B})$.
\\
\\Or, ${\rm card}(\mathcal{A})=q-1$. En effet, puisque $a$ est non nul, $xy=a$ si et seulement si $y \neq 0_{A}$ et $x=ay^{-1}$. Par conséquent, ${\rm card}(\mathcal{B})=q-1$.
\\
\\On revient maintenant au problème initial. Soit $\mathcal{C}:=\{x \in \mathbb{F}_{q},~x^{2}+a~{\rm est~un~carr\acute{e}~dans}~\mathbb{F}_{q}\}$. Le couple $(y,x) $ est dans $\mathcal{B}$ si et seulement si $x \in \mathcal{C}$ et $x^{2}+a=y^{2}$.
\\
\\Supposons que $-a$ est un carré dans $\mathbb{F}_{q}$. Il existe $z \in \mathbb{F}_{q}$ tel que $-a=z^{2}=(-z)^{2}$. Ainsi, $(0_{\mathbb{F}_{q}},z)$ et $(0_{\mathbb{F}_{q}},-z)$ sont dans $\mathcal{B}$. Si $x \in \mathcal{C}-\{\pm z\}$, il existe exactement deux éléments de $\mathbb{F}_{q}^{*}$, $y_{x}$ et $-y_{x}$, tels que $x^{2}+a=(\pm y_{x})^{2}$. Donc, 
\[\mathcal{B}=\{(0_{\mathbb{F}_{q}},z),(0_{\mathbb{F}_{q}},-z)\} \bigsqcup_{x \in \mathcal{C}-\{\pm z\}} \{(y_{x},x),(-y_{x},x)\}.\] 
\noindent Ainsi, $q-1={\rm card}(\mathcal{B})=2+2({\rm card}(\mathcal{C})-2)$. Donc, ${\rm card}(\mathcal{C})=\frac{q+1}{2}$.
\\
\\Supposons que $-a$ n'est pas un carré dans $\mathbb{F}_{q}$. Si $x \in \mathcal{C}$, il existe exactement deux éléments de $\mathbb{F}_{q}^{*}$, $y_{x}$ et $-y_{x}$, tels que $x^{2}+a=(\pm y_{x})^{2}$. Donc, 
\[\mathcal{B}= \bigsqcup_{x \in \mathcal{C}} \{(y_{x},x),(-y_{x},x)\}.\] 
\noindent Ainsi, $q-1={\rm card}(\mathcal{B})=2{\rm card}(\mathcal{C})$. Donc, ${\rm card}(\mathcal{C})=\frac{q-1}{2}$.

\end{proof}

On va maintenant se placer sur les corps $\mathbb{Z}/p\mathbb{Z}$. Soit $p$ un nombre premier impair et $a$ est un entier premier avec $p$, on note $\left(\dfrac{a}{p}\right)$ le symbole de Legendre, c'est-à-dire :
\[\left(\dfrac{a}{p}\right)=\left\{
    \begin{array}{ll}
        1 & \mbox{si } \overline{a}~{\rm est~un~carr\acute{e}~dans}~\mathbb{Z}/p\mathbb{Z}; \\
        -1 & \mbox{sinon}.
    \end{array}
\right.  \\ \]

\noindent Le symbole de Legendre vérifie les propriétés suivantes :

\begin{proposition}[\cite{G}, proposition XII.20]
\label{15}
Soient $p$ un nombre premier impair et $a$ et $b$ deux entiers premiers avec $p$. On a :
\\i)~(critère d'Euler)~$\left(\dfrac{a}{p}\right) \equiv a^{\frac{p-1}{2}}~[p]$;
\\ii)~(multiplicativité)~$\left(\dfrac{ab}{p}\right)=\left(\dfrac{a}{p}\right)\left(\dfrac{b}{p}\right)$;
\\iii)~$\left(\dfrac{-1}{p}\right)=(-1)^{\frac{p-1}{2}}$.

\end{proposition}

\noindent Pour savoir si un élément donné de $\mathbb{Z}/p\mathbb{Z}$ est un carré, on dispose des deux résultats majeurs ci-dessous :

\begin{theorem}[Loi de réciprocité quadratique de Gauss; \cite{G}, Théorème XII.25]
\label{16}

Soient $p$ et $q$ deux nombres premiers impairs distincts. On a : \[\left(\dfrac{p}{q}\right)\left(\dfrac{q}{p}\right)= (-1)^{\frac{p-1}{2}\frac{q-1}{2}}.\]

\end{theorem}

\begin{theorem}[Loi complémentaire; \cite{G}, proposition XII.27]
\label{17}

Soit $p$ un nombre premier impair. On a $\left(\dfrac{2}{p}\right)= (-1)^{\frac{p^{2}-1}{8}}$, c'est-à-dire $\left(\dfrac{2}{p}\right)=1$ si et seulement si $p \equiv \pm 1 [8]$.

\end{theorem}

\noindent On peut notamment en déduire le résultat ci-dessous :

\begin{proposition}
\label{18}

Soit $p$ un nombre premier. 
\[\overline{3}~{\rm est~un~carr\acute{e}~dans}~\mathbb{Z}/p\mathbb{Z} \Longleftrightarrow \left\{
    \begin{array}{ll}
        p=2;  \\
        p=3; \\
				p \equiv \pm 1~[12]. 
    \end{array}
\right.  \\ \]
	
\end{proposition}

\begin{proof}

Si $p=2$ alors $\overline{3}=\overline{1}=\overline{1}^{2}$ et si $p=3$ alors $\overline{3}=\overline{0}=\overline{0}^{2}.$ On suppose maintenant $p$ supérieur à 5. D'après la loi de réciprocité quadratique de Gauss on a : \[\left(\dfrac{3}{p}\right)=\left(\dfrac{p}{3}\right)(-1)^{\frac{p-1}{2}\frac{3-1}{2}}=\left(\dfrac{p}{3}\right)(-1)^{\frac{p-1}{2}}.\]
Comme $p$ est premier, $p \equiv 1~[3]$ ou $p \equiv -1~[3]$. De plus, on a $\left(\dfrac{1}{3}\right)=1$ et $\left(\dfrac{-1}{3}\right)=-1$. On distingue les deux cas.
\begin{itemize}
\item Supposons que $p=3k+1$ avec $k$ un entier naturel pair.
\[\left(\dfrac{3}{p}\right)=1 \Longleftrightarrow \left(\dfrac{1}{3}\right)(-1)^{\frac{3k}{2}}=1 \Longleftrightarrow (-1)^{\frac{3k}{2}}=1 \Longleftrightarrow 4~{\rm divise}~k \Longleftrightarrow p \equiv 1~[12].\]
\item Supposons que $p=3k-1$ avec $k$ un entier naturel pair.
\[\left(\dfrac{3}{p}\right)=1 \Longleftrightarrow \left(\dfrac{-1}{3}\right)(-1)^{\frac{3k-2}{2}}=1 \Longleftrightarrow (-1)^{\frac{3k}{2}}=1 \Longleftrightarrow 4~{\rm divise}~k \Longleftrightarrow p \equiv -1~[12].\]
\end{itemize}

\end{proof}

Pour conclure, on va donner un résultat, basé sur le calcul de la norme des éléments des corps finis (voir \cite{Sa} chapitre II, section 2.6, proposition 1), qui permet d'utiliser les éléments précédents pour étudier les carrés dans les corps $\mathbb{F}_{q}$.

\begin{proposition}
\label{19}

Soient $q>1$ la puissance d'un nombre premier impair $p$ et $a \in \mathbb{F}_{q}$. $a$ est un carré dans $\mathbb{F}_{q}$ si et seulement si $a^{\frac{q-1}{p-1}}$ est un carré dans $\mathbb{Z}/p\mathbb{Z}$.   

\end{proposition}

\section{Solutions monomiales minimales et dynomiales minimales}
\label{mondyn}

L'objectif de cette section est de généraliser les résultats déjà connus sur les solutions monimiales minimales et dynomiales minimales.

\subsection{Solutions monomiales minimales}

Les solutions monomiales minimales ont été définies dans le cas des anneaux $\mathbb{Z}/N\mathbb{Z}$ (voir \cite{M1} section 3.3) et possèdent d'excellentes propriétés. En particulier, toutes les solutions monomiales minimales non nulles sur $\mathbb{Z}/p\mathbb{Z}$ (avec $p$ premier) sont irréductibles. On peut facilement généraliser ce résultat à n'importe quel  corps, puisque cette propriété découle de l'intégrité de $\mathbb{Z}/p\mathbb{Z}$ (notons que, pour un anneau commutatif unitaire fini, être intègre est équivalent au fait d'être un corps). On donne rapidement ci-dessous les grandes lignes de cette généralisation.
\\
\\On commence par définir formellement la famille en question.

\begin{definition}
\label{21}

Soient $A$ un anneau commutatif unitaire fini et $x \in A$. On appelle solution $(x,A)$-monomiale minimale une solution de \eqref{p} dont toutes les composantes sont égales à $x$ et qui est de taille minimale. On appelle solution monomiale minimale sur $A$ une solution $(x,A)$-monomiale minimale pour un certain $x \in A$.

\end{definition} 

\begin{remark}
{\rm 
On pourrait bien entendu se contenter de parler de solutions $x$-monomiale minimale de \eqref{p}. Toutefois, le concept de solutions monomiales minimales étant avant tout utilisé dans le cas des anneaux $\mathbb{Z}/N\mathbb{Z}$, il semble préférable de rajouter dans le cas général la précision de l'anneau pour bien indiquer qu'on ne travaille pas sur les classes d'entiers.
}

\end{remark}

Notons qu'une telle solution existe toujours puisque $M_{1}(x)$ est d'ordre fini dans $SL_{2}(A)/\{-Id,Id\}$, puisque ce groupe est de cardinal fini. 

\begin{examples}
{\rm
\begin{itemize}
\item Soit $\mathbb{K}$ un corps fini, $(1_{\mathbb{K}},1_{\mathbb{K}},1_{\mathbb{K}})$ est une solution monomiale minimale de $(E_{\mathbb{K}})$.
\item $(\overline{3},\overline{3},\overline{3},\overline{3},\overline{3},\overline{3})$ est une solution monomiale minimale de $(E_{\mathbb{Z}/6\mathbb{Z}})$.
\item $(X,X,X,X,X)$ est une solution monomiale minimale sur $\mathbb{F}_{4}=\frac{\mathbb{F}_{2}[X]}{\left\langle X^{2}+X+\overline{1}\right\rangle}$.
\end{itemize}
}
\end{examples}

Le résultat clé pour la suite est la proposition ci-dessous dont la preuve est identique à celle menée dans le cas des anneaux $\mathbb{Z}/N\mathbb{Z}$ (voir \cite{M1} proposition 3.15 et la preuve du lemme \ref{33} qui reprend toutes les idées de la preuve).

\begin{proposition}
\label{22}

Soient $A$ un anneau commutatif unitaire fini, $(a,b,x) \in A^{3}$ et $n$ un entier supérieur à 3. Si $(a,x,\ldots,x,b) \in A^{n}$ est une solution de \eqref{p} alors $a=b$ et $a(a-x)=0_{A}$.

\end{proposition}

\noindent Ce résultat permet d'avoir le théorème suivant.

\begin{theorem}
\label{23}

Soit $\mathbb{K}$ un corps fini. Les solutions monomiales minimales non nulles de $(E_{\mathbb{K}})$ sont irréductibles.

\end{theorem}

\begin{proof}

Soit $x \in \mathbb{K}$, $x \neq 0_{\mathbb{K}}$. Soit $n$ la taille de la solution $(x,\mathbb{K})$-monomiale minimale de $(E_{\mathbb{K}})$. Supposons par l'absurde que cette dernière soit réductible. Il existe un $l$-uplet $(a_{1},\ldots,a_{l})$ et une solution $(b_{1},\ldots,b_{l'})$ tels que $l,l' \geq 3$, $l+l'=n+2$ et
\[(x,\ldots,x)=(a_{1},\ldots,a_{l}) \oplus (b_{1},\ldots,b_{l'})=(a_{1}+b_{l'},a_{2},\ldots,a_{l-1},a_{l}+b_{1},b_{2},\ldots,b_{l'-1}).\]

\noindent Ainsi, $b_{2}=\ldots=b_{l'-1}=x$. Donc, par la proposition précédente, $b_{1}=b_{l'}:=b$ et $b(b-x)=0_{\mathbb{K}}$. Comme $\mathbb{K}$ est un corps, $\mathbb{K}$ est intègre et donc $b=x$ ou $b=0_{\mathbb{K}}$. Si $b=x$ alors nécessairement, par minimalité de $n$, $l' \geq n$ et donc $l \leq 2$ ce qui est absurde. Si $b=0_{\mathbb{K}}$ alors $(b_{1},\ldots,b_{l'}) \sim \underbrace{(x,\ldots,x)}_{l'-2} \oplus (0_{\mathbb{K}},0_{\mathbb{K}},0_{\mathbb{K}},0_{\mathbb{K}})$. Donc, le $l'-2$-uplet ne contenant que des $x$ est une solution de $(E_{\mathbb{K}})$. Cela implique, par minimalité de $n$, $l'-2 \geq n$ et donc $l'=n+2$ et $l=0$ ce qui est absurde.
\\
\\Donc, la solution $(x,\mathbb{K})$-monomiale minimale non nulle de $(E_{\mathbb{K}})$ est irréductible.

\end{proof}

On sait aussi que, pour tout anneau $A$, la solution $(2_{A},A)$-monomiale minimale est irréductible de taille ${\rm car}(A)$ (c'est l'origine de la minoration du théorème \ref{05}, voir \cite{M}). De plus, dans le cas des anneaux $\mathbb{Z}/N\mathbb{Z}$, on dispose également de nombreuses autres propriétés (voir par exemple \cite{M4,M5}). En particulier, on a classifié les entiers monomialement irréductibles, qui sont les entiers $N \geq 2$ pour lesquels toutes les solutions monomiales minimales non nulles de $(E_{\mathbb{Z}/N\mathbb{Z}})$ sont irréductibles (voir \cite{OEIS} A350242). L'entier $N$ est monomialment irréductible si et seulement si $N$ est premier ou est égal à 4, 6, 8, 12 ou 24 (voir \cite{M7} Théorème 2.6). On reviendra brièvement sur les solutions monomiales dans la dernière section. 

\subsection{Solutions dynomiales généralisées}
\label{dyn}

On va maintenant généraliser la famille des solutions dynomiales minimales, qui ont été originellement définies dans le cas des anneaux $\mathbb{Z}/p\mathbb{Z}$ et qui consistaient en les solutions de la forme $(\overline{k},-\overline{k},\ldots,\overline{k},-\overline{k})$ de taille minimale. Cette généralisation va prendre deux formes. D'une part, on va considérer l'ensemble des anneaux, et, d'autre part, on va élargir les composantes possibles. 

\subsubsection{Définition et premières propriétés} 

On commence par définir formellement la famille que nous allons étudier.

\begin{definition}
\label{31}

Soient $A$ un anneau commutatif unitaire fini et $(a,b) \in A^{2}$, $a \neq b$. On appelle solution $(a,b)$-dynomiale minimale une solution de \eqref{p} de la forme $(a,b,\ldots,a,b)$ et qui est de taille minimale. On appelle solution dynomiale minimale sur $A$ une solution $(a,b)$-dynomiale minimale pour un certain couple $(a,b) \in A^{2}$ ($a \neq b$).

\end{definition} 

Notons qu'une telle solution existe toujours puisque $M_{1}(a,b)$ est d'ordre fini dans $SL_{2}(A)/\{-Id,Id\}$, puisque ce groupe est de cardinal fini. De plus, la solution $(b,a)$-dynomiale minimale de \eqref{p} est équivalente à la solution $(a,b)$-dynomiale minimale de \eqref{p}. Par conséquent, ces deux solutions ont la même taille et sont toutes les deux réductibles ou bien toutes les deux irréductibles.

\begin{examples}
{\rm On donne ci-dessous quelques exemples de solutions dynomiales minimales.
\begin{itemize}
\item Soit $\mathbb{K}$ un corps fini, $(1_{\mathbb{K}},2_{\mathbb{K}},1_{\mathbb{K}},2_{\mathbb{K}})$ est une solution dynomiale minimale de $(E_{\mathbb{K}})$.
\item $(X,X+\overline{1},X,X+\overline{1},X,X+\overline{1})$ est une solution dynomiale minimale sur le corps $\mathbb{F}_{4}=\frac{(\mathbb{Z}/2\mathbb{Z})[X]}{\left\langle X^{2}+X+\overline{1}\right\rangle}$.
\item $(\overline{3},\overline{4},\overline{3},\overline{4},\overline{3},\overline{4},\overline{3},\overline{4})$ est la solution $(\overline{3},\overline{4})$-dynomiale minimale sur $\mathbb{Z}/7\mathbb{Z}$.
\end{itemize}
}
\end{examples} 

On commence par le résultat ci-dessous qui donne, pour presque tous les anneaux, des solutions dynomiales minimales irréductibles.

\begin{proposition}
\label{311}

Soit $A$ un anneau commutatif unitaire fini tel que $U(A) \neq \{\pm 1_{A}\}$. On suppose qu'il existe $u \in U(A)-\{\pm 1_{A}\}$ tel que $u \neq u^{-1}$. La solution $(u,u^{-1})$-dynomiale minimale de \eqref{p} est irréductible de taille 6.

\end{proposition}

\begin{proof}

Par la proposition \ref{12}, on a :
\begin{eqnarray*}
M &=& M_{2}(u,u^{-1})^{3} \\
  &=& M_{6}(u,u^{-1},u,u^{-1},u,u^{-1}) \\
  &=& \begin{pmatrix}
    K_{6}(u,u^{-1},u,u^{-1},u,u^{-1}) & -K_{5}(u^{-1},u,u^{-1},u,u^{-1}) \\
    K_{5}(u,u^{-1},u,u^{-1},u)  & -K_{4}(u^{-1},u,u^{-1},u) 
   \end{pmatrix} \\
	&=& \begin{pmatrix}
    u(u^{-1})u(u^{-1})u(u^{-1})-5_{A}u(u^{-1})u(u^{-1})+6_{A}u(u^{-1})-1_{A} & -u^{-1}u(u^{-1})u(u^{-1}) \\
		                                                                         & +4_{A}u^{-1}u(u^{-1})-3_{A}u^{-1} \\
																																						 & \\
    u(u^{-1})u(u^{-1})u-4_{A}u(u^{-1})u+3_{A}u  & -u^{-1}u(u^{-1})u+3_{A}u^{-1}u-1_{A} \\
		\end{pmatrix} \\
 &=& \begin{pmatrix}
    1_{A}-5_{A}+6_{A}-1_{A} & -u^{-1}+4_{A}u^{-1}-3_{A}u^{-1} \\
    u-4_{A}u+3_{A}u  & -1_{A}+3_{A}-1_{A}
   \end{pmatrix} \\
	&=& \begin{pmatrix}
    1_{A} & 0_{A} \\
    0_{A}  & 1_{A}
   \end{pmatrix}.\\
\end{eqnarray*}

\noindent Donc, la taille $l$ de la solution $(u,u^{-1})$-dynomiale minimale de \eqref{p} est égale à $2m$ avec $m$ divisant 3, puisque $m$ est l'ordre de $M_{2}(u,u^{-1})$ dans $SL_{2}(A)/\{-Id,Id\}$. Comme $M_{2}(u,u^{-1}) \neq \pm Id$, $m \neq 1$. Ainsi, $l=6$.
\\
\\Si la solution $(u,u^{-1})$-dynomiale minimale de \eqref{p} était réductible alors elle serait nécessairement la somme de deux solutions de taille 4, puisque $u, u^{-1} \neq \pm 1_{A}$. Cela impliquerait qu'il existe des solutions de la forme $(x,u,u^{-1},y)$ ou $(x,u^{-1},u,y)$. En particulier, $1_{A}=u(u^{-1})=0_{A}$ ou $2_{A}$, ce qui est absurde. Donc, cette solution est irréductible.

\end{proof}

\begin{examples}
{\rm On donne ci-dessous quelques applications de cette proposition.
\begin{itemize}
\item $A=\mathbb{F}_{4}=\frac{(\mathbb{Z}/2\mathbb{Z})[X]}{\left\langle X^{2}+X+\overline{1}\right\rangle}$. $(X,X+\overline{1},X,X+\overline{1},X,X+\overline{1})$ est une solution irréductible sur $A$.
\item $A=\mathbb{Z}/9\mathbb{Z}$. $(\overline{2},\overline{5},\overline{2},\overline{5},\overline{2},\overline{5})$ est une solution irréductible sur $A$.
\end{itemize}
}
\end{examples} 

\noindent Avec, ce résultat, on peut obtenir une très légère amélioration de la minoration de $\ell_{A}$.

\begin{corollary}
\label{312}

Soit $A$ un anneau commutatif unitaire fini tel que $U(A) \neq \{\pm 1_{A}\}$. On suppose qu'il existe $u \in U(A)-\{\pm 1_{A}\}$ tel que $u \neq u^{-1}$. On a $\ell_{A} \geq {\rm max}(6,{\rm car}(A))$.

\end{corollary}

Notons que cette inégalité possède des cas d'égalités, par exemple lorsque $A=\mathbb{Z}/5\mathbb{Z}$ ou $A=\mathbb{Z}/6\mathbb{Z}$ (voir \cite{M1}).

\begin{corollary}
\label{313}

Soit $\mathbb{K}$ un corps fini de cardinal $q>3$ et de caractéristique $p$. Si $p=2$ alors $(E_{\mathbb{K}})$ a au moins $q-2$ solutions irréductibles de taille 6. Si $p \neq 2$ alors $(E_{\mathbb{K}})$ a au moins $q-3$ solutions irréductibles de taille 6.

\end{corollary}

\begin{proof}

Comme $q>3$, $U(\mathbb{K}) \neq \{\pm 1_{\mathbb{K}}\}$. De plus, si $u \in \mathbb{K}$ est tel que $u=u^{-1}$ alors $u=\pm 1_{\mathbb{K}}$. Donc, pour tout $u \in \mathbb{K}-\{0_{\mathbb{K}},\pm 1_{\mathbb{K}}\}$, la solution $(u,u^{-1})$-dynomiale minimale est irréductible, ce qui donne $q-2$ solutions irréductibles de taille 6 si $p=2$ et $q-3$ si $p>2$.

\end{proof}

\subsubsection{Résultats d'irréductibilité}

On cherche dans cette sous-partie à obtenir un certain nombre de conditions suffisantes pour l'irréductibilité des solutions dynomiales minimales. Plus précisément, notre objectif principal ici est de généraliser le résultat suivant obtenu dans le cadre restreint utilisé auparavant.

\begin{theorem}[\cite{M2} Théorème 2.7]
\label{32}

Soient $p$ un nombre premier supérieur à 5 et $\overline{k} \in \mathbb{Z}/p\mathbb{Z}$. On suppose que les deux conditions suivantes sont vérifiées : 
\begin{itemize}
\item  $\overline{k} \neq \overline{0}$;
\item $\overline{k}^{2}+\overline{8}$ n'est pas un carré dans $\mathbb{Z}/p\mathbb{Z}$.
\end{itemize}
La solution $(\overline{k},-\overline{k})$-dynomiale minimale de $(E_{\mathbb{Z}/p\mathbb{Z}})$ est irréductible.

\end{theorem}

\noindent Pour faire cela, on va démontrer l'analogue ci-dessous de la proposition \ref{22}.

\begin{lemma}
\label{33}

Soient $A$ un anneau commutatif unitaire fini, $n \in \mathbb{N}^{*}$, $n \geq 5$ et $(x,y,a,b) \in A^{4}$ avec $a,b$ inversibles et $a \neq b$.
\\
\\i) Si $(x,a,b,\ldots,a,b,y) \in A^{n}$ ($n$ pair) est solution de \eqref{p} alors $y=ab^{-1}x$ et $0_{A}=x(b-x)$. Si $(x,b,a,\ldots,b,a,y) \in A^{n}$ ($n$ pair) est solution de \eqref{p} alors $y=a^{-1}bx$ et $0_{A}=x(a-x)$.
\\
\\ii) Si $(x,a,b,\ldots,a,b,a,y) \in A^{n}$ ($n$ impair) est solution de \eqref{p} alors $y=x$ et $x(a-ab^{-1}x)=1_{A}-ab^{-1}$. Si $(x,b,a,\ldots,b,a,b,y) \in A^{n}$ ($n$ impair) est solution de \eqref{p} alors $y=x$ et $x(b-a^{-1}bx)=1_{A}-a^{-1}b$.
	
\end{lemma}	
	
\begin{proof}

i) Supposons que $n$ est pair. Comme $(x,a,b,\ldots,a,b,y)$ est une solution de \eqref{p}, il existe $\epsilon$ dans $\{-1_{A},1_{A}\}$ tel que :
\[\epsilon Id=M_{n}(x,a,b,\ldots,a,b,y)=\begin{pmatrix}
   K_{n}(x,a,b,\ldots,a,b,y)   & -K_{n-1}(a,b,\ldots,a,b,y) \\
   K_{n-1}(x,a,b,\ldots,a,b) & -K_{n-2}(a,b,\ldots,a,b)
\end{pmatrix}.\]
\noindent Ainsi, \[K_{n-1}(x,a,b,\ldots,a,b)=K_{n-1}(a,b,\ldots,a,b,y)=0_{A}~~{\rm et}~~K_{n-2}(a,b,\ldots,a,b)=-\epsilon.\] 

\noindent Donc, en développant le déterminant, on obtient :
\begin{eqnarray*}
0_{A} &=& K_{n-1}(x,a,b,\ldots,a,b) \\
      &=& x K_{n-2}(a,b,\ldots,a,b)-K_{n-3}(b,a,b,\ldots,a,b) \\
      &=& -\epsilon x-K_{n-3}(b,a,b,\ldots,a,b). \\
\end{eqnarray*}

Ainsi, comme $\epsilon^{2}=1_{A}$, on a $x=-\epsilon K_{n-3}(b,a,b,\ldots,a,b)$. De même, on a :
\begin{eqnarray*}
0_{A} &=& K_{n-1}(a,b,\ldots,a,b,y) \\
      &=& y K_{n-2}(a,b,\ldots,a,b)-K_{n-3}(a,b,\ldots,a,b,a) \\
      &=& -\epsilon y-K_{n-3}(a,b,\ldots,a,b,a). \\
\end{eqnarray*}
Il en découle :
\begin{eqnarray*}
y &=& -\epsilon K_{n-3}(a,b,\ldots,a,b,a) \\
  &=& -\epsilon ab^{-1} K_{n-3}(b,a,b,\ldots,a,b)~~~({\rm proposition}~\ref{13}~{\rm iii)~avec}~n-3~{\rm impair}) \\
	&=&  ab^{-1}x.\\
\end{eqnarray*}
						
\noindent De plus, on dispose des deux égalités ci-dessous : 
\[-\epsilon=K_{n-2}(a,b,\ldots,a,b)=a K_{n-3}(b,a,b,\ldots,a,b)-K_{n-4}(a,b,\ldots,a,b),\]							
et 
\begin{eqnarray*}
M &=& M_{n-2}(a,b,\ldots,a,b) \\
  &=& \begin{pmatrix}
   K_{n-2}(a,b,\ldots,a,b)  & -K_{n-3}(b,a,b,\ldots,a,b) \\
   K_{n-3}(a,b,\ldots,a,b,a) & -K_{n-4}(b,a\ldots,b,a)
\end{pmatrix} \\
  &=& \begin{pmatrix}
   K_{n-2}(a,b,\ldots,a,b)  & -K_{n-3}(b,a,b\ldots,a,b) \\
   ab^{-1} K_{n-3}(b,a,b\ldots,a,b) & -K_{n-4}(b,a\ldots,b,a)
\end{pmatrix} \in SL_{2}(A).
\end{eqnarray*}

\noindent On en déduit la relation suivante : 
\begin{eqnarray*}
1_{A} &=& -K_{n-2}(a,b,\ldots,a,b)K_{n-4}(b,a,\ldots,b,a)+ab^{-1} K_{n-3}(b,a,b\ldots,a,b)^{2} \\
      &=& \epsilon K_{n-4}(b,a,\ldots,b,a)+ab^{-1} x^{2}.
\end{eqnarray*}
\noindent Donc, par la proposition \ref{13} i), on a :

\[K_{n-4}(a,b,\ldots,a,b)=K_{n-4}(b,a,\ldots,b,a)=\epsilon(1_{A}-ab^{-1} x^{2}).\]
																																												
\noindent Ainsi, on a les égalités suivantes :
\begin{eqnarray*}
-\epsilon &=& a K_{n-3}(b,a,b,\ldots,a,b)-K_{n-4}(a,b,\ldots,a,b) \\
          &=& -\epsilon ax-\epsilon(1_{A}-ab^{-1} x^{2}) \\ 
					&=& -\epsilon x(a-ab^{-1} x)-\epsilon.
\end{eqnarray*}

\noindent Donc, $y=ab^{-1}x$, et, comme $\epsilon$ est inversible, on a $0_{A}=x(a-ab^{-1}x)=ab^{-1}x(b-x)$. Comme $ab^{-1}$ est inversible, on obtient $0_{A}=x(b-x)$.
\\
\\Comme $a$ et $b$ ont un rôle interchangeable, on obtient aisément le résultat suivant. Si $(x,b,a,\ldots,b,a,y)$ est un $n$-uplet d'éléments de $A$ solution de \eqref{p} alors $y=a^{-1}bx$ et $0_{A}=x(a-x)$.
\\
\\ii) Supposons que $n$ est impair. Comme $(x,a,b,\ldots,a,b,a,y)$ est une solution de \eqref{p}, il existe $\epsilon$ dans $\{-1_{A},1_{A}\}$ tel que :
\[\epsilon Id=M_{n}(x,a,b,\ldots,a,b,a,y)=\begin{pmatrix}
   K_{n}(x,a,b,\ldots,a,b,a,y)   & -K_{n-1}(a,b,\ldots,a,b,a,y) \\
   K_{n-1}(x,a,b,\ldots,a,b,a) & -K_{n-2}(a,b,\ldots,a,b,a)
\end{pmatrix}.\]
\noindent Donc, \[K_{n-1}(a,b,\ldots,a,b,a,y)=K_{n-1}(x,a,b,\ldots,a,b,a)=0_{A}~~{\rm et}~~K_{n-2}(a,b,\ldots,a,b,a)=-\epsilon.\] 
\noindent En développant le déterminant, on obtient les égalités ci-dessous :
\begin{eqnarray*}
0_{A} &=& K_{n-1}(x,a,b,\ldots,a,b,a) \\
      &=& x K_{n-2}(a,b,\ldots,a,b,a)-K_{n-3}(b,a,\ldots,b,a) \\
      &=& -\epsilon x-K_{n-3}(b,a,\ldots,b,a) \\
			&=& -\epsilon x-K_{n-3}(a,b,\ldots,a,b)~~({\rm proposition}~\ref{13}~{\rm i)}).\\
\end{eqnarray*}
Ainsi, comme $\epsilon^{2}=1_{A}$, on a $x=-\epsilon K_{n-3}(a,b,\ldots,a,b)$. De même, on a : 
\begin{eqnarray*}
0_{A} &=& K_{n-1}(a,b,\ldots,a,b,a,y) \\
      &=& y K_{n-2}(a,b,\ldots,a,b,a)-K_{n-3}(a,b,\ldots,a,b) \\
      &=& -\epsilon y-K_{n-3}(a,b,\ldots,a,b).
\end{eqnarray*}
Donc, $y=-\epsilon K_{n-3}(a,b,\ldots,a,b)=x$.
\\
\\De plus, on a :
\begin{eqnarray*}
-\epsilon &=& K_{n-2}(a,b,\ldots,a,b,a) \\
          &=& a K_{n-3}(b,a,\ldots,b,a)-K_{n-4}(a,b,\ldots,a,b,a) \\
					&=& a K_{n-3}(a,b,\ldots,a,b)-K_{n-4}(a,b,\ldots,a,b,a).
\end{eqnarray*}
\noindent et, par la proposition \ref{13} i) :

\begin{eqnarray*}
M &=& M_{n-2}(a,b,\ldots,a,b,a) \\
  &=& \begin{pmatrix}
   K_{n-2}(a,b,\ldots,a,b,a)  & -K_{n-3}(b,a,\ldots,b,a) \\
   K_{n-3}(a,b,\ldots,a,b) & -K_{n-4}(b,a,b,\ldots,a,b)
\end{pmatrix} \\
  &=& \begin{pmatrix}
   K_{n-2}(a,b,\ldots,a,b,a)  & -K_{n-3}(a,b,\ldots,a,b) \\
   K_{n-3}(a,b,\ldots,a,b)    & -K_{n-4}(b,a,b,\ldots,a,b)
\end{pmatrix} \in SL_{2}(A).
\end{eqnarray*}

\noindent On obtient les relations ci-dessous.
\begin{eqnarray*}
1_{A} &=& -K_{n-2}(a,b,\ldots,a,b,a)K_{n-4}(b,a,b,\ldots,a,b)+K_{n-3}(a,b,\ldots,a,b)^{2} \\
      &=& \epsilon K_{n-4}(b,a,b,\ldots,a,b)+x^{2}.
\end{eqnarray*}
\noindent En particulier, on en déduit, puisque $n-4$ est impair :
\begin{eqnarray*}
K_{n-4}(a,b,\ldots,a,b,a) &=& ab^{-1} K_{n-4}((a^{-1}b)a,(ab^{-1})b,(a^{-1}b)a,\ldots,(ab^{-1})b,(a^{-1}b)a) \\
                          &=& ab^{-1}K_{n-4}(b,a,b,\ldots,a,b) \\
													&=& \epsilon ab^{-1} (1_{A}-x^{2}).
\end{eqnarray*}					
								
\noindent Par conséquent, on a :
\begin{eqnarray*}
-\epsilon &=& a K_{n-3}(a,b,\ldots,a,b)-K_{n-4}(a,b,\ldots,a,b,a) \\
          &=& -\epsilon ax-\epsilon ab^{-1} (1_{A}-x^{2}) \\ 
			    &=& -\epsilon ax-\epsilon ab^{-1}+\epsilon ab^{-1} x^{2} \\
					&=& -\epsilon [x(a-ab^{-1}x)+ab^{-1}].
\end{eqnarray*}

\noindent Donc, puisque $\epsilon^{2}=1$, on a les égalités $x(a-ab^{-1}x)=1_{A}-ab^{-1}$ et $x=y$.
\\
\\Comme $a$ et $b$ ont un rôle interchangeable, on obtient aisément le résultat suivant. Si $(x,b,a,\ldots,b,a,b,y)$ est un $n$-uplet d'éléments de $A$ solution de \eqref{p} alors $y=x$ et $x(b-a^{-1}bx)=1_{A}-a^{-1}b$.

\end{proof}

\begin{lemma}
\label{34}

Soient $A$ un anneau commutatif unitaire fini et $(a,b) \in (U(A)-\{\pm 1_{A}\})^{2}$ avec $a \neq b$. On suppose que les conditions suivantes sont vérifiées :
\begin{itemize}
\item le polynôme $X(b-X)$ a pour seules racines $0_{A}$ et $b$ sur $A$;
\item le polynôme $X(a-ab^{-1}X)-(1_{A}-ab^{-1})$ n'a pas de racine sur $A$;
\item le polynôme $X(b-a^{-1}bX)-(1_{A}-a^{-1}b)$ n'a pas de racine sur $A$.
\end{itemize}

\noindent La solution $(a,b)$-dynomiale minimale de \eqref{p} est irréductible.
	
\end{lemma}	

\begin{proof}

Soient $A$ un anneau commutatif unitaire fini et $(a,b) \in (U(A)-\{\pm 1_{A}\})^{2}$ avec $a \neq b$. On note $n$ la taille de la solution $(a,b)$-dynomiale minimale de \eqref{p}. Supposons par l'absurde que celle-ci soit réductible. Il existe deux solutions $(x,c_{1},\ldots,c_{l-2},y) \in A^{l}$ et $(x',d_{1},\ldots,d_{l'-2},y') \in A^{l'}$ tels que $l,l' \geq 3$ et
\begin{eqnarray*}
(a,b,\ldots,a,b) &\sim& (x,c_{1},\ldots,c_{l-2},y) \oplus (x',d_{1},\ldots,d_{l'-2},y') \\
                 &=& (x+y',c_{1},\ldots,c_{l-2},y+x',d_{1},\ldots,d_{l'-2}).
\end{eqnarray*}
\noindent Si $l=3$ ou $l'=3$ alors $a=\pm 1_{A}$ ou $b=\pm 1_{A}$, ce qui n'est pas le cas. Si $l=4$ ou $l'=4$ alors soit $ab=0_{A}$ ce qui est impossible puisque $a$ et $b$ sont inversibles soit $ab=2_{A}$. Si $ab=2_{A}$ alors $(a,b,a,b)$ est une solution de \eqref{p}. C'est donc la solution $(a,b)$-dynomiale minimale de \eqref{p} et elle est irréductible (puisqu'une solution de taille 4 est réductible si et seulement si elle contient $\pm 1_{A}$), ce qui contredit l'hypothèse de départ. Donc, on suppose $l,l' \geq 5$. On distingue deux cas :
\\
\\i) Si $l'$ est pair. Comme $n$ est pair et $n=l+l'-2$, $l$ est pair. On distingue deux cas :
\begin{itemize}
\item $(a,b,\ldots,a,b)=(x+y',c_{1},\ldots,c_{l-2},y+x',d_{1},\ldots,d_{l'-2})$. Dans ce cas, $(x',d_{1},\ldots,d_{l'-2},y')$ est égal à $(x',a,b,\ldots,a,b,y')$. Par le lemme précédent, on a $y'=ab^{-1}x'$ et $0_{A}=x'(b-x')$. Or, par hypothèse, $0_{A}=X(b-X)$ a pour seules racines $0_{A}$ et $b$ sur $A$. Ainsi, $x'=y'=0_{A}$ ou $x'=b$ et $y'=a$. Si $x'=y'=0_{A}$ alors le $l'-2$-uplet $(a,b,\ldots,a,b)$ est solution (car $(x',a,b,\ldots,a,b,y') \sim (a,b,\ldots,a,b) \oplus (0_{A},0_{A},0_{A},0_{A})$). Par minimalité de $n$, on a $l' \geq n+2$ et donc $l=0$, ce qui est absurde. Si $x'=b$ et $y'=a$ alors le $l'$-uplet $(b,a\ldots,b,a)$ est solution. Par minimalité de $n$, on a $l' \geq n$ et donc $l \leq 2$, ce qui est absurde. 
\\
\item $(b,a,\ldots,b,a)=(x+y',c_{1},\ldots,c_{l-2},y+x',d_{1},\ldots,d_{l'-2})$. Dans ce cas, $(x,c_{1},\ldots,c_{l-2},y)$ est égal à $(x,a,b,\ldots,a,b,y)$ et on procède comme dans le cas précédent.
\\
\end{itemize}

\noindent Si $l'$ est impair. Comme $n$ est pair et $n=l+l'-2$, $l$ est impair. On distingue deux cas :
\begin{itemize}
\item $(a,b,\ldots,a,b)=(x+y',c_{1},\ldots,c_{l-2},y+x',d_{1},\ldots,d_{l'-2})$. Dans ce cas, $(x,c_{1},\ldots,c_{l-2},y)$ est égal à $(x,b,a,b,\ldots,a,b,y)$ (et $(x',d_{1},\ldots,d_{l'-2},y')=(x',b,a,b,\ldots,a,b,y')$). Par le lemme précédent, $y=x$ et $x(b-a^{-1}bx)=1_{A}-a^{-1}b$. Or, par hypothèse, $X(b-a^{-1}bX)-(1_{A}-a^{-1}b)$ n'a pas de racine sur $A$. Donc, on arrive à une absurdité.
\\
\item $(b,a,\ldots,b,a)=(x+y',c_{1},\ldots,c_{l-2},y+x',d_{1},\ldots,d_{l'-2})$. Dans ce cas, $(x,c_{1},\ldots,c_{l-2},y)$ est égal à $(x,a,b,a,\ldots,b,a,y)$ (et $(x',d_{1},\ldots,d_{l'-2},y')=(x',a,b,a,\ldots,b,a,y')$). Par le lemme précédent, $y=x$ et $x(a-ab^{-1}x)=1_{A}-ab^{-1}$. Or, par hypothèse, le polynôme $X(a-ab^{-1}X)-(1_{A}-ab^{-1})$ n'a pas de racine sur $A$. Donc, on arrive à une absurdité.
\\
\end{itemize}

\noindent Ainsi, dans tous les cas, on arrive à une absurdité. Par conséquent, la solution $(a,b)$-dynomiale minimale de \eqref{p} est irréductible.

\end{proof}

\noindent On peut maintenant démontrer le résultat principal de cette section.

\begin{theorem}
\label{35}

Soient $\mathbb{K}$ un corps fini de caractéristique différente de 2 et $(a,b) \in (\mathbb{K}^{*})^{2}$, avec $a \neq b$ et $a,b \neq \pm 1_{\mathbb{K}}$. Si $\Delta_{1}:=a^{2}+4_{\mathbb{K}}ab^{-1}(ab^{-1}-1_{\mathbb{K}})$ et $\Delta_{2}:=b^{2}+4_{\mathbb{K}}a^{-1}b(a^{-1}b-1_{\mathbb{K}})$ ne sont pas des carrés sur $\mathbb{K}$ alors la solution $(a,b)$-dynomiale minimale de $(E_{\mathbb{K}})$ est irréductible.
	
\end{theorem}

\begin{proof}

Comme $\mathbb{K}$ est un corps, $\mathbb{K}$ est intègre. Donc, le polynôme $X(b-X)$ a pour seules racines $0_{\mathbb{K}}$ et $b$ sur $\mathbb{K}$. Le polynôme $P(X)=X(b-a^{-1}bX)-(1_{\mathbb{K}}-a^{-1}b)$ a pour discriminant :
\[\Delta_{2}=b^{2}-4_{\mathbb{K}}(-a^{-1}b(-(1_{\mathbb{K}}-a^{-1}b)))=b^{2}+4_{\mathbb{K}}a^{-1}b(a^{-1}b-1_{\mathbb{K}}).\]

\noindent Par hypothèse, $\Delta_{2}$ n'est pas un carré sur $\mathbb{K}$. Donc, $P$ n'a pas de racine sur $\mathbb{K}$.
\\
\\Le polynôme $Q(X)=X(a-ab^{-1}X)-(1_{\mathbb{K}}-ab^{-1})$ a pour discriminant :
\[\Delta_{1}=a^{2}-4_{\mathbb{K}}(-ab^{-1}(-(1_{\mathbb{K}}-ab^{-1})))=a^{2}+4_{\mathbb{K}}ab^{-1}(ab^{-1}-1_{\mathbb{K}}).\]

\noindent Par hypothèse, $\Delta_{1}$ n'est pas un carré sur $\mathbb{K}$. Donc, $Q$ n'a pas de racine sur $\mathbb{K}$. 
\\
\\Par le lemme précédent, la solution $(a,b)$-dynomiale minimale de $(E_{\mathbb{K}})$ est irréductible.

\end{proof}

\begin{remark}

{\rm Soient $\mathbb{K}=\mathbb{Z}/p\mathbb{Z}$ avec $p$ premier différent de 2 et $\overline{k} \in \mathbb{K}^{*}$. Si $a=\overline{k}$ et $b=-\overline{k}$ alors $a^{2}+4_{\mathbb{K}}ab^{-1}(ab^{-1}-1_{\mathbb{K}})=b^{2}+4_{\mathbb{K}}a^{-1}b(a^{-1}b-1_{\mathbb{K}})=\overline{k}^{2}+\overline{8}$. Le théorème précédent généralise donc bien le théorème \ref{32}.}
	
\end{remark}		
		
\begin{examples}
{\rm Voici quelques exemples d'utilisation du théorème \ref{35}.
\begin{itemize}
\item Soient $\mathbb{K}=\mathbb{Z}/23\mathbb{Z}$, $\overline{a}=\overline{2}$ et $\overline{b}=\overline{5}$. On a $\overline{a}^{-1}=\overline{12}$ et $\overline{b}^{-1}=\overline{14}$. Calculons $\Delta_{1}$ et $\Delta_{2}$. On a $\Delta_{1}=\overline{2}^{2}+\overline{4} \times \overline{2} \times \overline{14}(\overline{2} \times \overline{14}-\overline{1})=\overline{15}$ et $\Delta_{2}=\overline{5}^{2}+\overline{4} \times \overline{12}\times \overline{5}(\overline{12}\times \overline{5}-\overline{1})=\overline{17}$. Or, les carrés modulo 23 sont : $\{\overline{0},\overline{1},\overline{4},\overline{9},\overline{16},\overline{2},\overline{13},\overline{3},\overline{18},\overline{12},\overline{8},\overline{6}\}$. Donc, par le théorème \ref{35}, la solution $(\overline{2},\overline{5})$-dynomiale minimale de $(E_{\mathbb{K}})$ est irréductible.
\\
\item Soient $\mathbb{K}=\mathbb{F}_{25}=\frac{(\mathbb{Z}/5\mathbb{Z})[X]}{\left\langle X^{2}+X+\overline{1}\right\rangle}$, $a=\overline{2}$ et $b=X$. On a $a^{-1}=\overline{3}$ et $b^{-1}=-X-\overline{1}$. Par ailleurs :
\[\Delta_{1}=\overline{2}^{2}+\overline{4} \times \overline{2} \times (-X-\overline{1})(\overline{2} \times (-X-\overline{1})-\overline{1})=-X+\overline{2};\] 
\[\Delta_{2}=X^{2}+\overline{4} \times \overline{3}\times X(\overline{3}\times X-\overline{1})=X+\overline{3}.\]
\noindent $\Delta_{1}$ et $\Delta_{2}$ ne sont pas des carrés dans $\mathbb{F}_{25}$ (voir annexe \ref{A}). Donc, par le théorème \ref{35}, la solution $(\overline{2},X)$-dynomiale minimale de $(E_{\mathbb{K}})$ est irréductible.

\end{itemize}
}
\end{examples} 		
		
\begin{remarks}

{\rm i) Pour certaines valeurs de $a$ et $b$, $\Delta_{1}$ et $\Delta_{2}$ peuvent prendre des formes plus simples. Par exemple, si $\mathbb{K}$ est un corps fini de caractéristique différente de 2, $a \in \mathbb{K}^{*}$ et $b=2_{\mathbb{K}}a$, on a $\Delta_{1}=a^{2}-1_{\mathbb{K}}$ et $\Delta_{2}=4_{\mathbb{K}}(a^{2}+2_{\mathbb{K}})$. Comme $4_{\mathbb{K}}$ est un carré, l'hypothèse du théorème \ref{35} se réduit à montrer que $a^{2}-1_{\mathbb{K}}$ et $a^{2}+2_{\mathbb{K}}$ ne sont pas des carrés.
 \\
\\ii) Les conditions du lemme \ref{34} (et donc du théorème \ref{35}) sont suffisantes mais pas nécessaires. En effet, comme on l'avait déjà observé après la preuve du théorème \ref{32}, on peut avoir des racines pour les polynômes présents dans l'énoncé et avoir tout de même une solution irréductible. Par exemple, prenons $N=59$. La solution $(\overline{2},-\overline{2})$-dynomiale minimale de $(E_{\mathbb{Z}/59\mathbb{Z}})$ est irréductible et le polynôme $X(X-\overline{2})-\overline{2}$ a pour racines $\overline{12}$ et $\overline{49}$.
}
	
\end{remarks}	

Une autre possibilité est de chercher, pour un couple $(a,b)$ donné, des corps sur lesquels la solution $(a,b)$-dynomiale minimale est irréductible. Pour cela, la loi de réciprocité quadratique de Gauss va nous être très utile. En suivant cette direction, on peut obtenir, par exemple, le résultat ci-dessous :
		
\begin{proposition}
\label{36}

Soit $p$ un nombre premier. Si $p \equiv \pm 7 [24]$ alors la solution $(\overline{2},\overline{4})$-dynomiale minimale de $(E_{\mathbb{Z}/p\mathbb{Z}})$ est irréductible.

\end{proposition}

\begin{proof}

Soit $p$ un nombre premier vérifiant $p \equiv \pm 7 [24]$. Commençons par remarquer que $\overline{2}$ et $\overline{4}$ sont des éléments inversibles de $\mathbb{Z}/p\mathbb{Z}$ différents de $\{\pm \overline{1}\}$.
\\
\\Soient $\overline{a}=\overline{2}$, $\overline{b}=\overline{4}$, $\Delta_{1}=\overline{a}^{2}-\overline{1}=\overline{3}$ et $\Delta_{2}=\overline{4}(\overline{a}^{2}+\overline{2})=\overline{24}$. Comme $p \equiv \pm 7 [24]$, $p \equiv \pm 5 [12]$. Donc, par la proposition \ref{18}, $\overline{3}$ n'est pas un carré modulo $p$. De plus, puisque $p \equiv \pm 7 [24]$ on a $p \equiv \pm 1 [8]$ et donc $\overline{2}$ est un carré modulo $p$ (Théorème \ref{17}). Ainsi, 
\[\left(\dfrac{24}{p}\right)=\left(\dfrac{2}{p}\right)^{3}\left(\dfrac{3}{p}\right)=1^{3} \times (-1)=-1.\]
\noindent Comme $\mathbb{Z}/p\mathbb{Z}$ est un corps, la solution $(\overline{2},\overline{4})$-dynomiale minimale de $(E_{\mathbb{Z}/p\mathbb{Z}})$ est irréductible, par le théorème \ref{35}.

\end{proof}

\begin{example}
{\rm La solution $(\overline{2},\overline{4})$-dynomiale minimale de $(E_{\mathbb{Z}/17\mathbb{Z}})$ est irréductible.
}
\end{example} 

\subsubsection{Irréductibilité sur les anneaux $\mathbb{Z}/N\mathbb{Z}$} 

On a obtenu dans la sous-section précédente des résultats intéressants lorsque l'on considère des corps finis. Cela dit, on peut également chercher des conditions suffisantes pour l'irréductibilité sur certains anneaux qui ne sont pas des corps. Pour cela, on donne déjà un résultat préliminaire.

\begin{lemma}
\label{37}

Soient $p$ un nombre premier, $n \in \mathbb{N}^{*}$ et $N=p^{n}$. Soit $k \in \mathbb{Z}$ premier avec $p$. Le polynôme $X(X-\overline{k})$ a pour seule racines $\overline{0}$ et $\overline{k}$.

\end{lemma}

\begin{proof}

On vérifie aisément que $\overline{0}$ et $\overline{k}$ sont des racines de $X(X-\overline{k})$. Soit $\overline{x} \in \mathbb{Z}/p^{n}\mathbb{Z}$ tel que $\overline{x}(\overline{x}-\overline{k})=\overline{0}$. Si $\overline{x}$ est inversible alors on a $\overline{x}=\overline{k}$. Supposons maintenant que $\overline{x}$ n'est pas inversible. Il existe $1 \leq j \leq n$ et $u$ un entier premier avec $p$ tels que $x=up^{j}$. 
\\
\\Supposons par l'absurde que $j<n$. Comme $p^{n}$ divise $up^{j}(up^{j}-k)$, $p$ divise $u(up^{j}-k)$. Par le lemme de Gauss, $p$ divise $up^{j}-k$ (puisque $u$ et $p$ sont premiers entre eux). Ainsi, $p$ divise $k$, ce qui est absurde. Par conséquent, $j=n$ et $\overline{x}=\overline{0}$.

\end{proof}
	
\begin{proposition}
\label{38}

Soient $p$ un nombre premier supérieur à 5 et $n \geq 1$. Soient $a$ et $b$ deux entiers non divisibles par $p$ tels que $a+p\mathbb{Z} \neq b+p\mathbb{Z}$. Si la solution $(a+p\mathbb{Z},b+p\mathbb{Z})$-dynomiale minimale de $(E_{\mathbb{Z}/p\mathbb{Z}})$ est irréductible alors la solution $(a+p^{n}\mathbb{Z},b+p^{n}\mathbb{Z})$-dynomiale minimale de $(E_{\mathbb{Z}/p^{n}\mathbb{Z}})$ est irréductible.

\end{proposition}		

\begin{proof}

Soient $a$ et $b$ deux entiers non divisibles par $p$ tels que $a+p\mathbb{Z} \neq b+p\mathbb{Z}$. On raisonne par contraposée. On suppose donc qu'il existe un $n \geq 1$ tel que la solution $(a+p^{n}\mathbb{Z},b+p^{n}\mathbb{Z})$-dynomiale minimale de $(E_{\mathbb{Z}/p^{n}\mathbb{Z}})$ est réductible. On va montrer que la solution $(a+p\mathbb{Z},b+p\mathbb{Z})$-dynomiale minimale de $(E_{\mathbb{Z}/p\mathbb{Z}})$ est réductible. Dans cette preuve, on pose, pour $u \in \mathbb{Z}$, $\overline{u}:=u+p^{n}\mathbb{Z}$.
\\
\\Soient $m$ la taille de la solution $(a+p\mathbb{Z},b+p\mathbb{Z})$-dynomiale minimale de $(E_{\mathbb{Z}/p\mathbb{Z}})$ et $h$ la taille de la solution $(\overline{a},\overline{b})$-dynomiale minimale de $(E_{\mathbb{Z}/p^{n}\mathbb{Z}})$. Par définition, les entiers $m$ et $h$ sont pairs. De plus, il existe $\epsilon \in \{-1,1\}$ tel que $M_{m}(a+p\mathbb{Z},b+p\mathbb{Z},\ldots,a+p\mathbb{Z},b+p\mathbb{Z})=(\epsilon+p\mathbb{Z})Id$.
\\
\\Il existe deux solutions de $(E_{\mathbb{Z}/p^{n}\mathbb{Z}})$, $(\overline{c_{1}},\ldots,\overline{c_{l'}})$ et $(\overline{d_{1}},\ldots,\overline{d_{l}})$ telles que $3 \leq l,l' \leq h-1$ et :
\[(\overline{a},\overline{b},\ldots,\overline{a},\overline{b}) \sim (\overline{c_{1}},\ldots,\overline{c_{l'}}) \oplus (\overline{d_{1}},\ldots,\overline{d_{l}}).\]

\noindent Sans perte de généralité, on peut supposer que cette relation est en fait une égalité. Pour obtenir la réductibilité souhaitée, on procède par étapes.
\\
\\i) Montrons pour commencer que $l$ et $l'$ sont impairs. 
\\
\\Supposons par l'absurde que $l$ est pair. On a $(\overline{d_{1}},\ldots,\overline{d_{l}})=(\overline{x},\overline{a},\overline{b},\ldots,\overline{a},\overline{b},\overline{y})$. Par le lemme \ref{33}, on a $\overline{y}=\overline{a}(\overline{b})^{-1}\overline{x}$ et $\overline{0}=\overline{x}(\overline{b}-\overline{x})$. Par le lemme \ref{37}, on obtient $\overline{x}=\overline{y}=\overline{0}$ ou $\overline{x}=\overline{b}$ et $\overline{y}=\overline{a}$. 
\\
\\Si $\overline{x}=\overline{y}=\overline{0}$ alors $(\overline{d_{1}},\ldots,\overline{d_{l}}) \sim (\overline{a},\overline{b},\ldots,\overline{a},\overline{b}) \oplus (\overline{0},\overline{0},\overline{0},\overline{0})$. En particulier, le $(l-2)$-uplet $(\overline{a},\overline{b},\ldots,\overline{a},\overline{b})$ est une solution, ce qui contredit la minimalité de $h$. 
\\
\\Si $\overline{x}=\overline{b}$ et $\overline{y}=\overline{a}$ alors on obtient une solution $(\overline{a},\overline{b})$ dynomiale de taille $l<h$, ce qui est absurde.
\\
\\Ainsi, $l$ est impair. De plus, comme $h$ est pair, $l'=h+2-l$ est également impair. En particulier, on a $(\overline{d_{1}},\ldots,\overline{d_{l}})=(\overline{x},\overline{b},\overline{a},\overline{b},\ldots,\overline{a},\overline{b},\overline{y})$.
\\
\\ii) Montrons que $l$ peut s'écrire sous la forme $qm+r$ avec $q$ entier et $3 \leq r \leq m$.
\\
\\On effectue la division euclidienne de $l$ par $m$. Il existe $(q,r) \in \mathbb{N}^{2}$ tel que $l=qm+r$ et $r<m$. Comme $l$ est impair et $m$ est pair, $r$ est impair. Donc, $r \geq 1$. Supposons par l'absurde que $r=1$.
\\
\\Comme $p$ divise $p^{n}$, $(d_{1}+p\mathbb{Z},\ldots,d_{l}+p\mathbb{Z})$ est une $\lambda$-quiddité. Donc, $(d_{l}+p\mathbb{Z},d_{1}+p\mathbb{Z},\ldots,d_{l-1}+p\mathbb{Z})$ est également une $\lambda$-quiddité. Ainsi, il existe $\alpha \in \{-1,1\}$ tel que :
\begin{eqnarray*}
(\alpha+p\mathbb{Z}) Id &=& M_{l}(d_{l}+p\mathbb{Z},d_{1}+p\mathbb{Z},\ldots,d_{l-1}+p\mathbb{Z}) \\
                        &=& M_{m}(a+p\mathbb{Z},b+p\mathbb{Z},\ldots,a+p\mathbb{Z},b+p\mathbb{Z})^{q-1} \\
											  &\times& M_{m+1}(d_{l}+p\mathbb{Z},d_{1}+p\mathbb{Z},b+p\mathbb{Z},a+p\mathbb{Z},\ldots,b+p\mathbb{Z},a+p\mathbb{Z},b+p\mathbb{Z})\\
											  &=& (\epsilon+p\mathbb{Z})^{q-1} M_{m+1}(d_{l}+p\mathbb{Z},d_{1}+p\mathbb{Z},b+p\mathbb{Z},a+p\mathbb{Z},\ldots,b+p\mathbb{Z},a+p\mathbb{Z},b+p\mathbb{Z}).
\end{eqnarray*}
\noindent En particulier, $M_{m+1}(d_{1}+p\mathbb{Z},b+p\mathbb{Z},a+p\mathbb{Z},\ldots,b+p\mathbb{Z},a+p\mathbb{Z},b+p\mathbb{Z},d_{l}+p\mathbb{Z})= \pm Id$. Ainsi, par la proposition \ref{12}, $K_{m-1}(b+p\mathbb{Z},a+p\mathbb{Z},\ldots,b+p\mathbb{Z},a+p\mathbb{Z},b+p\mathbb{Z})=\pm 1+p\mathbb{Z}$. Or, par la proposition \ref{12}, on sait que $K_{m-1}(b+p\mathbb{Z},a+p\mathbb{Z},\ldots,b+p\mathbb{Z},a+p\mathbb{Z},b+p\mathbb{Z})=0+p\mathbb{Z}$, puisqu'on dispose de l'égalité matricielle $M_{m}(a+p\mathbb{Z},b+p\mathbb{Z},\ldots,a+p\mathbb{Z},b+p\mathbb{Z})=(\epsilon+p\mathbb{Z})Id$. Comme $p \neq 1$, ceci est absurde. Ainsi, $r \neq 1$ et donc $r \geq 3$.
\\
\\iii) On conclut en construisant une solution permettant de réduire la solution $(a+p\mathbb{Z},b+p\mathbb{Z})$-dynomiale minimale de $(E_{\mathbb{Z}/p\mathbb{Z}})$.
\\
\\On dispose des égalités suivantes :
\begin{eqnarray*}
(\alpha+p\mathbb{Z}) Id &=& M_{l}(d_{l}+p\mathbb{Z},d_{1}+p\mathbb{Z},\ldots,d_{l-1}+p\mathbb{Z}) \\
                        &=& M_{m}(a+p\mathbb{Z},b+p\mathbb{Z},\ldots,a+p\mathbb{Z},b+p\mathbb{Z})^{q} \\
												&\times& M_{r}(d_{l}+p\mathbb{Z},d_{1}+p\mathbb{Z},b+p\mathbb{Z},a+p\mathbb{Z},\ldots,b+p\mathbb{Z},a+p\mathbb{Z},b+p\mathbb{Z}) \\
											  &=& (\epsilon+p\mathbb{Z})^{q} M_{r}(d_{l}+p\mathbb{Z},d_{1}+p\mathbb{Z},b+p\mathbb{Z},a+p\mathbb{Z},\ldots,b+p\mathbb{Z},a+p\mathbb{Z},b+p\mathbb{Z}).
\end{eqnarray*}
\noindent Ainsi, $M_{r}(d_{1}+p\mathbb{Z},b+p\mathbb{Z},a+p\mathbb{Z},\ldots,b+p\mathbb{Z},a+p\mathbb{Z},b+p\mathbb{Z},d_{l}+p\mathbb{Z})=(\alpha \times \epsilon^{q}+p\mathbb{Z}) Id$. Comme $3 \leq r \leq m-1$, on peut utiliser cette solution pour réduire la solution $(a+p\mathbb{Z},b+p\mathbb{Z})$-dynomiale minimale. Cette dernière est donc réductible.

\end{proof}

\begin{remarks}
{\rm i) Comme $p$ divise $p^{n}$, $M_{h}(a+p\mathbb{Z},b+p\mathbb{Z},\ldots,a+p\mathbb{Z},b+p\mathbb{Z})= \pm Id$. La taille de la solution $(a+p\mathbb{Z},b+p\mathbb{Z})$-dynomiale minimale étant égale au double de l'ordre de $M_{2}(a+p\mathbb{Z},b+p\mathbb{Z})$, $h$ est un multiple de $m$. En fait, on peut même montrer qu'il existe $i \geq 0$ tel que $h=mp^{i}$. En effet, il existe $B \in SL_{2}(\mathbb{Z})$ tel que $M_{2}(a,b)^{\frac{m}{2}}=\epsilon Id+pB$. Ainsi, en utilisant et le binôme de Newton, on constate qu'il existe $C \in SL_{2}(\mathbb{Z})$ tel que $M_{2}(a,b)^{\frac{pm}{2}}=\epsilon^{p} Id+p^{2}C$, ce qui implique que la taille de la solution $(a+p^{2}\mathbb{Z},b+p^{2}\mathbb{Z})$-dynomiale minimale est égale à $2m'$, avec $m'$ multiple de $\frac{m}{2}$ et diviseur de $\frac{pm}{2}$. Elle est donc égale à $m$ ou $pm$. Par récurrence, on obtient ensuite le résultat souhaité.
\\
\\ii) La réciproque de la proposition \ref{38} est fausse. 
Prenons par exemple, $N=49$, $a=2$ et $b=3$. Informatiquement, on trouve que la solution $(\overline{2},\overline{3})$-dynomiale minimale de $(E_{\mathbb{Z}/49\mathbb{Z}})$ est irréductible de taille 56 alors que la solution $(2+7\mathbb{Z},3+7\mathbb{Z})$-dynomiale minimale de $(E_{\mathbb{Z}/7\mathbb{Z}})$ est réductible de taille 8. On peut la réduire avec la solution $(5+7\mathbb{Z},2+7\mathbb{Z},3+7\mathbb{Z},2+7\mathbb{Z},5+7\mathbb{Z})$.
}
\end{remarks}

\begin{examples}
{\rm En combinant la proposition précédente aux cas déjà étudiés, on peut obtenir les résultats ci-dessous.
\begin{itemize}
\item La solution $(\overline{2},\overline{4})$-dynomiale minimale de $(E_{\mathbb{Z}/49\mathbb{Z}})$ est irréductible.
\item La solution $(\overline{2},\overline{5})$-dynomiale minimale de $(E_{\mathbb{Z}/529\mathbb{Z}})$ est irréductible, puisque $529=23^{2}$.
\end{itemize}
}
\end{examples}

On aimerait disposer d'un résultat analogue pour tous les anneaux $\mathbb{Z}/N\mathbb{Z}$. Malheureusement, il est possible, si $N$ est un entier qui n'est pas une puissance d'un nombre premier, d'avoir une solution dynomiale minimale qui soit irréductible modulo chaque facteur premier de $N$ mais réductible modulo $N$. Prenons par exemple, $N=77$, $a=2$ et $b=5$. Comme $7 \not\equiv \pm 1 [12]$, la solution $(2+7\mathbb{Z},5+7\mathbb{Z})$-dynomiale minimale de $(E_{\mathbb{Z}/7\mathbb{Z}})$ est irréductible (Théorème \ref{32} et proposition \ref{18}). Plaçons-nous maintenant modulo 11 et reprenons les notations du théorème \ref{35}. On a $\Delta_{1}=\Delta_{2}=7+11\mathbb{Z}$ et $7+11\mathbb{Z}$ n'est pas un carré de $\mathbb{Z}/11\mathbb{Z}$. Donc, par le théorème \ref{35}, la solution $(2+11\mathbb{Z},5+11\mathbb{Z})$-dynomiale minimale de $(E_{\mathbb{Z}/11\mathbb{Z}})$ est irréductible. En revanche, la solution $(\overline{2},\overline{5})$-dynomiale minimale de $(E_{\mathbb{Z}/77\mathbb{Z}})$ est réductible de taille 30. On peut par exemple la réduire avec la solution $(\overline{33},\overline{2},\overline{5},\overline{2},\overline{5},\overline{2},\overline{5},\overline{2},\overline{5},\overline{2},\overline{5},\overline{44})$.
\\
\\\indent En réalité, ce qui pose problème dans le cas général, c'est l'impossibilité d'avoir un résultat analogue au lemme \ref{37}. En effet, posons $N:=uv$ avec $u,v>1$ premiers entre eux, $k$ un entier premier avec $N$ et $x$ l'unique entier compris entre 0 et $N-1$ vérifiant $x \equiv 0 [u]$ et $x \equiv k [v]$. Par le lemme chinois, $\overline{x}(\overline{x}-\overline{k})=\overline{0}$ et $\overline{x} \neq \overline{0}, \overline{k}$. Cette existence de racines non triviales empêche d'écarter la possibilité qu'il existe des $\lambda$-quiddités de taille paire permettant de réduire la solution dynomiale minimale considérée, ce qui réduit à néant tout espoir d'obtenir un résultat semblable à celui de la proposition \ref{38}. En effet, reprenons les notations de la preuve. La démonstration du résultat repose sur la réduction, obtenue par suppression de blocs de taille $m$, des $\lambda$-quiddités de taille $l$ permettant de réduire la solution modulo $p^{n}$. On aboutit ainsi, après réduction, à une $\lambda$-quiddité dont la taille permet de réduire la solution dynomiale minimale modulo $p$. Cependant, si $l$ est de la forme $dm+2$ ou $d'm$ alors on obtient après réduction une $\lambda$-quiddité qui ne permet pas de réduire la solution dynomiale minimale modulo $p$.
		
\section{Solutions trinomiales}
\label{tri}

On va s'intéresser dans cette section à une nouvelle famille de $\lambda$-quiddités, définies par la répétition d'un motif de taille 3, qui nous permettront d'avoir des minorations intéressantes pour $\ell_{\mathbb{F}_{q}}$.

\subsection{Définition et caractérisation de l'irréductibilité}

\noindent L'objectif de cette partie est d'étudier la nouvelle famille de solutions définie ci-dessous :

\begin{definition}
\label{41}

Soient $A$ un anneau commutatif unitaire fini et $u \in U(A)$. On appelle solution $u$-trinomiale minimale une solution de \eqref{p} de la forme $(u,u^{-1},u^{-1},\ldots,u,u^{-1},u^{-1})$ et qui est de taille minimale. On appelle solution trinomiale minimale sur $A$ une solution $u$-trinomiale minimale pour un certain $u \in U(A)$.

\end{definition} 

Comme précédemment, la solution $u$-trinomiale minimale de \eqref{p} existe toujours, puisque la matrice $M_{3}(u,u^{-1},u^{-1})$ est d'ordre fini dans $SL_{2}(A)/\{-Id,Id\}$. De plus, la solution $u$-trinomiale minimale de \eqref{p} a la même taille que la solution $-u$-trinomiale minimale de \eqref{p}.

\begin{examples}
{\rm On donne ci-dessous quelques exemples de solutions trinomiales minimales.
\begin{itemize}
\item $(1_{A},1_{A},1_{A})$ est une solution trinomiale minimale de \eqref{p}.
\item $A=\mathbb{Z}/7\mathbb{Z}$. $(\overline{4},\overline{2},\overline{2},\overline{4},\overline{2},\overline{2},\overline{4},\overline{2},\overline{2})$ est la solution $\overline{4}$-trinomiale minimale de \eqref{p}.
\item $(X+\overline{1},X,X,X+\overline{1},X,X,X+\overline{1},X,X)$ est la solution $(X+\overline{1})$-trinomiale minimale sur le corps $\mathbb{F}_{4}=\frac{\mathbb{F}_{2}[X]}{\left\langle X^{2}+X+\overline{1}\right\rangle}$.
\end{itemize}
}
\end{examples}

\noindent Pour étudier ces solutions, on va démontrer le résultat important énoncé ci-dessous :

\begin{theorem}
\label{42}

Soient $\mathbb{K}$ un corps fini et $u \in U(\mathbb{K})$. Soient $m$ ta taille de la solution $u$-trinomiale minimale et $o(u)$ l'ordre de $u$ dans $\mathbb{K}^{*}$. On a trois cas :
\begin{itemize}
\item si ${\rm car}(\mathbb{K})=2$ alors $m=3o(u)$;
\item si ${\rm car}(\mathbb{K})\neq 2$ et $o(u)$ est pair alors $m=\frac{3o(u)}{2}$;
\item si ${\rm car}(\mathbb{K})\neq 2$ et $o(u)$ est impair alors $m=3o(u)$.
\end{itemize}

\noindent De plus, si $u \neq \pm 1_{\mathbb{K}}$, la solution $u$-trinomiale minimale de $(E_{\mathbb{K}})$ est irréductible si et seulement si $u$ n'est pas racine des polynômes $X^{2l}+X^{l+1}-1_{\mathbb{K}}$ et $X^{2l}-X^{l+1}-1_{\mathbb{K}}$ pour $1 \leq l \leq E\left[\frac{m}{6}\right]$.

\end{theorem}

\begin{proof}

Le $3n$-uplet $(u,u^{-1},u^{-1},\ldots,u,u^{-1},u^{-1})$ est une solution de $(E_{\mathbb{K}})$ si et seulement si le $3n$-uplet $(u^{-1},u,u^{-1}\ldots,u^{-1},u,u^{-1})$ l'est également. Or, 
\begin{eqnarray*}
M_{3n}(u^{-1},u,u^{-1},\ldots,u^{-1},u,u^{-1}) &=& M_{3}(u^{-1},u,u^{-1})^{n} \\
                                               &=& \left(\begin{pmatrix}
    -u^{-1} &  0_{\mathbb{K}} \\
     0_{\mathbb{K}}  & -u
\end{pmatrix}\right)^{n} \\
                                               &=& \begin{pmatrix}
    (-u^{-1})^{n} &  0_{\mathbb{K}} \\
     0_{\mathbb{K}}  & (-u)^{n}
\end{pmatrix}.
\end{eqnarray*}

\noindent Donc, $m=3v$ où $v$ est le plus petit entier tel que $u^{v}=\pm 1_{\mathbb{K}}$. Si ${\rm car}(\mathbb{K})=2$ alors nécessairement $v=o(u)$ et $m=3o(u)$. Supposons maintenant ${\rm car}(\mathbb{K})\neq 2$. Si $o(u)$ est impair alors $v=o(u)$. En effet, supposons par l'absurde qu'il existe $k<o(u)$ tel que $u^{k}=\pm 1_{\mathbb{K}}$. Comme $k<o(u)$, on a, par définition de l'ordre, $u^{k}=-1_{\mathbb{K}}$ et donc $u^{2k}=1_{\mathbb{K}}$. Ainsi, $2k$ est un multiple de $o(u)$, c'est-à-dire $k$ est un multiple de $o(u)$ (puisque $o(u)$ est impair). Cela implique $k \geq o(u)$, ce qui est absurde. Si $o(u)$ est pair alors $v=\frac{o(u)}{2}$. En effet, puisque $u^{o(u)}=1_{\mathbb{K}}$, on a $(u^{\frac{o(u)}{2}}-1_{\mathbb{K}})(u^{\frac{o(u)}{2}}+1_{\mathbb{K}})=0_{\mathbb{K}}$. Or, $\mathbb{K}$ est intègre et $u^{\frac{o(u)}{2}} \neq 1_{\mathbb{K}}$ (puisque $\frac{o(u)}{2} <o(u)$). Donc, $u^{\frac{o(u)}{2}}=-1_{\mathbb{K}}$. De plus, s'il existait $k<\frac{o(u)}{2}$ tel que $u^{k}=\pm 1_{\mathbb{K}}$ alors on aurait $u^{2k}=1_{\mathbb{K}}$ et $2k<o(u)$, ce qui est impossible.
\\
\\Ainsi, $m$ vérifie bien les égalités annoncées.
\\
\\On va maintenant chercher des conditions nécessaires et suffisantes pour que cette solution soit irréductible. Pour cela, on suppose $u \neq \pm 1_{\mathbb{K}}$, ce qui implique $m \geq 6$, et on note $(a_{1},\ldots,a_{m})$ la solution $(u,u^{-1},u^{-1},\ldots,u,u^{-1},u^{-1})$.
\\
\\Si $(a_{1},\ldots,a_{m})$ est réductible alors elle est équivalente à la somme d'une solution de taille $j+2$ avec une solution de taille $g+2$, avec $j,g \geq 1$. Comme on a $m+2=j+g+4$ ($m \geq 4$ car $u \neq \pm 1_{\mathbb{K}}$), on a nécessairement $j \leq \frac{m-2}{2}$ ou $g \leq \frac{m-2}{2}$. Donc, il existe $1 \leq j \leq \frac{m-2}{2}$ et $i$ tels que $K_{j}(a_{i},\ldots,a_{j+i-1})=\pm 1_{\mathbb{K}}$ (avec les indices des $a_{k}$ vus modulo $m$). Réciproquement, s'il existe $1 \leq j \leq \frac{m-2}{2}$ et $i$ tels que $K_{j}(a_{i},\ldots,a_{j+i-1})=\pm 1_{\mathbb{K}}$ (avec les indices des $a_{k}$ vus modulo $m$) alors on peut construire grâce au lemme \ref{131} une solution de taille $j+2$ qui nous permettra de réduire $(a_{1},\ldots,a_{m})$. Autrement dit, l'étude de la réductibilité dépend du calcul de tous les sous-continuants de $(a_{1},\ldots,a_{m})$ possibles, c'est-à-dire des neufs calculs suivants pour $l \geq 1$ :
\begin{enumerate}
\item $K_{3l}(u,u^{-1},u^{-1},\ldots,u,u^{-1},u^{-1})$;
\item $K_{3l}(u^{-1},u^{-1},u,\ldots,u^{-1},u^{-1},u)$;
\item $K_{3l}(u^{-1},u,u^{-1},\ldots,u^{-1},u,u^{-1})$;
\item $K_{3l-1}(u,u^{-1},u^{-1},\ldots,u,u^{-1},u^{-1},u,u^{-1})$;
\item $K_{3l-1}(u^{-1},u,u^{-1},\ldots,u^{-1},u,u^{-1},u^{-1},u)$;
\item $K_{3l-1}(u^{-1},u^{-1},u,\ldots,u^{-1},u^{-1},u,u^{-1},u^{-1})$;
\item $K_{3l-2}(u^{-1},u,u^{-1},\ldots,u^{-1},u,u^{-1},u^{-1})$;
\item $K_{3l-2}(u^{-1},u^{-1},u,\ldots,u^{-1},u^{-1},u,u^{-1})$;
\item $K_{3l-2}(u,u^{-1},u^{-1},\ldots,u,u^{-1},u^{-1},u)$.

\end{enumerate}

\noindent Avec la proposition \ref{13} i), on constate que 2 est égal à 1, 4 est égal à 5 et 7 est égal à 8. Il ne reste donc que six continuants à calculer. Pour cela, on va utiliser des calculs matriciels. On a, par ce qui précède, pour tout $l \geq 1$ :
\begin{eqnarray*}
M &=& M_{3l}(u^{-1},u,u^{-1},\ldots,u^{-1},u,u^{-1}) \\
  &=& \begin{pmatrix}
   K_{3l}(u^{-1},u,u^{-1},\ldots,u^{-1},u,u^{-1}) &  -K_{3l-1}(u,u^{-1},\ldots,u^{-1},u,u^{-1}) \\
   K_{3l-1}(u^{-1},u,u^{-1},\ldots,u^{-1},u)  & -K_{3l-2}(u,u^{-1},\ldots,u^{-1},u)
		\end{pmatrix} \\
	&=& \begin{pmatrix}
    (-u^{-1})^{l} &  0_{\mathbb{K}} \\
     0_{\mathbb{K}}  & (-u)^{l}
		\end{pmatrix}.
\end{eqnarray*}
\noindent On va maintenant montrer par récurrence la formule ci-dessous :
\[M_{3l}(u,u^{-1},u^{-1},\ldots,u,u^{-1},u^{-1})=\begin{pmatrix}
   (-u)^{l} & (-1_{\mathbb{K}})^{l}((u^{-1})^{l+1}-u^{l-1})  \\
   0_{\mathbb{K}}  & (-u^{-1})^{l}
		\end{pmatrix}.\]

\noindent On vérifie aisément, avec la proposition \ref{12}, que la formule est vraie pour $l=1$. Supposons qu'il existe un $l \geq 1$ tel que la formule est vraie au rang $l$. On a :
\begin{eqnarray*}
M &=& M_{3(l+1)}(u,u^{-1},u^{-1},\ldots,u,u^{-1},u^{-1}) \\
  &=& M_{3l}(u,u^{-1},u^{-1},\ldots,u,u^{-1},u^{-1})M_{3}(u,u^{-1},u^{-1}) \\
	&=& \begin{pmatrix}
   (-u)^{l} & (-1_{\mathbb{K}})^{l}((u^{-1})^{l+1}-u^{l-1})  \\
   0_{\mathbb{K}}  & (-u^{-1})^{l}
		\end{pmatrix}\begin{pmatrix}
   -u & 1_{\mathbb{K}}-(u^{-1})^{2}  \\
   0_{\mathbb{K}}  & -u^{-1}
		\end{pmatrix} \\
	&=& \begin{pmatrix}
   (-u)^{l+1} & (-u)^{l}+(-1_{\mathbb{K}})^{l+1}u^{l-2}-u^{-1}(-1_{\mathbb{K}})^{l}((u^{-1})^{l+1}-u^{l-1})  \\
   0_{\mathbb{K}}  & (-u^{-1})^{l+1}
		\end{pmatrix} \\
	&=& \begin{pmatrix}
   (-u)^{l+1} & (-1_{\mathbb{K}})^{l+1}(-u^{l}+u^{l-2}+(u^{-1})^{l+2}-u^{l-2}) \\
   0_{\mathbb{K}}  & (-u^{-1})^{l+1}
		\end{pmatrix} \\
	&=& \begin{pmatrix}
   (-u)^{l+1} & (-1_{\mathbb{K}})^{l+1}((u^{-1})^{(l+1)+1}-u^{(l+1)-1}) \\
   0_{\mathbb{K}}  & (-u^{-1})^{l+1}
		\end{pmatrix}.
\end{eqnarray*}

\noindent Ainsi, la formule est vraie par récurrence.
\\
\\En combinant ces calculs matriciels avec la proposition \ref{12}, on obtient les égalités suivantes :
\begin{enumerate}
\item $K_{3l}(u,u^{-1},u^{-1},\ldots,u,u^{-1},u^{-1})=(-u)^{l}$;
\item $K_{3l}(u^{-1},u^{-1},u,\ldots,u^{-1},u^{-1},u)=(-u)^{l}$;
\item $K_{3l}(u^{-1},u,u^{-1},\ldots,u^{-1},u,u^{-1})=(-u^{-1})^{l}$;
\item $K_{3l-1}(u,u^{-1},u^{-1},\ldots,u,u^{-1},u^{-1},u,u^{-1})=0_{\mathbb{K}}$;
\item $K_{3l-1}(u^{-1},u,u^{-1},\ldots,u^{-1},u,u^{-1},u^{-1},u)=0_{\mathbb{K}}$;
\item $K_{3l-1}(u^{-1},u^{-1},u,\ldots,u^{-1},u^{-1},u,u^{-1},u^{-1})=(-1_{\mathbb{K}})^{l+1}((u^{-1})^{l+1}-u^{l-1})$;
\item $K_{3l-2}(u^{-1},u,u^{-1},\ldots,u^{-1},u,u^{-1},u^{-1})=(-1)^{l+1}(u^{-1})^{l}$;
\item $K_{3l-2}(u^{-1},u^{-1},u,\ldots,u^{-1},u^{-1},u,u^{-1})=(-1)^{l+1}(u^{-1})^{l}$;
\item $K_{3l-2}(u,u^{-1},u^{-1},\ldots,u,u^{-1},u^{-1},u)=(-1)^{l+1}u^{l}$.

\end{enumerate}

\noindent S'il existe $1 \leq j \leq \frac{m-2}{2}$ et $i$ tels que $K_{j}(a_{i},\ldots,a_{j+i-1})=\pm 1_{\mathbb{K}}$ alors $j$ est nécessairement de la forme $3l-1$ et on a :
\[(a_{i},\ldots,a_{j+i-1})=(u^{-1},u^{-1},u,\ldots,u^{-1},u^{-1},u,u^{-1},u^{-1}).\]
\noindent En effet, si $j=3l$ on aurait $u^{l}=\pm 1_{\mathbb{K}}$ et donc $m \leq 3l$ (puisque $m=3v$ où $v$ est le plus petit entier tel que $u^{v}=\pm 1_{\mathbb{K}}$), ce qui est absurde puisque $j=3l \leq \frac{m-2}{2} \leq m-1$. Si $j=3l-2$ alors on aurait $u^{l}=\pm 1_{\mathbb{K}}$ et donc $m \leq 3l$, ce qui est absurde puisque $m-2 \leq j=3l-2 \leq \frac{m-2}{2} \leq m-3$ (car $m \geq 6$). L'entier $j$ est donc nécessairement de la forme $3l-1$. Donc, la solution $u$-trinomiale minimale de $(E_{\mathbb{K}})$ est réductible si et seulement s'il existe un entier naturel non nul $l$ tel que $(-1_{\mathbb{K}})^{l+1}((u^{-1})^{l+1}-u^{l-1}) \neq \pm 1_{\mathbb{K}}$ et $3l-1 \leq \frac{m-2}{2}$. Ainsi, la solution $u$-trinomiale minimale de $(E_{\mathbb{K}})$ est réductible si et seulement s'il existe un entier naturel non nul $1 \leq l \leq E\left[\frac{m}{6}\right]$ tel que $u^{2l} \pm u^{l+1}-1_{\mathbb{K}}=0_{\mathbb{K}}$.

\end{proof}

\begin{remarks}
{\rm i) Dans le cas où ${\rm car}(\mathbb{K})=2$, le théorème prend la forme simplifiée suivante. Soit $u \neq 1_{\mathbb{K}}$, la solution $u$-trinomiale minimale de $(E_{\mathbb{K}})$ est irréductible si et seulement si $u$ n'est pas racine des polynômes $X^{2l}+X^{l+1}+1_{\mathbb{K}}$ pour $2 \leq l \leq E\left[\frac{m}{6}\right]$. De plus, si $l$ est impair, $X^{2l}+X^{l+1}+1_{\mathbb{K}}=(X^{l}+X^{\frac{l+1}{2}}+1_{\mathbb{K}})^{2}$ et il suffit donc de regarder si $u$ est racine de $X^{l}+X^{\frac{l+1}{2}}+1_{\mathbb{K}}$.
\\
\\ii) Si $l$ est pair alors $u$ est racine de $X^{2l}+X^{l+1}-\overline{1}$ si et seulement si $-u$ est racine de $X^{2l}-X^{l+1}-\overline{1}$.
}
\end{remarks}

\begin{examples}
{\rm On donne ci-dessous quelques exemples d'application du théorème précédent.
\begin{itemize}
\item Considérons $\mathbb{K}=\mathbb{Z}/7\mathbb{Z}$ et $u=\overline{2}$. Comme $\overline{2}^{2 \times 1}+\overline{2}^{1+1}-\overline{1}=\overline{0}$, la solution $\overline{2}$-trinomiale minimale de $(E_{\mathbb{K}})$ est réductible.
\\
\item Considérons $\mathbb{K}=\mathbb{F}_{4}=\frac{(\mathbb{Z}/2\mathbb{Z})[X]}{\left\langle X^{2}+X+\overline{1}\right\rangle}$ et $u=X$. L'élément $u$ est d'ordre 3. La solution $X$-trinomiale minimale de $(E_{\mathbb{K}})$ est donc de taille 9. Comme $E\left[\frac{9}{6}\right]=1$, la condition du théorème précédent est automatiquement vérifiée. Ainsi, la solution $X$-trinomiale minimale de $(E_{\mathbb{K}})$ est irréductible.
\\
\item Considérons $\mathbb{K}=\mathbb{Z}/7\mathbb{Z}$ et $u=\overline{4}$. L'élément $u$ est d'ordre 3. La solution $u$-trinomiale minimale de $(E_{\mathbb{K}})$ est donc de taille 9. De plus, on a $\overline{4}^{2 \times 1}+\overline{4}^{1+1}-\overline{1}=\overline{3}$ et $\overline{4}^{2 \times 1}-\overline{4}^{1+1}-\overline{1}=\overline{6}$. Ainsi, par le théorème \ref{42}, la solution $\overline{4}$-trinomiale minimale de $(E_{\mathbb{K}})$ est irréductible.
\\
\item Considérons $\mathbb{K}=\mathbb{Z}/31\mathbb{Z}$ et $u=\overline{3}$. L'élément $u$ est d'ordre 30. La solution $u$-trinomiale minimale de $(E_{\mathbb{K}})$ est donc de taille 45. On a : 
\begin{itemize}[label=$\circ$]
\item $\overline{3}^{2 \times 1}+\overline{3}^{1+1}-\overline{1}=\overline{17}$,~~~~~~~~~~~~$\overline{3}^{2 \times 1}-\overline{3}^{1+1}-\overline{1}=\overline{-1}$;
\item $\overline{3}^{2 \times 2}+\overline{3}^{2+1}-\overline{1}=\overline{14}$,~~~~~~~~~~~~$\overline{3}^{2 \times 2}-\overline{3}^{2+1}-\overline{1}=\overline{-9}$;
\item $\overline{3}^{2 \times 3}+\overline{3}^{3+1}-\overline{1}=\overline{3}$,$\hphantom{0}$~~~~~~~~~~~~$\overline{3}^{2 \times 3}-\overline{3}^{3+1}-\overline{1}=\overline{-4}$;
\item $\overline{3}^{2 \times 4}+\overline{3}^{4+1}-\overline{1}=\overline{14}$,~~~~~~~~~~~~$\overline{3}^{2 \times 4}-\overline{3}^{4+1}-\overline{1}=\overline{-7}$;
\item $\overline{3}^{2 \times 5}+\overline{3}^{5+1}-\overline{1}=\overline{9}$,$\hphantom{0}$~~~~~~~~~~~~$\overline{3}^{2 \times 5}-\overline{3}^{5+1}-\overline{1}=\overline{8}$;
\item $\overline{3}^{2 \times 6}+\overline{3}^{6+1}-\overline{1}=\overline{24}$,~~~~~~~~~~~~$\overline{3}^{2 \times 6}-\overline{3}^{6+1}-\overline{1}=\overline{-10}$;
\item $\overline{3}^{2 \times 7}+\overline{3}^{7+1}-\overline{1}=\overline{29}$,~~~~~~~~~~~~$\overline{3}^{2 \times 7}-\overline{3}^{7+1}-\overline{1}=\overline{-11}$.
\end{itemize}
Ainsi, par le théorème \ref{42}, la solution $\overline{3}$-trinomiale minimale de $(E_{\mathbb{K}})$ est irréductible.

\end{itemize}
}
\end{examples}

\begin{remark}
{\rm Les exemples précédents montrent qu'il peut exister des éléments $u$ dans $\mathbb{K}$ pour lesquels la solution $u$-trinomiale minimale de $(E_{\mathbb{K}})$ est irréductible et la solution $u^{-1}$-trinomiale minimale de $(E_{\mathbb{K}})$ est réductible. En revanche, la solution $u$-trinomiale minimale de $(E_{\mathbb{K}})$ est irréductible si et seulement si la solution $-u$-trinomiale minimale de $(E_{\mathbb{K}})$ est irréductible, puisque, d'une part, $(-u)^{-1}=-u^{-1}$, et que, d'autre part, une solution $(a_{1},\ldots,a_{n})$ est irréductible si et seulement si $(-a_{1},\ldots,-a_{n})$ l'est.
\\
}
\end{remark}

\noindent Le petit résultat ci-dessous nous sera utile dans la suite.

\begin{lemma}
\label{43}

Soient $\mathbb{K}$ un corps fini et $u \in \mathbb{K}^{*}$, $u \neq \pm 1_{\mathbb{K}}$. La solution $u$-trinomiale minimale de $(E_{\mathbb{K}})$ est irréductible si et seulement si $u$ n'est pas racine des polynômes de la forme $X^{2l}+X^{l+1}-1_{\mathbb{K}}$ et $X^{2l}-X^{l+1}-1_{\mathbb{K}}$ pour $l \geq 1$.

\end{lemma}

\begin{proof}

Soit $m$ la taille de la solution $u$-trinomiale minimale. On a $m>3$ car $u \neq \pm 1_{\mathbb{K}}$.
\\
\\Si $u$ n'est pas racine des polynômes de la forme $X^{2l}+X^{l+1}-1_{\mathbb{K}}$ et $X^{2l}-X^{l+1}-1_{\mathbb{K}}$ pour $l \geq 1$ alors, par le théorème \ref{42}, la solution $u$-trinomiale minimale de $(E_{\mathbb{K}})$ est irréductible.
\\
\\S'il existe $l \geq 1$ tel que $u$ est racine de $X^{2l}+X^{l+1}-1_{\mathbb{K}}$ ou de $X^{2l}-X^{l+1}-1_{\mathbb{K}}$. Par les calculs effectués dans la preuve précédente et le lemme \ref{131}, on peut construire une solution de $(E_{\mathbb{K}})$ de la forme $(a,u^{-1},u^{-1},u,\ldots,u^{-1},u^{-1},u,u^{-1},u^{-1},b)$ de taille $l'=3l+1$. On effectue la division euclidienne de $l'$ par $m$. On a $l'=qm+r$ avec $r<m$. De plus, $l' \equiv 1 [3]$ et $m \equiv 0 [3]$. Donc, $r \equiv 1 [3]$. Par ailleurs, si $r=1$ on a $q \geq 1$ puisque $l' \geq 4$.
\\
\\Supposons par l'absurde que $r=1$. Comme $(a,u^{-1},u^{-1},u,\ldots,u^{-1},u^{-1},u,u^{-1},u^{-1},b)$ est solution, $(u,u^{-1},u^{-1},\ldots,u,u^{-1},u^{-1},b,a,u^{-1},u^{-1})$ l'est aussi. Il existe $(\alpha,\beta) \in \{\pm 1_{\mathbb{K}}\}^{2}$ tel que :
\[M_{m}(u,u^{-1},u^{-1},\ldots,u,u^{-1},u^{-1})=\alpha Id,\]
\[M_{qm+1}(u,u^{-1},u^{-1},\ldots,u,u^{-1},u^{-1},b,a,u^{-1},u^{-1})=\beta Id.\]

\begin{eqnarray*}
\beta Id &=& M_{qm+1}(u,u^{-1},u^{-1},\ldots,u,u^{-1},u^{-1},b,a,u^{-1},u^{-1}) \\
         &=& M_{m+1}(\underbrace{u,u^{-1},u^{-1},\ldots,u,u^{-1},u^{-1}}_{m-3},b,a,u^{-1},u^{-1}) \times \\
				 & & M_{(q-1)m}(u,u^{-1},u^{-1},\ldots,u,u^{-1},u^{-1}) \\
				 &=& M_{m+1}(u,u^{-1},u^{-1},\ldots,u,u^{-1},u^{-1},b,a,u^{-1},u^{-1}) (\alpha^{q-1} Id) \\
				 &=& \alpha^{q-1} M_{4}(b,a,u^{-1},u^{-1})M_{m-3}(u,u^{-1},u^{-1},\ldots,u,u^{-1},u^{-1}).
\end{eqnarray*}

\noindent Donc, il existe $\mu \in \{\pm 1_{\mathbb{K}}\}$ tel que :
\begin{eqnarray*}
M &=& M_{4}(b,a,u^{-1},u^{-1})M_{m-3}(u,u^{-1},u^{-1},\ldots,u,u^{-1},u^{-1}) \\
  &=& M_{2}(u^{-1},u^{-1})M_{2}(b,a)M_{m-3}(u,u^{-1},u^{-1},\ldots,u,u^{-1},u^{-1}) \\
  &=& \mu M_{m}(u,u^{-1},u^{-1},\ldots,u,u^{-1},u^{-1}) \\
	&=& \mu M_{2}(u^{-1},u^{-1})M_{1}(u)M_{m-3}(u,u^{-1},u^{-1},\ldots,u,u^{-1},u^{-1}).
\end{eqnarray*}

\noindent Ainsi, en multipliant à droite par l'inverse de $M_{m-3}(u,u^{-1},u^{-1},\ldots,u,u^{-1},u^{-1})$, et à gauche par l'inverse de $M_{2}(u^{-1},u^{-1})$, on obtient :
\[M_{2}(b,a)=\begin{pmatrix}
   ab-1_{\mathbb{K}} & -b  \\
   a  & -1_{\mathbb{K}}
		\end{pmatrix}=\mu \begin{pmatrix}
   u & -1_{\mathbb{K}}  \\
   1_{\mathbb{K}}  & 0_{\mathbb{K}}
		\end{pmatrix}=\mu M_{1}(u).\]
\noindent On en déduit $0_{\mathbb{K}}=-1_{\mathbb{K}}$, ce qui est absurde. Donc, $r \neq 1$ et $r \geq 4$.
\\
\\En supprimant $q$ fois le $m$-uplet $(u,u^{-1},u^{-1},\ldots,u,u^{-1},u^{-1})$ dans le $l'$-uplet d'éléments du corps $\mathbb{K}$ $(a,u^{-1},u^{-1},u,\ldots,u^{-1},u^{-1},u,u^{-1},u^{-1},b)$, on obtient une solution de taille $r \geq 4$ avec laquelle on peut réduire la solution $u$-trinomiale minimale.
                             
\end{proof}

On dispose donc de conditions nécessaires et suffisantes simples pour tester l'irréductibilité de ces solutions (voir l'annexe \ref{C} pour quelques applications numériques). Notons par ailleurs que ce résultat fait apparaître des trinômes sur des corps finis. On dispose de très nombreux résultats sur ces polynômes particuliers et on en cite quelques-uns ci-dessous :
\begin{itemize}
\item (Théorème de Swan-Stickelberger, \cite{S} corollaire 5) Soit $n>k>0$. On suppose que parmi $n$ et $k$ il y a exactement un entier impair. $X^{n}+X^{k}+\overline{1}$ a un nombre pair de facteurs irréductibles (et est donc en particulier réductible) sur $\mathbb{Z}/2\mathbb{Z}$ si et seulement si on est dans une des situations ci-dessous :
\begin{itemize}[label=$\circ$] 
\item $n$ est pair, $k$ est impair, $n \neq 2k$ et $\frac{nk}{2} \equiv 0,1 [4]$;
\item $n$ est impair, $k$ est pair, $k$ ne divise pas $2n$ et $n \equiv \pm 3 [8]$;
\item $n$ est impair, $k$ est pair, $k$ divise $2n$ et $n \equiv \pm 1 [8]$.
\\
\end{itemize}
\item (Théorème de Kelley-Owen, \cite{KO} Théorème 1.2) Soient $q>1$ la puissance d'un nombre premier et $(a,b) \in (\mathbb{F}_{q}^{*})^{2}$. Soient $n>s>0$ et $\delta:={\rm pgcd}(n,s,q-1)$. Le trinôme $X^{n}+aX^{s}+b$ a au plus $\delta~E\left[\frac{1}{2}+\sqrt{\frac{q-1}{\delta}}\right]$ racines sur $\mathbb{F}_{q}$.
\\
\item (\cite{KK} corollaire 4) Le polynôme $X^{n}+X^{k}+\overline{1} \in (\mathbb{Z}/2\mathbb{Z})[X]$ ($n>k>0$) est divisible par $X^{2}+X+\overline{1}$ si et seulement si $n$, $k$ et $n-k$ ne sont pas divisibles par 3.
\\
\end{itemize}

Le théorème \ref{42} apporte des éléments intéressants. Toutefois, les conditions requises par l'énoncé deviennent rapidement complexes à manier et ne sont donc pas particulièrement adaptées à la recherche de minorations sur $\ell_{\mathbb{K}}$, qui est un des objectifs majeurs de cette section. Par conséquent, on va retravailler les propriétés demandées dans l'énoncé afin de rendre ces dernières plus facilement exploitables. Plus précisément, on va démontrer le résultat ci-dessous :

\begin{theorem}
\label{44}

Soient $\mathbb{K}$ un corps fini de caractéristique différente de 2 et $u \in \mathbb{K}^{*}$. On suppose que $u \neq \pm 1_{\mathbb{K}}$.
\\
\\i) Si $u^{2}+4_{\mathbb{K}}$ n'est pas un carré sur $\mathbb{K}$ alors la solution $u$-trinomiale minimale de $(E_{\mathbb{K}})$ est irréductible.
\\
\\ii) On suppose que $u$ ou $-u$ est un générateur du groupe $\mathbb{K}^{*}$. La solution $u$-trinomiale minimale de $(E_{\mathbb{K}})$ est irréductible si et seulement si $u^{2}+4_{\mathbb{K}}$ n'est pas un carré sur $\mathbb{K}$.

\end{theorem} 

\begin{proof}

Soit $m$ la taille de la solution $u$-trinomiale minimale de $(E_{\mathbb{K}})$.
\\
\\i) On raisonne par contraposée. Supposons que la solution $u$-trinomiale minimale de $(E_{\mathbb{K}})$ est réductible. Par le théorème \ref{42}, il existe $1 \leq l \leq E\left[\frac{m}{6}\right]$ tel que $u^{2l}+u^{l+1}-1_{\mathbb{K}}=(u^{l})^{2}+u(u^{l})-1_{\mathbb{K}}=0_{\mathbb{K}}$ ou $u^{2l}-u^{l+1}-1_{\mathbb{K}}=(u^{l})^{2}-u(u^{l})-1_{\mathbb{K}}=0_{\mathbb{K}}$. Donc, l'un des deux trinômes $X^{2}+uX-1_{\mathbb{K}}$ ou $X^{2}-uX-1_{\mathbb{K}}$ possède une racine sur $\mathbb{K}$. Donc, puisque $\mathbb{K}$ est de caractéristique différente de 2, le discriminant de ces polynômes $\Delta=(\pm u)^{2}+4_{\mathbb{K}}=u^{2}+4_{\mathbb{K}}$ est un carré sur $\mathbb{K}$.
\\
\\ii) On suppose que $u$ est un générateur du groupe $\mathbb{K}^{*}$.
\\
\\ Si $u^{2}+4_{\mathbb{K}}$ n'est pas un carré sur $\mathbb{K}$ alors, par i), la solution $u$-trinomiale minimale de $(E_{\mathbb{K}})$ est irréductible. 
\\
\\Si $u^{2}+4_{\mathbb{K}}$ est un carré sur $\mathbb{K}$. Le discriminant du trinôme $P(X)=X^{2}+uX-1_{\mathbb{K}}$ est un carré sur $\mathbb{K}$. Donc, $P$ possède des racines sur $\mathbb{K}$, c'est-à-dire qu'il existe $y \in \mathbb{K}$ tel que $y^{2}+uy-1_{\mathbb{K}}=0_{\mathbb{K}}$ (puisque le corps est de caractéristique différente de 2). De plus, $y \neq 0_{\mathbb{K}}$ car $P(0_{\mathbb{K}})=-1_{\mathbb{K}}$ et $y \neq \pm 1_{\mathbb{K}}$ car $P(\pm 1_{\mathbb{K}})=\pm u$. Ainsi, comme $u$ est un générateur de $\mathbb{K}^{*}$, il existe $1 \leq l \leq {\rm card}(\mathbb{K})-1$ tel que $y=u^{l}$ et donc $u^{2l}+u^{l+1}-1_{\mathbb{K}}=0_{\mathbb{K}}$. 
\\
\\Par le lemme \ref{43}, la solution $u$-trinomiale minimale de $(E_{\mathbb{K}})$ est réductible.
\\
\\Si $-u$ est un générateur du groupe $\mathbb{K}^{*}$ alors, par ce qui précède, la solution $-u$-trinomiale minimale de $(E_{\mathbb{K}})$ est irréductible si et seulement si $(-u)^{2}+4_{\mathbb{K}}=u^{2}+4_{\mathbb{K}}$ n'est pas un carré. Or, la solution $u$-trinomiale minimale de $(E_{\mathbb{K}})$ est irréductible si et seulement si la solution $-u$-trinomiale minimale de $(E_{\mathbb{K}})$ est irréductible. On obtient donc le résultat souhaité.

\end{proof}

\begin{examples}
{\rm On donne ci-dessous quelques exemples d'application du théorème précédent.
\begin{itemize}
\item Considérons $\mathbb{K}=\mathbb{Z}/13\mathbb{Z}$ et $u=\overline{2}$. L'ensemble des carrés modulo 13 est $\{\overline{0},\overline{1},\overline{4},\overline{9},\overline{3},\overline{-1},\overline{10}\}$. Comme $\overline{2}^{2}+\overline{4}=\overline{8}$ n'est pas un carré modulo 13, la solution $\overline{2}$-trinomiale minimale de $(E_{\mathbb{K}})$ est irréductible.
\\
\item Considérons $\mathbb{K}=\mathbb{Z}/29\mathbb{Z}$ et $u=\overline{9}$. $\overline{9}^{2}+\overline{4}=\overline{85}=\overline{27}$ n'est pas un carré modulo 29. En effet, $\left(\dfrac{27}{29}\right)=\left(\dfrac{3}{29}\right)^{2}\left(\dfrac{3}{29}\right)=(-1)$ (proposition \ref{18}). Ainsi, la solution $\overline{9}$-trinomiale minimale de $(E_{\mathbb{K}})$ est irréductible.
\\
\item Considérons $\mathbb{K}=\mathbb{Z}/19\mathbb{Z}$ et $u=\overline{10}$. $u$ est un générateur de $\mathbb{K}^{*}$. $\overline{10}^{2}+\overline{4}=\overline{104}=\overline{9}=\overline{3}^{2}$. Donc, la solution $\overline{10}$-trinomiale minimale de $(E_{\mathbb{K}})$ est réductible.
\\
\item Considérons $\mathbb{K}=\mathbb{Z}/23\mathbb{Z}$ et $u=\overline{18}$. $u$ n'est pas un générateur de $\mathbb{K}^{*}$ mais $-u$ en est un. $\overline{5}^{2}+\overline{4}=\overline{29}=\overline{6}=\overline{11}^{2}$. Donc, la solution $\overline{18}$-trinomiale minimale de $(E_{\mathbb{K}})$ est réductible.
\\
\item Considérons $\mathbb{K}=\mathbb{F}_{9}=\frac{(\mathbb{Z}/3\mathbb{Z})[X]}{\left\langle X^{2}+\overline{1}\right\rangle}$ et $u=X+\overline{1}$. $u^{2}+\overline{4}=-X+\overline{1}$ n'est pas un carré dans $\mathbb{F}_{9}$ (voir annexe \ref{A}). Ainsi, la solution $u$-trinomiale minimale de $(E_{\mathbb{K}})$ est irréductible.
\\
\item Considérons $\mathbb{K}=\mathbb{F}_{49}=\frac{(\mathbb{Z}/7\mathbb{Z})[X]}{\left\langle X^{2}+\overline{1}\right\rangle}$ et $u=\overline{3}X+\overline{2}$. On a $u^{2}+\overline{4}=\overline{5}X-\overline{1}$ et $(\overline{5}X-\overline{1})^{8}=-\overline{2}$. Comme $-\overline{2}$ n'est pas un carré dans $\mathbb{Z}/7\mathbb{Z}$, $\overline{5}X-\overline{1}$ n'est pas un carré dans $\mathbb{F}_{49}$, par la proposition \ref{19}. Ainsi, la solution $u$-trinomiale minimale de $(E_{\mathbb{K}})$ est irréductible.
\end{itemize}
}
\end{examples}

\begin{remarks}
{\rm i) Ce théorème est faux si le corps considéré est de caractéristique 2, puisque dans ce cas tous les éléments sont des carrés.
\\
\\ii) La réciproque de i) est fausse. Considérons par exemple $\mathbb{K}=\mathbb{Z}/13\mathbb{Z}$ et $u=\overline{5}$. $u^{2}+\overline{4}=\overline{3}=\overline{4}^{2}$. Or, la solution $\overline{5}$-trinomiale minimale est irréductible sur $\mathbb{Z}/13\mathbb{Z}$. En effet, celle-ci est de taille 6 et $\overline{5}^{2 \times 1}-\overline{5}^{1+1}-\overline{1}=-\overline{1}$ et $\overline{5}^{2 \times 1}+\overline{5}^{1+1}-\overline{1}=-\overline{3}$. Par le théorème \ref{42}, elle est donc bien irréductible.
}
\end{remarks}

Le théorème \ref{44} est beaucoup plus maniable que le précédent, comme le montre les exemples précédents et les quelques applications simples présentées ci-dessous :

\begin{proposition}
\label{45}

Soit $\mathbb{K}$ un corps fini de cardinal $q$ et de caractéristique différente de 2. Il y a au moins $\frac{q-1}{2}$ éléments $u$ de $\mathbb{K}$ pour lesquels la solution $u$-trinomiale minimale de $(E_{\mathbb{K}})$ est irréductible. Si $\mathbb{K}=\mathbb{Z}/p\mathbb{Z}$ avec $p \equiv 3 [4]$ alors il y a au moins $\frac{p+1}{2}$ éléments $\overline{u}$ pour lesquels la solution $u$-trinomiale minimale de $(E_{\mathbb{K}})$ est irréductible.

\end{proposition}

\begin{proof}

Par la proposition \ref{14}, il y a au moins $\frac{q-1}{2}$ éléments $u$ de $\mathbb{K}$ pour lesquels $u^{2}+4_{\mathbb{K}}$ n'est pas un carré sur $\mathbb{K}$. Parmi ceux-ci, il n'y a pas $0_{\mathbb{K}}$. Donc, il y a au moins $\frac{q-1}{2}$ éléments inversibles $u$ de $\mathbb{K}$ pour lesquels $u^{2}+4_{\mathbb{K}}$ n'est pas un carré sur $\mathbb{K}$. Par le théorème \ref{44} (et l'irréductibilité de la solution $\pm 1_{\mathbb{K}}$-trinomiale minimale), chacun de ces $u$ donne une solution trinomiale minimale de $(E_{\mathbb{K}})$ irréductible.
\\
\\Si $p \equiv 3 [4]$ alors, par la proposition \ref{15}, $\left(\dfrac{-4}{p}\right)=\left(\dfrac{-1}{p}\right)\left(\dfrac{2}{p}\right)^{2}=(-1)^{\frac{p-1}{2}}=-1$. Donc, par la proposition \ref{14}, il y a $\frac{p+1}{2}$ éléments $\overline{u}$ de $\mathbb{Z}/p\mathbb{Z}$ pour lesquels $\overline{u}^{2}+\overline{4}$ n'est pas un carré modulo $p$. Tous ces éléments sont inversibles car $\overline{0}$ n'est pas parmi eux. Donc, par le théorème \ref{44} (et l'irréductibilité de la solution $\pm \overline{1}$-trinomiale minimale), chacun de ces éléments donne une solution trinomiale minimale de $(E_{\mathbb{Z}/p\mathbb{Z}})$ irréductible.

\end{proof}

\begin{proposition}
\label{46}

Soit $p$ un nombre premier impair différent de 3.
\\
\\i) Si $p \equiv \pm 3 [8]$ alors la solution $\pm \overline{2}$-trinomiale minimale de $(E_{\mathbb{Z}/p\mathbb{Z}})$ est irréductible.
\\
\\ii) Si $p=5$ ou si $p \equiv \pm 2 [5]$ alors la solution $\pm \overline{4}$-trinomiale minimale de $(E_{\mathbb{Z}/p\mathbb{Z}})$ est irréductible.

\end{proposition}

\begin{proof}

i) On suppose $p \neq 3$ et $p \equiv \pm 3 [8]$. On a $\pm \overline{2} \neq \pm \overline{1}$. Par ailleurs, $(\pm \overline{2})^{2}+\overline{4}=\overline{8}$ et $\left(\dfrac{8}{p}\right)=\left(\dfrac{2}{p}\right)\left(\dfrac{2}{p}\right)^{2}=(-1)^{\frac{p^{2}-1}{8}}=-1$ (Loi complémentaire de réciprocité quadratique, Théorème \ref{17}). Donc, par le théorème \ref{44}, la solution $\pm \overline{2}$-trinomiale minimale de $(E_{\mathbb{Z}/p\mathbb{Z}})$ est irréductible.
\\
\\ii) On suppose $p \neq 5$. On a $\pm \overline{4} \neq \pm \overline{1}$ et $(\pm \overline{4})^{2}+\overline{4}=\overline{20}$. Par la loi de réciprocité quadratique de Gauss (Théorème \ref{16}), on a : \[\left(\dfrac{20}{p}\right)=\left(\dfrac{2}{p}\right)^{2}\left(\dfrac{5}{p}\right)=(-1)^{\frac{5-1}{2}\frac{p-1}{2}}\left(\dfrac{p}{5}\right)=(-1)^{p-1}\left(\dfrac{p}{5}\right)=\left(\dfrac{p}{5}\right).\]
\noindent L'ensemble des carrés modulo 5 est $\{\overline{0},\overline{1},\overline{4}\}$. Donc, si $p \equiv \pm 2 [5]$, on a $\left(\dfrac{20}{p}\right)=-1$ et, par le théorème \ref{44}, la solution $\pm \overline{4}$-trinomiale minimale de $(E_{\mathbb{Z}/p\mathbb{Z}})$ est irréductible.
\\
\\Si $p=5$ alors $\overline{4}=\overline{-1}$ et la solution $\pm \overline{4}$-trinomiale minimale de $(E_{\mathbb{Z}/p\mathbb{Z}})$ est irréductible de taille 3.

\end{proof}

\begin{example}

{\rm  Comme $101 \equiv -3 [8]$, la solution $\overline{2}$-trinomiale minimale de $(E_{\mathbb{Z}/101\mathbb{Z}})$ est irréductible.
}

\end{example}

\subsection{Recherche de minorants}

On va maintenant revenir à la question sous-tendant toute cette partie : peut-on utiliser la famille des solutions trinomiales minimales pour obtenir une amélioration de la minoration de $\ell_{\mathbb{K}}$ ? On commence par démontrer le premier point du théorème \ref{07}. Toutefois, avant de faire cela, on a besoin de plusieurs résultats préliminaires.

\begin{lemma}
\label{47}

Soient $p \geq 3$, $q \geq 5$ et $(n,m) \in (\mathbb{N}^{*})^{2}$. $2(p^{n}-p^{n-1})(q^{m}-q^{m-1})>p^{n}q^{m}$.

\end{lemma}

\begin{proof}

On a :
\begin{eqnarray*}
2(p^{n}-p^{n-1})(q^{m}-q^{m-1})-p^{n}q^{m} &=& p^{n-1}q^{m-1}[2(p-1)(q-1)-pq] \\
																					 &=& p^{n-1}q^{m-1}[pq-2p-2q+2] \\
																					 &=& p^{n-1}q^{m-1}[(p-2)(q-2)-2].
\end{eqnarray*}

\noindent Or, $p-2 \geq 1$ et $q-2 \geq 3$. Donc, $(p-2)(q-2)-2 \geq 1$. Par conséquent, on a :
\[2(p^{n}-p^{n-1})(q^{m}-q^{m-1})-p^{n}q^{m}>0.\]
\noindent Donc, $2(p^{n}-p^{n-1})(q^{m}-q^{m-1})>p^{n}q^{m}$.

\end{proof}

\begin{proposition}
\label{48}

i) Si $p$ est un nombre premier de la forme $2q^{n}+1$ avec $q$ premier impair et $n \in \mathbb{N}^{*}$ alors $\ell_{\mathbb{Z}/p\mathbb{Z}} \geq 3 \times \frac{p-1}{2}$.
\\
\\ii) Si $p$ est un nombre premier de la forme $2q^{n}r^{m}+1$ avec $q$ et $r$ des nombres premiers impairs distincts et $(n,m) \in (\mathbb{N}^{*})^{2}$ alors $\ell_{\mathbb{Z}/p\mathbb{Z}} \geq 3 \times \frac{p-1}{2}$.

\end{proposition}

\begin{proof}

Notons $\mathcal{G'}_{p}$ l'ensemble des éléments du groupe $(\mathbb{Z}/p\mathbb{Z})^{*}$ d'ordre $p-1$ ou $\frac{p-1}{2}$ et $\mathcal{C}_{p}$ l'ensemble $\{\overline{x} \in \mathbb{Z}/p\mathbb{Z}, \overline{x}^{2}+\overline{4}~{\rm est~un~carr\acute{e}~modulo}~p\}$.
\\
\\i) Comme $(\mathbb{Z}/p\mathbb{Z})^{*}$ est cyclique, ${\rm card}(\mathcal{G'}_{p})=\varphi(p-1)+\varphi\left(\frac{p-1}{2}\right)=2(q^{n}-q^{n-1})$. De plus, puisque $q$ est premier, $q \equiv \epsilon [4]$ avec $\epsilon=\pm 1$ et $2q^{n}+1 \equiv 2 \epsilon^{n}+1 \equiv 3 [4]$. Donc, $-\overline{4}$ n'est pas un carré modulo $p$ (proposition \ref{15}) et, par la proposition \ref{14}, ${\rm card}(\mathcal{C}_{p})=\frac{p-1}{2}=q^{n}$. Ainsi, puisque $q \geq 3$, on a :
\[{\rm card}(\mathcal{G'}_{p})=q^{n}+q^{n-1}(q-2)>q^{n}={\rm card}(\mathcal{C}_{p}).\]
\noindent Donc, il existe un élément $\overline{u}$ de $\mathbb{Z}/p\mathbb{Z}$, différent de $\pm \overline{1}$, tel que $o(\overline{u}) \in \{p-1, \frac{p-1}{2}\}$ et $\overline{u}^{2}+\overline{4}$ n'est pas un carré modulo $p$. Par les théorèmes \ref{42} et \ref{44}, la solution $\overline{u}$-trinomiale minimale de $(E_{\mathbb{Z}/p\mathbb{Z}})$ est irréductible de taille $3 \times \frac{p-1}{2}$ (car $\frac{p-1}{2}$ est impair). On obtient ainsi la minoration souhaitée.
\\
\\ii) Comme $(\mathbb{Z}/p\mathbb{Z})^{*}$ est cyclique, ${\rm card}(\mathcal{G'}_{p})=\varphi(p-1)+\varphi\left(\frac{p-1}{2}\right)=2(q^{n}-q^{n-1})(r^{m}-r^{m-1})$. De plus, puisque $q$ et $r$ sont premiers, $q,r \equiv \pm 1 [4]$. Donc, $q^{n}r^{m} \equiv \pm 1 [4]$ et $p \equiv 3 [4]$. Donc, $-\overline{4}$ n'est pas un carré modulo $p$ (proposition \ref{15}) et, par la proposition \ref{14}, ${\rm card}(\mathcal{C}_{p})=\frac{p-1}{2}=q^{n}r^{m}$. Ainsi, par le lemme \ref{47}, on a :
\[{\rm card}(\mathcal{G'}_{p})>{\rm card}(\mathcal{C}_{p}).\]
\noindent Donc, il existe un élément $\overline{u}$ de $\mathbb{Z}/p\mathbb{Z}$, différent de $\pm \overline{1}$, tel que $o(\overline{u}) \in \{p-1, \frac{p-1}{2}\}$ et $\overline{u}^{2}+\overline{4}$ n'est pas un carré modulo $p$. Par les théorèmes \ref{42} et \ref{44}, la solution $\overline{u}$-trinomiale minimale de $(E_{\mathbb{Z}/p\mathbb{Z}})$ est irréductible de taille $3 \times \frac{p-1}{2}$ (car $\frac{p-1}{2}$ est impair). On obtient ainsi la minoration souhaitée.

\end{proof}

\begin{remark}
{\rm L'argument de i) ne peut pas se généraliser aux nombres premiers de la forme $p=2^{n}q^{m}+1$ avec $n \geq 2$, puisque dans ce cas $\frac{p-1}{2}$ est pair. L'argument de ii) ne peut pas non plus se généraliser aux nombres premiers $p$ pour lesquels $\frac{p-1}{2}$ a au moins trois facteurs premiers car dans ce cas on n'a plus nécessairement ${\rm card}(\mathcal{G'}_{p})>{\rm card}(\mathcal{C}_{p})$. En effet, prenons, par exemple, $p=211=2\times(3 \times 5 \times 7)+1$. On a ${\rm card}(\mathcal{G'}_{p})=96$ et ${\rm card}(\mathcal{C}_{p})=105$.
}
\end{remark}

\begin{examples}
{\rm \begin{itemize}
\item $487=2\times 3^{5}+1$ est premier et $\ell_{\mathbb{Z}/487\mathbb{Z}} \geq 729$;
\item $12~251=2\times 5^{3} \times 7^{2}+1$ est premier et $\ell_{\mathbb{Z}/12~251\mathbb{Z}} \geq 18~375$.
\end{itemize}
}
\end{examples}

\noindent On peut également s'intéresser au cas très particulier où $\overline{2}$ est un générateur de $(\mathbb{Z}/p\mathbb{Z})^{*}$.

\begin{proposition}
\label{481}

Soit $p$ un nombre premier impair. Si $\overline{2}$ est un générateur de $(\mathbb{Z}/p\mathbb{Z})^{*}$ alors la solution $\overline{2}$-trinomiale minimale de $(E_{\mathbb{Z}/p\mathbb{Z}})$ est irréductible. En particulier : 
\[\ell_{\mathbb{Z}/p\mathbb{Z}} \geq 3 \times \frac{p-1}{2}.\]

\end{proposition}

\begin{proof}

Si $p=3$ alors $\overline{2}=-\overline{1}$ et la solution $\overline{2}$-trinomiale minimale est irréductible de taille $3$. On suppose maintenant $p>3$. On a $\overline{2}^{2}+\overline{4}=\overline{2} \times (\overline{2})^{2}$. Or, $\overline{2}$ étant un générateur de $(\mathbb{Z}/p\mathbb{Z})^{*}$, $\overline{2}$ n'est pas un carré. Par conséquent, $\overline{2} \times (\overline{2})^{2}$ n'est pas un carré. Ainsi, par le théorème \ref{44}, la solution $\overline{2}$-trinomiale minimale est irréductible de taille $3 \times \frac{p-1}{2}$, ce qui donne la minoration souhaitée.

\end{proof}

\begin{examples}
{\rm \begin{itemize}
\item Prenons $p=661$. On a $p-1=660=4\times 3 \times 5 \times 11$. Cette décomposition empêche d'utiliser la proposition \ref{48}. En revanche, $\overline{2}$ engendre $(\mathbb{Z}/p\mathbb{Z})^{*}$, ce qui implique, par la proposition précédente, $\ell_{\mathbb{Z}/661\mathbb{Z}} \geq 990$.\\
\item Prenons $p=773$. On a $p-1=772=4\times 193$. On ne peut donc pas utiliser la proposition \ref{48}. En revanche, $\overline{2}$ engendre $(\mathbb{Z}/p\mathbb{Z})^{*}$, ce qui implique, par la proposition précédente, $\ell_{\mathbb{Z}/773\mathbb{Z}} \geq 1158$.
\end{itemize}
}
\end{examples}

\begin{remark}
{\rm Il n'existe pas de résultat permettant de savoir, pour un nombre premier $p$ donné, si $\overline{2}$ est ou non un générateur de $(\mathbb{Z}/p\mathbb{Z})^{*}$. On ignore même s'il existe une infinité de nombres premiers $p$ pour lesquels 2 serait un générateur du groupe des inversibles modulo $p$. Une célèbre conjecture d'Artin (appliquée ici au cas de 2) énonce qu'il y aurait une infinité de tels nombres, mais celle-ci n'a jusqu'à présent pas été démontrée. Notons toutefois qu'une preuve  utilisant l'hypothèse de Riemann généralisée existe (voir \cite{H}). Par ailleurs, D. R. Heath-Brown a démontré le résultat suivant (voir \cite{HB}). Soit $S_{p}$ l'ensemble $\{q \in \mathbb{P},~\overline{p}~{\rm engendre}~(\mathbb{Z}/q\mathbb{Z})^{*}\}$. On a $\left|\{p \in \mathbb{P},~\left|S_{p}\right|<+\infty\}\right| \leq 2$. En particulier, il y a un entier parmi 2, 3 et 5  qui engendre $(\mathbb{Z}/p\mathbb{Z})^{*}$ pour une infinité de nombres premiers $p$. Notons pour terminer qu'une liste des nombres premiers $p$ pour lesquels $\overline{2}$ engendre $(\mathbb{Z}/p\mathbb{Z})^{*}$ est disponible dans \cite{OEIS} A001122.
}
\end{remark}

\noindent On aura également besoin du résultat général suivant.

\begin{theorem}[\cite{BCST} Théorème 1]
\label{49}

Soient $q>1$ la puissance d'un nombre premier, $\mathbb{K}:=\mathbb{F}_{q}$ et $\Omega$ l'ensemble ci-dessous :
\[\Omega:=\{2, 3, 4, 5, 7, 9, 11, 13, 16, 19, 23, 25, 29, 31, 37, 41, 43, 49, 61, 67, 71, 73, 79, 121, 127, 151, 211\}.\]
On suppose que $q \notin \Omega$. Soit $(a,b,c) \in (\mathbb{F}_{q})^{3}$ tel que $a \neq 0_{\mathbb{K}}$ et $b^{2}-4_{\mathbb{K}}ac \neq 0_{\mathbb{K}}$. Il existe un générateur $g$ du groupe $\mathbb{K}^{*}$ tel que $ag^{2}+bg+c$ est un générateur de $\mathbb{K}^{*}$.

\end{theorem}

\noindent On peut maintenant énoncer et démontrer la minoration annoncée.

\begin{theorem}
\label{410}

Soient $q>1$ la puissance d'un nombre premier impair et $\mathbb{K}:=\mathbb{F}_{q}$. On a :
\[\ell_{\mathbb{K}} \geq \frac{3(q-1)}{2}.\]

\end{theorem}

\begin{proof}

Soit $\Theta:=\{{\color{red} 3}, {\color{red} 5}, {\color{blue} 7}, 9, {\color{blue} 11}, {\color{red} 13}, {\color{red} 19}, {\color{blue} 23}, 25, {\color{red} 29}, {\color{blue} 31}, {\color{red} 37}, 41, {\color{blue} 43}, 49, {\color{red} 61}, {\color{blue} 67}, {\color{blue} 71}, 73, {\color{blue} 79}, 121, {\color{blue} 127}, {\color{blue} 151}, {\color{red} 211}\}$. 
\\
\\Supposons pour le moment que $q \notin \Theta$. On a $0_{\mathbb{K}}^{2}-4_{\mathbb{K}} \times 1_{\mathbb{K}} \times 4_{\mathbb{K}}=-16_{\mathbb{K}} \neq 0_{\mathbb{K}}$, puisque $q$ est impair. Donc, par le théorème \ref{49}, il existe un générateur $u$ de $\mathbb{K}^{*}$ tel que $u^{2}+4_{\mathbb{K}}$ est un générateur de $\mathbb{K}^{*}$. En particulier, $u^{2}+4_{\mathbb{K}}$ n'est pas un carré et $u \neq \pm 1_{\mathbb{K}}$ (car $q>3$). Ainsi, par les théorèmes \ref{42} et \ref{44}, la solution $u$-trinomiale minimale de $(E_{\mathbb{K}})$ est irréductible de taille $\frac{3(q-1)}{2}$, ce qui implique le résultat souhaité.
\\
\\Supposons maintenant que  $q \in \Theta$. Pour les nombres premiers $p$ écrits en {\color{blue} bleu}, la minoration de $\ell_{\mathbb{Z}/p\mathbb{Z}}$ résulte directement la proposition \ref{48}. Pour les nombres premiers $p$ écrits en {\color{red} rouge}, on vérifie que $\overline{2}$ est un générateur de $(\mathbb{Z}/p\mathbb{Z})^{*}$ et on applique la proposition \ref{481}. Pour les cas restants, on utilise le tableau ci-dessous :

\begin{center}
\begin{tabular}{|c|c|c|c|c|c|c|}
\hline
  $\mathbb{K}$     & $\mathbb{F}_{9}=\frac{(\mathbb{Z}/3\mathbb{Z})[X]}{\left\langle X^{2}+\overline{1}\right\rangle}$ & $\mathbb{F}_{25}=\frac{(\mathbb{Z}/5\mathbb{Z})[X]}{\left\langle X^{2}+X+\overline{1}\right\rangle}$  & $\mathbb{Z}/41\mathbb{Z}$ & $\mathbb{F}_{49}=\frac{(\mathbb{Z}/7\mathbb{Z})[X]}{\left\langle X^{2}+\overline{1}\right\rangle}$ & $\mathbb{Z}/73\mathbb{Z}$ & $\mathbb{F}_{121}=\frac{(\mathbb{Z}/11\mathbb{Z})[X]}{\left\langle X^{2}-\overline{2}\right\rangle}$  \rule[-7pt]{0pt}{18pt} \\
	\hline
	  $u$ & $X+\overline{1}$ & $X+\overline{2}$ & $\overline{7}$ & $X+\overline{3}$ & $\overline{5}$ &  $X+\overline{2}$ \rule[-7pt]{0pt}{18pt} \\
	\hline
	
\end{tabular}
\end{center}

Dans chacun des cas, $u$ est un générateur de $\mathbb{K}^{*}$ et, grâce aux résultats de la section \ref{carré} et de l'annexe \ref{A}, on montre que $u^{2}+4_{\mathbb{K}}$ n'est pas un carré. En utilisant le théorème \ref{44}, on obtient la minoration souhaitée.

\end{proof}

Notons que cette minoration est optimale puisqu'on a obtenu $\ell_{\mathbb{Z}/5\mathbb{Z}}=6$ et, informatiquement, $\ell_{\mathbb{Z}/7\mathbb{Z}}=9$ (voir \cite{M1}) . En revanche, il n'y a en général pas égalité. Par exemple, on a obtenu informatiquement $\ell_{\mathbb{F}_{9}}=15>\frac{3(9-1)}{2}=12$ (voir annexe \ref{F}).

\begin{example}
{\rm On a $\ell_{\mathbb{F}_{625}} \geq 936$.
}
\end{example}

\subsection{Le cas des corps de caractéristique 2}
\label{cardeux}

Les éléments développés dans la sous-partie précédente nous ont permis d'obtenir une minoration intéressante pour les corps de caractéristique impaire. Pour achever la preuve du théorème \ref{07}, il ne nous reste donc plus qu'à considérer les corps de caractéristique 2. Cependant, pour ces derniers on ne peut plus utiliser les méthodes précédentes basées sur les carrés. On va donc adopter une autre approche que l'on va présenter dans cette sous-section. Celle-ci nécessite de rappeler les deux résultats ci-dessous.

\begin{proposition}
\label{412}

Soient $n \in \mathbb{N}^{*}$ et $(a,b) \in ((\mathbb{F}_{2^{n}})^{*})^{2}$. L'application suivante est bien définie.
\[\begin{array}{ccccc} 
T : & \mathbb{F}_{2^{n}} & \longrightarrow & \mathbb{F}_{2} \\
  & x & \longmapsto & \sum_{i=0}^{n-1} x^{2^{i}}  \\
\end{array}.\]
$X^{2}+aX+b$ est scindé sur $\mathbb{F}_{2^{n}}$ si et seulement si $T(ba^{-2})=\overline{0}$.

\end{proposition}

Ce résultat semble plutôt classique. Toutefois, nous n'avons pas pu localiser de références précises dans la littérature. Aussi, on va effectuer la preuve complète de celui-ci.

\begin{proof}

Soit $x \in \mathbb{F}_{2^{n}}$. On a $x^{2^{n}}=x$. En effet, cette égalité est immédiate pour $x=\overline{0}$ et pour $x \neq \overline{0}$ elle découle du fait que $(\mathbb{F}_{2^{n}})^{*}$ est un groupe de cardinal $2^{n}-1$. Donc, en utilisant le morphisme de Frobenius, on a $(\sum_{i=0}^{n-1} x^{2^{i}})^{2}=\sum_{i=0}^{n-1} x^{2^{i+1}}=\sum_{j=0}^{n-1} x^{2^{j}}$. Ainsi, $T(x)$ est une racine de $Y^{2}-Y$, c'est-à-dire $T(x) \in \mathbb{F}_{2}$. $T$ est donc bien définie.
\\
\\On considère l'application : \[\begin{array}{ccccc} 
f : & \mathbb{F}_{2^{n}} & \longrightarrow & \mathbb{F}_{2^{n}} \\
  & x & \longmapsto & x+x^{2}  \\
\end{array}.\]

\noindent En utilisant le morphisme de Frobenius, on constate que $T$ et $f$ sont des morphismes de groupe et que $T(x^{2})=T(x)$. De plus, ${\rm Ker}(f)=\mathbb{F}_{2}$. Par le premier théorème d'isomorphisme, le groupe ${\rm Im}(f)$ est isomorphe au quotient de $\mathbb{F}_{2^{n}}$ par $\mathbb{F}_{2}$. En particulier, ${\rm card(Im}(f))=2^{n-1}$.
\\
\\Soit $x \in  \mathbb{F}_{2^{n}}$. $T(f(x))=T(x)+T(x^{2})=\overline{2}T(x)=\overline{0}$. Donc, ${\rm Im}(f) \subset {\rm Ker}(T)$. Or, ${\rm Ker}(T)$ est constitué des racines sur $\mathbb{F}_{2^{n}}$ du polynôme $P(Y)=Y+Y^{2}+Y^{4}+\ldots+Y^{2^{n-1}}$. Par conséquent, on a l'inégalité ${\rm card(Ker}(T)) \leq {\rm deg}(P)=2^{n-1}$. Ainsi, ${\rm Im}(f)={\rm Ker}(T)$.
\\
\\Considérons maintenant le polynôme $Q(X)=X^{2}+aX+b$. Comme $a$ est inversible, $Q$ est scindé sur $\mathbb{F}_{2^{n}}$ si et seulement si l'équation $(a^{-1}x)^{2}+a^{-1}x=ba^{-2}$ a des solutions sur $\mathbb{F}_{2^{n}}$. En posant $y=a^{-1}x$, on constate que cela est équivalent à $y^{2}+y=ba^{-2}$ a des solutions. Donc, $Q$ est scindé sur $\mathbb{F}_{2^{n}}$ si et seulement si $ba^{-2} \in {\rm Im}(f)={\rm Ker}(T)$, c'est-à-dire si et seulement si $T(ba^{-2})=\overline{0}$.

\end{proof}

\begin{theorem}[Théorème de la base normale primitive de Davenport, \cite{D}]
\label{419}

Soient $p$ un nombre premier et $m \geq 1$. Il existe un générateur $\alpha$ du groupe $(\mathbb{F}_{p^{m}})^{*}$ tel que $(\alpha,\alpha^{p},\alpha^{p^{2}},\ldots,\alpha^{p^{m-1}})$ est une base du $(\mathbb{Z}/p\mathbb{Z})$-espace vectoriel $\mathbb{F}_{p^{m}}$.

\end{theorem}

\noindent On peut maintenant obtenir la minoration manquante.

\begin{theorem}
\label{420}

Soit $n>1$ un entier naturel. On a : $\ell_{\mathbb{F}_{2^{n}}} \geq 3(2^{n}-1)$.

\end{theorem}

\begin{proof}

Par le théorème de la base normale primitive de Davenport, il existe un générateur $\alpha$ du groupe $(\mathbb{F}_{2^{n}})^{*}$ tel que $(\alpha,\alpha^{2},\alpha^{2^{2}},\ldots,\alpha^{2^{n-2}},\alpha^{2^{n-1}})$ est une base du $(\mathbb{Z}/2\mathbb{Z})$-espace vectoriel $\mathbb{F}_{2^{n}}$. Posons $T(X):=\sum_{i=0}^{n-1} X^{2^{i}}$. L'élément $\alpha$ n'est pas racine de $T$, car sinon $(\alpha,\alpha^{2},\ldots,\alpha^{2^{n-1}})$ serait liée sur $\mathbb{Z}/2\mathbb{Z}$. Par la proposition \ref{412}, $X^{2}+\alpha^{-1}X+\overline{1}$ n'a pas de racine sur $\mathbb{F}_{2^{n}}$ (car $T((\alpha^{-1})^{-2})=T(\alpha^{2})=T(\alpha) \neq \overline{0}$). En particulier, $\alpha^{-1}$ n'est racine d'aucun trinôme de la forme $X^{2l}+X^{l+1}+\overline{1}$. De plus, $\alpha^{-1}$ est d'ordre $2^{n}-1$ dans $((\mathbb{F}_{2^{n}})^{*},\times)$. Ainsi, par le théorème \ref{42}, la solution $\alpha^{-1}$-trinomiale minimale de $(E_{\mathbb{F}_{2^{n}}})$ est irréductible de taille $3(2^{n}-1)$, ce qui donne l'inégalité souhaitée.

\end{proof}

Notons que la borne $3({\rm card}(\mathbb{F}_{2^{n}})-1)$ est optimale puisqu'elle est atteinte pour $\mathbb{F}_{4}$ (voir par exemple \cite{M} Théorème 2.8). En revanche, il n'y a en général pas égalité. On a par exemple obtenu informatiquement $\ell_{\mathbb{F}_{8}}=25>21$ (en particulier le 25-uplet ci-dessous est une $\lambda$-quiddité irréductible sur $\mathbb{F}_{8}=\frac{(\mathbb{Z}/2\mathbb{Z})[x]}{\left\langle x^3+x+\overline{1}\right\rangle}$ :
\[(x^{2}+x+\overline{1},x,x,x,x^{2}+x+\overline{1},x^{2}+x,x^{2},x^{2}+x,x,x^{2}+\overline{1},x^{2}+x,x^{2}+x+\overline{1},x,x+\overline{1},\]
\[x+\overline{1},x+\overline{1},x+\overline{1},x,x^{2}+x+\overline{1},x^{2}+x,x^{2}+\overline{1},x,x^{2}+x,x^{2},x^{2}+x)).\]

\begin{example}
{\rm On a $\ell_{\mathbb{F}_{256}} \geq 765$.
}
\end{example}

On possède donc une minoration valable pour tous les corps de caractéristique 2. Cela dit, celle-ci repose en partie sur le théorème de Davenport qui est loin d'être un résultat élémentaire. Il serait donc intéressant d'avoir une preuve plus simple pour certaines valeurs particulières de $n$, notamment dans les cas où $2^{n}-1$ est premier. Pour mener à bien cette tâche, on va utiliser un théorème dont la démonstration est beaucoup plus simple que celui de Davenport.

\begin{theorem}[Théorème de Szymiczek, \cite{Sz} Théorème 1]
\label{411}

Soient $q>1$ la puissance d'un nombre premier, $\mathcal{G}_{q}$ l'ensemble des générateurs de $\mathbb{F}_{q}^{*}$, $m$ un entier naturel et $\mu$ la fonction de Möbius. On a :
\[\sum_{g \in \mathcal{G}_{q}} g^{m}=\overline{\mu(e)\frac{\varphi(q-1)}{\varphi(e)}},~~~~~~{\rm avec}~e=\frac{q-1}{{\rm pgcd}(m,q-1)}.\]

\end{theorem}

\noindent Ce résultat va nous permettre de démontrer le théorème suivant :

\begin{theorem}
\label{413}

Soit $n>1$ un entier naturel impair tel que $2^{n}-1$ ne possède pas de facteur carré. 
\[\ell_{\mathbb{F}_{2^{n}}} \geq 3(2^{n}-1).\]

\end{theorem}

\begin{proof}

Soient $n \neq 1$ un entier naturel impair tel que $2^{n}-1$ ne possède pas de facteur carré et $\mathcal{G}_{2^{n}}$ l'ensemble des générateurs de $(\mathbb{F}_{2^{n}})^{*}$. 
\\
\\Supposons par l'absurde que, pour tout $\alpha \in \mathcal{G}_{2^{n}}$, la solution $\alpha$-trinomiale minimale de $(E_{\mathbb{F}_{2^{n}}})$ soit réductible. Soit $\alpha \in \mathcal{G}_{2^{n}}$. Par le théorème \ref{42}, il existe $l$ tel que $\alpha^{2l}+\alpha^{l+1}+\overline{1}=\overline{0}$. En particulier, $\alpha^{l}$ est une racine sur $\mathbb{F}_{2^{n}}$ de $P_{\alpha}(X)=X^{2}+\alpha X+\overline{1}$. 
\\
\\Donc, par la proposition \ref{412}, $\sum_{i=0}^{n-1} (\alpha^{-2})^{2^{i}}=\sum_{i=0}^{n-1} (\alpha^{-1})^{2^{i}}=\overline{0}$. Puisque cette égalité est vérifiée pour tous les générateurs, on a :
\begin{eqnarray*}
\overline{0} &=& \sum_{\alpha \in \mathcal{G}_{2^{n}}} \sum_{i=0}^{n-1} (\alpha^{-1})^{2^{i}} \\
             &=& \sum_{i=0}^{n-1}\left(\sum_{\alpha \in \mathcal{G}_{2^{n}}} (\alpha^{-1})^{2^{i}}\right) \\
						 &=& \sum_{i=0}^{n-1}\left(\sum_{\alpha \in \mathcal{G}_{2^{n}}} \alpha^{2^{i}}\right)~~({\rm car}~\{\alpha^{-1},~\alpha \in \mathcal{G}_{2^{n}}\}=\mathcal{G}_{2^{n}}) \\
						 &=& \sum_{i=0}^{n-1}\overline{\mu\left(\frac{2^{n}-1}{{\rm pgcd}(2^{i},2^{n}-1)}\right) \frac{\varphi(2^{n}-1)}{\varphi \left(\frac{2^{n}-1}{{\rm pgcd}(2^{i},2^{n}-1)}\right)}}~~({\rm Th\acute{e}or\grave{e}me~de~Szymiczek)} \\
             &=& \sum_{i=0}^{n-1}\overline{\mu(2^{n}-1) \frac{\varphi(2^{n}-1)}{\varphi(2^{n}-1)}} \\		
						 &=& \sum_{i=0}^{n-1}\overline{\mu(2^{n}-1)} \\
						 &=& \sum_{i=0}^{n-1}\overline{1}~~({\rm car}~2^{n}-1~{\rm est~sans~facteur~carr\acute{e}~{\rm et}~\overline{\mu(2^{n}-1)} \in \mathbb{Z}/2\mathbb{Z})} \\
             &=& \overline{n} \\
						 &=& \overline{1}~~({\rm car}~n~{\rm est~impair}).\\
\end{eqnarray*}		

\noindent On arrive ainsi à une absurdité. Donc, il existe $\beta \in \mathcal{G}_{2^{n}}$ tel que la solution $\beta$-trinomiale minimale est irréductible. Par le théorème \ref{42}, cette dernière est de taille $3(2^{n}-1)$, ce qui donne l'inégalité souhaitée.

\end{proof}

\begin{examples}
{\rm \begin{itemize}
\item $2^{9}=512$ et $2^{9}-1=511=7 \times 73$. Puisque 9 est impair et puisque 511 n'a pas de facteur carré, $\ell_{\mathbb{F}_{512}} \geq 1533$;
\item $2^{15}=32~768$ et $2^{15}-1=7 \times 31 \times 151$. Puisque 15 est impair et puisque $2^{15}-1$ n'a pas de facteur carré, $\ell_{\mathbb{F}_{32~768}} \geq 98~301$.
\end{itemize}
}
\end{examples}

\begin{corollary}
\label{414}

Si $M_{p}:=2^{p}-1$ est un nombre premier de Mersenne alors $\ell_{\mathbb{F}_{M_{p}+1}} \geq 3M_{p}$.

\end{corollary}

\begin{proof}

Comme $M_{p}$ est premier, $p$ est premier (voir par exemple \cite{CP} Théorème 1.3.1). Si $p=2$ alors $M_{p}=3$ et, par le théorème \ref{42}, la solution $X$-trinomiale minimale de $(E_{\mathbb{F}_{4}})$ est irréductible ce qui donne $\ell_{\mathbb{F}_{4}} \geq 9$. Supposons maintenant $p$ impair. Puisque $M_{p}$ est premier, il est sans facteur carré. Donc, par le théorème \ref{413}, $\ell_{\mathbb{F}_{M_{p}+1}} \geq 3M_{p}$.

\end{proof}

\begin{remarks}
{\rm i) On connaît actuellement une cinquantaine de nombres de Mersenne premiers (voir \cite{CP} tableau 1.2 et \cite{OEIS} A000043 et A000668) et on conjecture qu'il en existe une infinité. Pour savoir si, pour un nombre premier $p$, $M_{p}$ est premier, on dispose d'un test efficace : le test de Lucas-Lehmer (voir \cite{L} et \cite{CP} Théorème 4.2.6). On considère la suite $(v_{n})$ définie par $v_{0}=4$ et $v_{n+1}=v_{n}^{2}-2$. Si $p$ est un nombre premier impair alors $M_{p}$ est premier si et seulement si $M_{p}$ divise $v_{p-2}$. C'est notamment grâce à une version légèrement différente de ce test qu'Édouard Lucas a établi en 1876 le record (toujours d'actualité) du plus grand nombre premier de Mersenne trouvé sans avoir recours à une machine (il s'agit de $M_{127}$).
\\
\\ii) Dans le théorème \ref{413}, on utilise des nombres de Mersenne $M_{n}=2^{n}-1$ sans facteur carré. Cela pose naturellement la question de savoir si ces nombres peuvent avoir ou non des facteurs carrés. Sur ce point, on dispose du résultat important suivant. Soient $p$ et $q$ deux nombres premiers impairs. Si $p^{2}$ divise $M_{q}$ alors $2^{p-1} \equiv 1 [p^{2}]$ (voir \cite{WB} Théorème 1). Les nombres premiers $p$ vérifiant la condition $2^{p-1} \equiv 1 [p^{2}]$ sont appelés nombres premiers de Wieferich. On ne connaît que deux nombres de ce type, 1093 et 3511 et 1093 et 3511 ne divisent aucun $M_{q}$, avec $q$ premier (voir \cite{WB} Théorème 2). Pour d'autres informations sur les diviseurs carrés des nombres de Mersenne, on peut également consulter \cite{Ro}.
}
\end{remarks}

\section{Quelques autres familles particulières de solutions}
\label{autres}

L'objectif de cette section est de définir rapidement d'autres classes de solutions afin d'avoir d'autres résultats d’irréductibilité.

\subsection{Solutions définies par la répétition d'une même séquence}

Le but de cette partie est de constater rapidement que le procédé utilisé jusqu'ici pour définir des solutions, à savoir la répétition d'un motif de taille finie, peut bien entendu se généraliser mais peut difficilement fournir des résultats intéressants d'irréductibilité, même en choisissant les bords du motif de façon à obtenir une matrice diagonale associée à la séquence de base. Pour cela, on va utiliser la définition suivante :

\begin{definition}
\label{51}

Soient $A$ un anneau commutatif unitaire fini et $(a,b) \in A^{2}$ tel que $ab-1_{A} \in U(A)$. Soit $\alpha$ l'unique solution de l'équation $K_{3}(x,a,b)=0_{A}$. Soit $\beta$ l'unique solution de l'équation $K_{3}(a,b,x)=0_{A}$. On appelle solution $(a,b)$-quadrinomiale minimale une solution de \eqref{p} de la forme $(\alpha,a,b,\beta,\ldots,\alpha,a,b,\beta)$ et qui est de taille minimale. On appelle solution quadrinomiale minimale sur $A$ une solution $(a,b)$-quadrinomiale minimale pour un certain couple $(a,b)$.

\end{definition} 

\begin{proof}[Justification]

Cette définition nécessite quelques vérifications pour être valide. On a :
\[K_{3}(x,a,b)=xab-x-b=x(ab-1_{A})-b.\] 
\noindent Donc, $K_{3}(x,a,b)=0_{A}$ si et seulement si $x=b(ab-1_{A})^{-1}$. De même, 
\[K_{3}(a,b,x)=abx-x-a=x(ab-1_{A})-a.\] 
\noindent Donc, $K_{3}(a,b,x)=0_{A}$ si et seulement si $x=a(ab-1_{A})^{-1}$. Ainsi, $\alpha$ et $\beta$ sont bien définis.
\\
\\Comme la matrice $M_{4}(\alpha,a,b,\beta)$ est d'ordre fini dans $SL_{2}(A)/\{-Id,Id\}$, la solution $(a,b)$-quadrinomiale minimale de \eqref{p} existe toujours.

\end{proof}

\begin{examples}
{\rm 
\begin{itemize}
\item $(0_{A},0_{A},0_{A},0_{A})$ est une solution quadrinomiale minimale de \eqref{p}.
\item $A=\mathbb{Z}/8\mathbb{Z}$. $(\overline{2},\overline{2},\overline{6},\overline{6},\overline{2},\overline{2},\overline{6},\overline{6})$ est la solution $(\overline{2},\overline{6})$-quadrinomiale minimale de \eqref{p}.
\end{itemize}
}
\end{examples}

\begin{remarks}
{\rm 
i) Si $ab=2_{A}$ alors la solution $(a,b)$-quadrinomiale minimale est une solution dynomiale.
\\
\\ii) Si $A$ est un anneau local et si $a$ ou $b$ n'est pas inversible alors $ab-1_{A}$ est inversible.
}
\end{remarks}

Notons $(a_{1},\ldots,a_{m}):=(\alpha,a,b,\beta,\ldots,\alpha,a,b,\beta)$ et supposons que $(a_{1},\ldots,a_{m})$ est une solution. Comme dans le théorème \ref{42}, on a $(a_{1},\ldots,a_{m})$ réductible si et seulement s'il existe $1 \leq j \leq \frac{m-2}{2}$ et $i$ tels que $K_{j}(a_{i},\ldots,a_{j+i-1})=\pm 1_{A}$ (avec les indices des $a_{k}$ vus modulo $m$). Autrement dit, l'étude de la réductibilité dépend du calcul de tous les sous-continuants de $(a_{1},\ldots,a_{m})$ possibles. Cela dit, contrairement aux solutions trinomiales, il semble beaucoup plus difficile d'avoir des formules complètes pour ces derniers, comme le laisse entrevoir le calcul ci-dessous :
\begin{eqnarray*}
M &=& M_{4n}(b,\beta,\alpha,a,\ldots,b,\beta,\alpha,a) \\
  &=& M_{4}(b,\beta,\alpha,a)^{n} \\
                                                 &=& \left(\begin{pmatrix}
     ab(ab-1_{A})^{-1}[ab(ab-1_{A})^{-1}-2_{A}]+1_{A}-ab     & a(ab-1_{A})^{-1}-a[ab(ab-1_{A})^{-2}-1_{A}] \\
     b[ab(ab-1_{A})^{-2}-1_{A}]-b(ab-1_{A})^{-1}                          & 1_{A}-ab(ab-1_{A})^{-2}
\end{pmatrix}\right)^{n}. \\
\end{eqnarray*}

\noindent Par conséquent, il paraît difficile d'obtenir des conditions nécessaires et/ou suffisantes d'irréductibilité des solutions quadrinomiales qui soient réellement exploitables. À fortiori, il semble hautement improbable que l'on puisse obtenir des résultats intéressants en généralisant le procédé, c'est-à-dire en augmentant la taille du motif à répéter. Par conséquent, il est inutile de définir des solutions pentanomiales, hexanomiales, heptanomiales, etc. Par ailleurs, l'augmentation progressive de la taille du motif rend la possibilité d'avoir des solutions irréductibles de plus en plus faible. En effet, plus la taille de la séquence de base est importante plus celle-ci risque de contenir un sous-continuant valant $\pm 1_{A}$. De plus, une solution construite par ce procédé a une taille qui est un multiple de la taille de la séquence de base. Par conséquent, on se rapproche de plus en plus vite de la taille maximale des solutions irréductibles. Par exemple, pour $A=\mathbb{Z}/7\mathbb{Z}$, on a $\ell_{A}=9$. Donc, une solution quadrinomiale irréductible ne pourrait être que de taille 4 ou 8 et une solution pentanomiale irréductible ne pourrait être que de taille 5, ce qui en fait en réalité une solution quelconque puisqu'on perd la répétition du motif de base.
\\
\\ \indent On pourrait évidemment envisager des motifs différents, c'est-à-dire des motifs qui ne conduisent pas à des matrices diagonales. Cependant, les perspectives d'avoir des résultats intéressants semblent plutôt faibles. Pour s'en convaincre, on peut par exemple regarder le cas de $A=\mathbb{Z}/43\mathbb{Z}$. Pour avoir une base de comparaison, on commence par donner la liste, obtenue informatiquement, des triplets $[[a,b],m]$, avec $2 \leq a \leq 21$, pour lesquels la solution $(\overline{a},\overline{b})$-quadrinomiale minimale est irréductible (de taille $m$).
\\
\\(((2,4),12), ((2,19),12), ((3,3),28), ((3,12),28), ((3,14),28), ((3,17),12), ((3,20),28), ((3,-16),12), ((3,-5),28), ((4,2),12), ((4,-12),12), ((5,12),28), ((5,17),28),((5,-3),28), ((6,6),28), ((6,7),28), ((6,19),28), ((6,-20),28), ((7,6),28), ((7,-9),28), ((9,9),12), ((9,14),28), ((9,18),28), ((9,19),28), ((9,20),12), ((9,-14),28), ((9,-9),28), ((9,-7),28), ((10,18),12), ((10,21),12), ((10,-12),28), ((11,-20),12), ((11,-11),12), ((12,3),28), ((12,5),28), ((12,15),12), ((12,-12),28), ((12,-10),28), ((12,-4),12), ((13,13),28), ((13,15),28), ((14,3),28), ((14,9),28), ((14,15),12), ((14,19),12), ((14,-15),28), ((14,-9),28), ((15,12),12), ((15,13),28), ((15,14),12), ((15,18),28), ((15,20),28), ((15,-17),28), ((15,-15),28), ((15,-14),28), ((16,-3),12), ((17,3),12), ((17,5),28), ((17,-17),28), ((17,-15),28), ((18,9),28), ((18,10),12), ((18,15),28), ((18,18),28), ((19,2),12), ((19,6),28), ((19,-34),28), ((19, 14),12), ((19,19),28), ((19,20),28), ((20,3),28), ((20,9),12), ((20,15),28), ((20,19),28), ((20,-11),12), ((20,-6),28), ((21,10),12)).
\\
\\On constate que, malgré l'absence de résultats théoriques, on obtient de nombreuses solutions quadrinomiales minimales irréductibles, même si ces dernières ne permettent pas d'avoir d'informations sur $\ell_{A}$. En revanche, si on cherche informatiquement, pour $\mathbb{Z}/43\mathbb{Z}$, les solutions irréductibles de la forme $(\overline{a},\overline{a},\overline{a}^{-1},\overline{a}^{-1},\ldots,\overline{a},\overline{a},\overline{a}^{-1},\overline{a}^{-1})$ et minimales pour cette propriété on obtient une liste vide.
Si on s'intéresse aux solutions minimales de le forme $(\overline{a},\overline{a},\overline{-a},\overline{-a},\ldots,\overline{a},\overline{a},\overline{-a},\overline{-a})$ (avec $1 \leq a \leq 21$) on récupère la liste très réduite suivante : $(((4,-4),8),((18,-18),28),((20,-20),12))$.

\subsection{Solutions $x$-quasi-monomiales}

On va maintenant s'intéresser à une nouvelle famille de solutions, intimement liée à celle des solutions monomiales minimales, qui n'est pas définie par la répétition stricte d'un motif fixé.

\begin{definition}
\label{52}

Soient $A$ un anneau commutatif unitaire fini et $x \in A$. Une solution $x$-quasi-monomiale est un $l$-uplet solution de \eqref{p} de la forme $(a,x,\ldots,x,a)$ avec $a \in A$. Lorsque $l$ est le plus petit possible, on parlera de solution $x$-quasi-monomiale minimale de \eqref{p}. On dit qu'une solution de \eqref{p} est quasi monomiale minimale si elle est égale à une solution $x$-quasi-monomiale minimale pour un $x \in A$.

\end{definition}

Notons qu'il existe toujours une solution $x$-quasi-monomiale pour tout $x \in A$. Il suffit en effet de considérer la solution $(x,A)$-monomiale minimale. De plus, pour tout $x \in A$, il existe une unique solution $x$-quasi-monomiale minimale. En effet, 
puisqu'il existe des solutions $x$-quasi-monomiales, l'ensemble des tailles de celles-ci est une partie non vide de $\mathbb{N}^{*}$. Ce dernier possède donc un plus petit élément $m$ qui est supérieur à 3. Considérons $(a,x,\ldots,x,a)$, une solution $x$-quasi-monomiale de taille $m$. Il existe un unique $\epsilon$ dans $\{-1_{A},1_{A}\}$ tel que $K_{m-2}(x,\ldots,x)=\epsilon$ et $a=\epsilon K_{m-3}(x,\ldots,x)$.

\begin{examples}
{\rm On donne ci-dessous quelques exemples de solutions $x$-quasi-monomiales minimales.
\begin{itemize}
\item $(\overline{6},\overline{3},\overline{3},\overline{6})$ est la solution $\overline{3}$-quasi-monomiale minimale de $(E_{\mathbb{Z}/9\mathbb{Z}})$;
\item $(\overline{8},\overline{3},\overline{3},\overline{3},\overline{8})$ est la solution $\overline{3}$-quasi-monomiale minimale de $(E_{\mathbb{Z}/10\mathbb{Z}})$;
\item $(\overline{7},\overline{4},\overline{4},\overline{4},\overline{4},\overline{7})$ est la solution $\overline{4}$-quasi-monomiale minimale de $(E_{\mathbb{Z}/21\mathbb{Z}})$.
\end{itemize}
}
\end{examples}

\noindent L'intérêt de ces solutions réside dans le résultat suivant.

\begin{proposition}
\label{53}

Soient $A$ un anneau commutatif unitaire et $x \in A-\{0_{A}\}$.
\\
\\i) La solution $x$-quasi-monomiale minimale de \eqref{p} est irréductible.
\\
\\ii) La solution $(x,A)$-monomiale minimale de \eqref{p} est irréductible si et seulement si elle est égale à la solution $x$-quasi-monomiale minimale de \eqref{p}.

\end{proposition}

\begin{proof}

i) Soit $m$ la taille de la solution $x$-quasi-monomiale minimale de \eqref{p}. Il existe $a \in A$ tel que cette dernière est égale à $(a,x,\ldots,x,a) \in A^{m}$. Supposons par l'absurde que la solution $x$-quasi-monomiale minimale de \eqref{p} est réductible. Il existe deux solutions de \eqref{p}, $(b_{1},\ldots,b_{l})$ et $(c_{1},\ldots,c_{l'})$, telles que $l,l' \geq 3$ et :
\[(a,x,\ldots,x,a) \sim (b_{1},\ldots,b_{l}) \oplus (c_{1},\ldots,c_{l'})=(b_{1}+c_{l'},b_{2},\ldots,b_{l-1},b_{l}+c_{1},c_{2},\ldots,c_{l'-1}).\]

\noindent Nécessairement, $(b_{1},\ldots,b_{l})$ ou bien $(c_{1},\ldots,c_{l'})$ est de la forme $(y,x,\ldots,x,z)$. Par la proposition \ref{22}, $y=z$. Comme $l, l'<m$ (puisque $l+l'=m+2$ et $l,l' \geq 3$), cela contredit la minimalité de la solution. Donc, la solution $x$-quasi-monomiale minimale de \eqref{p} est irréductible.
\\
\\ii) Si la solution $(x,A)$-monomiale minimale de \eqref{p} est égale à la solution $x$-quasi-monomiale minimale de \eqref{p} alors, par i), elle est irréductible. Supposons maintenant que la solution $(x,A)$-monomiale minimale de \eqref{p} n'est pas égale à la solution $x$-quasi-monomiale minimale. La taille de la solution $x$-quasi-monomiale minimale est nécessairement strictement inférieure à celle de la solution $(x,A)$-monomiale minimale, puisque cette dernière est une solution $x$-quasi-monomiale. Par conséquent, on peut utiliser la solution $x$-quasi-monomiale minimale pour réduire la solution $(x,A)$-monomiale minimale.

\end{proof}

\begin{example}
{\rm $(\overline{7},\overline{4},\overline{4},\overline{4},\overline{4},\overline{7})$ est une $\lambda$-quiddité irréductible de $(E_{\mathbb{Z}/21\mathbb{Z}})$.
}
\end{example}

Les solutions quasi-monomiales minimales non nulles possèdent donc la propriété assez intéressante, et relativement rare, d'être toujours irréductibles quel que soit l'anneau considéré. Cependant, leur intérêt demeure modeste. En effet, la taille de la solution $x$-quasi-monomiale minimale est systématiquement inférieure à celle de la solution $(x,A)$-monomiale minimale, ce qui brise d'une façon nette tout espoir de les utiliser pour obtenir des améliorations pour les bornes de $\ell_{A}$. De plus, dans le cas des corps finis, les familles des solutions monomiales minimales et quasi-monomiales minimales sont rigoureusement identiques. Par conséquent, les solutions quasi-monomiales minimales doivent plus être vues comme des compagnons de route des solutions monomiales minimales que comme une famille à part entière de $\lambda$-quiddités.

\subsection{Solutions $x$-polarisées}

On va maintenant définir une dernière classe de $\lambda$-quiddités, liée aux solutions monomiales minimales.

\begin{definition}
\label{53}

Soient $A$ un anneau fini de caractéristique différente de 2 et $x \in A-\{0_{A}\}$. Une solution $x$-polarisée est une solution de \eqref{p} de taille $2l$ (avec $l \in \mathbb{N}^{*}$) de la forme $(\underbrace{x,\ldots,x}_{l},\underbrace{-x,\ldots,-x}_{l})$. Lorsque $l$ est le plus petit possible, on parlera de solution $x$-polarisée minimale de \eqref{p}.

\end{definition}

Notons qu'il existe une solution $x$-polarisée pour tout $x \in A-\{0_{A}\}$. Il suffit en effet de concaténer la solution $x$-monomiale minimale avec la solution $-x$-monomiale minimale. En revanche, ce procédé ne donne pas en général la solution $x$-polarisée minimale, comme l'illustre l'exemple suivant. Prenons $A=\mathbb{Z}/6\mathbb{Z}$. La solution $\pm \overline{2}$-monomiale minimale de \eqref{p} est de taille 6 et la solution $\overline{2}$-polarisée minimale de \eqref{p} est également de taille 6.

\begin{examples}
{\rm \begin{itemize}
\item $(\overline{3},\overline{3},\overline{3},\overline{5},\overline{5},\overline{5})$ est la solution $\overline{3}$-polarisée minimale de $(E_{\mathbb{Z}/8\mathbb{Z}})$.
\item $(X,X,X,X,-X,-X,-X,-X)$ est la solution $X$-polarisée minimale de $(E_{\mathbb{F}_{9}})$, avec $\mathbb{F}_{9}=\frac{(\mathbb{Z}/3\mathbb{Z})[X]}{\left\langle X^{2}+\overline{1}\right\rangle}$.
\end{itemize}
}
\end{examples}

\noindent Pour étudier ces solutions, on aura besoin du petit résultat calculatoire ci-dessous :

\begin{lemma}
\label{541}

Soient $A$ un anneau fini de caractéristique différente de 2, $x \in A-\{0_{A}\}$ et $l \in \mathbb{N}^{*}$.
\begin{eqnarray*}
M &=& M_{2l}(\underbrace{x,\ldots,x}_{l},\underbrace{-x,\ldots,-x}_{l}) \\
  &=& (-1_{A})^{l}\begin{pmatrix}
    K_{l}(x,\ldots,x)^{2}+K_{l-1}(x,\ldots,x)^{2} & -x K_{l-1}(x,\ldots,x)^{2} \\
    -x K_{l-1}(x,\ldots,x)^{2}  & K_{l-1}(x,\ldots,x)^{2}+K_{l-2}(x,\ldots,x)^{2} 
   \end{pmatrix}.
\end{eqnarray*}	

\end{lemma}

\begin{proof}

\begin{eqnarray*}
M &=& M_{2l}(\underbrace{x,\ldots,x}_{l},\underbrace{-x,\ldots,-x}_{l}) \\
  &=& M_{l}(-x,\ldots,-x)M_{l}(x,\ldots,x) \\
	&=& \begin{pmatrix}
    K_{l}(-x,\ldots,-x) & -K_{l-1}(-x,\ldots,-x) \\
    K_{l-1}(-x,\ldots,-x)  & -K_{l-2}(-x,\ldots,-x) 
   \end{pmatrix}\begin{pmatrix}
    K_{l}(x,\ldots,x) & -K_{l-1}(x,\ldots,x) \\
    K_{l-1}(x,\ldots,x)  & -K_{l-2}(x,\ldots,x) 
   \end{pmatrix} \\
	          &=& (-1_{A})^{l}\begin{pmatrix}
    K_{l}(x,\ldots,x) & K_{l-1}(x,\ldots,x) \\
    -K_{l-1}(x,\ldots,x)  & -K_{l-2}(x,\ldots,x) 
   \end{pmatrix}\begin{pmatrix}
    K_{l}(x,\ldots,x) & -K_{l-1}(x,\ldots,x) \\
    K_{l-1}(x,\ldots,x)  & -K_{l-2}(x,\ldots,x) 
   \end{pmatrix}\\
	          & & ~~({\rm proposition}~\ref{13}) \\
	          &=& (-1_{A})^{l} \begin{pmatrix}
     a  & b \\
     b  & c 
   \end{pmatrix}.
\end{eqnarray*}						
						
\noindent avec $a=K_{l}(x,\ldots,x)^{2}+K_{l-1}(x,\ldots,x)^{2}$, $c=K_{l-1}(x,\ldots,x)^{2}+K_{l-2}(x,\ldots,x)^{2}$ et 
\begin{eqnarray*}
b &=& -K_{l}(x,\ldots,x)K_{l-1}(x,\ldots,x)-K_{l-2}(x,\ldots,x)K_{l-1}(x,\ldots,x) \\
  &=& -K_{l-1}(x,\ldots,x)(K_{l}(x,\ldots,x)+K_{l-2}(x,\ldots,x)) \\
	&=& -K_{l-1}(x,\ldots,x)(xK_{l-1}(x,\ldots,x)-K_{l-2}(x,\ldots,x)+K_{l-2}(x,\ldots,x)) \\
	&=& -x K_{l-1}(x,\ldots,x)^{2}.
\end{eqnarray*}

\end{proof}

\noindent On peut maintenant étudier l'irréductibilité de ces solutions sur les corps finis.

\begin{proposition}
\label{54}

Soient $\mathbb{K}$ un corps fini de caractéristique différente de 2, $x \in \mathbb{K}^{*}$ et $m$ la taille de la solution $(x,\mathbb{K})$-monomiale minimale de $(E_{\mathbb{K}})$. La solution $x$-polarisée minimale de $(E_{\mathbb{K}})$ est réductible de taille $2m$.

\end{proposition}

\begin{proof}

Notons $2l$ la taille de la solution $x$-polarisée minimale de $(E_{\mathbb{K}})$. Il existe $\epsilon \in \{-1_{\mathbb{K}},1_{\mathbb{K}}\}$ tel que :
\begin{eqnarray*}
\epsilon Id &=& M_{2l}(\underbrace{x,\ldots,x}_{l},\underbrace{-x,\ldots,-x}_{l}) \\
            &=& (-1_{\mathbb{K}})^{l}\begin{pmatrix}
    K_{l}(x,\ldots,x)^{2}+K_{l-1}(x,\ldots,x)^{2} & -x K_{l-1}(x,\ldots,x)^{2} \\
    -x K_{l-1}(x,\ldots,x)^{2}  & K_{l-1}(x,\ldots,x)^{2}+K_{l-2}(x,\ldots,x)^{2} 
   \end{pmatrix} \\
	          & & ({\rm lemme~pr\acute{e}c\acute{e}dent}).
\end{eqnarray*}

\noindent Comme $\mathbb{K}$ est intègre (car $\mathbb{K}$ est un corps) et $x \neq 0_{\mathbb{K}}$, on a nécessairement $K_{l-1}(x,\ldots,x)=0_{\mathbb{K}}$. Puisque $K_{l}(x,\ldots,x)=xK_{l-1}(x,\ldots,x)-K_{l-2}(x,\ldots,x)$, on a $K_{l}(x,\ldots,x)=-K_{l-2}(x,\ldots,x)$. Posons $y:=K_{l}(x,\ldots,x)$. Comme $M_{l}(x,\ldots,x)=\begin{pmatrix}
    y & 0_{\mathbb{K}} \\
    0_{\mathbb{K}} & y 
   \end{pmatrix} \in SL_{2}(\mathbb{K})$ (par la proposition \ref{12}), on a $y^{2}={\rm det}(M_{l}(x,\ldots,x))=1_{\mathbb{K}}$ et donc $y=\pm 1_{\mathbb{K}}$ (puisque $\mathbb{K}$ est intègre). 
\\
\\Ainsi, $M_{l}(x,\ldots,x)=\pm Id$, ce qui implique que $l$ est un multiple de $m$. Or $2l \leq 2m$. Donc, on a nécessairement $2l=2m$ et, en utilisant la solution $(x,\mathbb{K})$-monomiale minimale, on peut réduire la solution $x$-polarisée minimale.

\end{proof}

\begin{examples}
{\rm \begin{itemize}
\item La solution $\overline{4}$-monomiale minimale de $(E_{\mathbb{Z}/7\mathbb{Z}})$ est de taille 4 (puisque $\overline{4}^{2}=\overline{2}$). Ainsi, la solution $\overline{4}$-polarisée minimale de $(E_{\mathbb{Z}/7\mathbb{Z}})$ est réductible de taille 8.
\item La solution $X$-monomiale minimale de $(E_{\mathbb{F}_{25}})$ $\left({\rm avec}~\mathbb{F}_{25}=\frac{\mathbb{Z}/5\mathbb{Z}[X]}{\left\langle X^{2}+X+\overline{1}\right\rangle}\right)$ est de taille 12. Ainsi, la solution $X$-polarisée minimale de $(E_{\mathbb{F}_{25}})$ est réductible de taille 24.
\\
\end{itemize}
}
\end{examples}

\noindent \textbf{Remerciements :} Je remercie un rapporteur anonyme des Annales Mathématiques du Québec dont les remarques et suggestions ont contribué d'une façon notable à l'amélioration de ce texte.

\appendix

\section{Carrés dans quelques corps finis}
\label{A}

On se place dans $\mathbb{F}_{9}=\frac{(\mathbb{Z}/3\mathbb{Z})[X]}{\left\langle X^{2}+\overline{1}\right\rangle}$. Il y a $\frac{9+1}{2}=5$ carrés dans ce corps. On peut les obtenir de la façon suivante :
\begin{itemize}
\item $\overline{0}^{2}=\overline{0}$;
\item $\overline{1}^{2}=\overline{1}$;
\item $X^{2}=-\overline{1}$;
\item $(X+\overline{1})^{2}=X^{2}+\overline{2}X+\overline{1}=\overline{2}X$;	
\item $(X-\overline{1})^{2}=X^{2}-\overline{2}X+\overline{1}=X$.
\\
\end{itemize}

On se place maintenant dans $\mathbb{F}_{25}=\frac{(\mathbb{Z}/5\mathbb{Z})[X]}{\left\langle X^{2}+X+\overline{1}\right\rangle}$. Il y a $\frac{25+1}{2}=13$ carrés dans ce corps. On peut les obtenir de la façon suivante :
\begin{itemize}
\item $\overline{0}^{2}=\overline{0}$;
\item $\overline{1}^{2}=\overline{1}$;
\item $\overline{2}^{2}=\overline{-1}$;
\item $X^{2}=-X-\overline{1}$;
\item $(X+\overline{1})^{2}=X^{2}+\overline{2}X+\overline{1}=X$;	
\item $(X+\overline{2})^{2}=X^{2}+\overline{4}X+\overline{4}=\overline{3}X+\overline{3}$;	
\item $(X+\overline{3})^{2}=X^{2}+\overline{6}X+\overline{9}=\overline{3}$;
\item $(X-\overline{1})^{2}=X^{2}+\overline{3}X+\overline{1}=\overline{2}X$;	
\item $(\overline{2}X)^{2}=\overline{4}X^{2}=X+\overline{1}$;		
\item $(\overline{2}X+\overline{1})^{2}=\overline{4}X^{2}+\overline{4}X+\overline{1}=\overline{2}$;		
\item $(\overline{2}X+\overline{2})^{2}=\overline{4}X^{2}+\overline{8}X+\overline{4}=-X$;
\item $(\overline{2}X+\overline{3})^{2}=\overline{4}X^{2}+\overline{12}X+\overline{9}=\overline{3}X$;
\item $(\overline{2}X+\overline{4})^{2}=\overline{4}X^{2}+\overline{16}X+\overline{16}=\overline{2}X+\overline{2}$.
\\
\end{itemize}

On se place pour finir dans $\mathbb{F}_{27}=\frac{(\mathbb{Z}/3\mathbb{Z})[X]}{\left\langle X^{3}+X^{2}-\overline{1}\right\rangle}$. Il y a $\frac{27+1}{2}=14$ carrés dans ce corps. On peut les obtenir de la façon suivante :
\begin{itemize}
\item $\overline{0}^{2}=\overline{0}$;
\item $\overline{1}^{2}=\overline{1}$;
\item $X^{2}$;
\item $(X+\overline{1})^{2}=X^{2}-X+\overline{1}$;
\item $(X-\overline{1})^{2}=X^{2}+X+\overline{1}$;
\item $(X^{2})^{2}=X^{4}=-X^{3}+X=X^{2}+X-\overline{1}$;
\item $(X^{2}+\overline{1})^{2}=X^{4}-X^{2}+\overline{1}=X$;
\item $(X^{2}-\overline{1})^{2}=X^{4}+X^{2}+\overline{1}=-X^{2}+X$; 
\item $(X^{2}+X)^{2}=X^{4}-X^{3}+X^{2}=X+\overline{1}$; 
\item $(X^{2}-X)^{2}=X^{4}+X^{3}+X^{2}=X^{2}+X$;
\item $(X^{2}+X+\overline{1})^{2}=X^{4}-X^{3}-X+\overline{1}=-X^{2}-\overline{1}$;
\item $(X^{2}-X+\overline{1})^{2}=X^{4}+X^{3}+X+\overline{1}=-X+\overline{1}$;
\item $(X^{2}+X-\overline{1})^{2}=X^{4}+X=X^{2}-X-\overline{1}$;
\item $(X^{2}-X-\overline{1})^{2}=X^{4}+X^{3}-X^{2}-X+\overline{1}=-X^{2}+\overline{1}$.
\\
\end{itemize}

\section{Application informatique du théorème \ref{42}}
\label{C}

Soit $p$ un nombre premier supérieur à 5. Le théorème \ref{42} nous permet de relier l'irréductibilité des solutions trinomiales minimales de $(E_{\mathbb{Z}/p\mathbb{Z}})$ à l'étude des racines des polynômes $X^{2l}+X^{l+1}-\overline{1}$ ou $X^{2l}-X^{l+1}-\overline{1}$. Cette tâche peut être facilement menée informatiquement. 
\\
\\\indent Le programme Maxima racinetrino(p) donne la liste des éléments $i$ de $\mathbb{Z}$ compris entre 2 et $\frac{p-1}{2}$ qui ne sont pas racines modulo $p$ d'un polynôme $X^{2l}+X^{l+1}-1$ ou $X^{2l}-X^{l+1}-1$ pour $l$ compris entre 1 et $E\left[\frac{m}{6}\right]$, avec $m$ la taille de la solution $\overline{i}$-trinomiale minimale (cette taille est calculée par le programme en déterminant l'ordre de $i$ dans $(\mathbb{Z}/p\mathbb{Z})^{*}$ et en utilisant les formules théoriques). Par le théorème \ref{42}, on obtient la liste des éléments non nuls $\overline{x}$, avec $2 \leq x \leq \frac{p-1}{2}$, de $\mathbb{Z}/p\mathbb{Z}$ pour lesquels la solution $\overline{x}$-trinomiale minimale de $(E_{\mathbb{Z}/p\mathbb{Z}})$ est irréductible.
\\
\\racinetrino(p):=block([t,k,z,h,W,H,L],
\\ \phantom{BC} modulus:p,
\\ \phantom{BC} H:setify(create\_list(i,i,2,(p-1)/2)),
\\ \phantom{BC} W:[],
\\ \phantom{BC} for i:2 thru (p-1)/2 do
\\ \phantom{BCDE} (t:0,
\\ \phantom{BCDE} k:1,
\\ \phantom{BCDE} z:rat(1),
\\ \phantom{BCDE} while t=0 and k<(p+1)/2 do
\\ \phantom{BCDEFG} (z:rat(z*i),
\\ \phantom{BCDEFG} if z=rat(1) then t:1,
\\ \phantom{BCDEFG} k:k+1),
\\ \phantom{BCDE} if t=1 then h:k-1,
\\ \phantom{BCDE} if t=0 then h:p-1,
\\ \phantom{BCDE} if divide(h,2)[2]=1 then 
\\ \phantom{BCDEFG} (for l:1 thru floor(h/2) do
\\ \phantom{BCDEFGHI} (if rat(i\textasciicircum (2*l)+i\textasciicircum (l+1)-1)=0 then W:append(W,[i]),
\\ \phantom{BCDEFGHI} if rat(i\textasciicircum (2*l)-i\textasciicircum (l+1)-1)=0 then W:append(W,[i]))),
\\ \phantom{BCDE} if divide(h,2)[2]=0 then 
\\ \phantom{BCDEFG} (for l:1 thru floor(h/4) do
\\ \phantom{BCDEFGHI} (if rat(i\textasciicircum (2*l)+i\textasciicircum (l+1)-1)=0 then W:append(W,[i]),
\\ \phantom{BCDEFGHI} if rat(i\textasciicircum (2*l)-i\textasciicircum (l+1)-1)=0 then W:append(W,[i])))),
\\ \phantom{BC} L:listify(setdifference(H,setify(W))),
\\return(L));
\\
\\Cela donne par exemple :
\\
\\racinetrino(17)=(2, 4, 5, 6), racinetrino(31)=(2, 3, 5, 6, 7, 8, 9, 10, 12, 15),
\\
\\racinetrino(43)=(2, 4, 5, 6, 7, 8, 9, 10, 11, 12, 14, 16, 18, 21),
\\
\\racinetrino(79)=(5, 7, 8, 9, 10, 12, 13, 15, 17, 18, 22, 23, 24, 26, 27, 28, 30, 31, 33, 34, 37, 38, 39),
\\
\\racinetrino(163)=(2, 3, 4, 5, 6, 8, 11, 12, 14, 15, 17, 18, 20, 21, 22, 23, 25, 26, 27, 28, 29, 30, 32, 34, 36, 37, 38, 40, 42, 44, 45, 47, 48, 49, 50, 53, 54, 58, 59, 60, 61, 62, 64, 65, 66, 67, 73, 74, 76, 77, 78, 81),
\\
\\racinetrino(389)=(2, 4, 5, 6, 7, 10, 11, 12, 13, 14, 15, 18, 20, 21, 24, 25, 26, 28, 30, 32, 35, 36, 38, 39, 40, 41, 42, 43, 44, 45, 48, 49, 50, 51, 52, 53, 54, 55, 58, 60, 61, 63, 65, 67, 68, 69, 71, 72, 74, 76, 77, 78, 79, 81, 82, 83, 85, 87, 88, 89, 91, 92, 94, 95, 96, 99, 100, 102, 103, 104, 106, 107, 108, 109, 111, 112, 114, 115, 116, 119, 120, 121, 122, 123, 125, 127, 128, 129, 130, 133, 136, 137, 140, 141, 142, 144, 145, 146, 148, 149, 151, 155, 156, 157, 158, 159, 160, 161, 162, 167, 169, 170, 171, 173, 179, 180, 181, 183, 185, 186, 187, 188, 193).
\\
\\On fournit dans le tableau ci-dessous la taille $l$ de la liste $L$ générée par le programme racinetrino.

\hfill\break

\begin{center}
\begin{tabular}{|c|c|c|c|c|c|c|c|c|c|c|c|c|}
\hline
  $p$     & 79 & 163  & 389 & 571 & 919 & 1279 & 1801 & 2039 & 2999 & 3461 & 4349 & 4951 \rule[-7pt]{0pt}{18pt} \\
	\hline
	  $l$ & 23 & 52 & 123 & 192 & 298 & 399 & 699 & 510 & 750 & 1176 & 1353 & 1780\rule[-7pt]{0pt}{18pt} \\
	\hline
	
\end{tabular}
\end{center}

\hfill\break

\section{Calcul informatique de $\ell_{\mathbb{F}_{9}}$}
\label{F}

\noindent Pour travailler sur le corps à 9 éléments sur Maxima, on utilise la commande : gf\textunderscore set\textunderscore data(3,x\textasciicircum 2+1).
\\
\\\indent On a également besoin pour ce qui va suivre du programme suivant qui permet de calculer pour des éléments $x_{1},\ldots,x_{n}$ dans $\mathbb{F}_{9}$ le continuant $K_{n}(x_{1},\ldots,x_{n})$.
\\
\\gfconti(L):=block([u,v,w],
\\ \phantom{BC} u:1,
\\ \phantom{BC} w:0,
\\ \phantom{BC} v:0,
\\ \phantom{BC} for i:1 thru length(L) do
\\ \phantom{BCDE }(w:u,
\\ \phantom{BCDE} u:gf\textunderscore sub(gf\textunderscore mult(L[i],u),v),
\\ \phantom{BCDE} v:w),
\\return(u));
\\
\\ \indent L'objectif de l'algorithme ci-dessous est de fournir une majoration pour $\ell_{\mathbb{F}_{9}}$. Pour cela, on cherche un $n$ tel que tout $n$-uplet $(x_{1},\ldots,x_{n})$ de $\mathbb{F}_{9}$ contient un sous-uplet $(x_{i},\ldots,x_{i+l-1})$ vérifiant l'égalité $K_{l}(x_{i},\ldots,x_{i+l-1})=\pm \overline{1}$. Une fois qu'on a trouvé le plus petit $n$ vérifiant cette condition, on a que toutes les solutions de $(E_{\mathbb{F}_{9}})$ de taille supérieure à $n+3$ sont réductibles, ce qui implique $\ell_{\mathbb{F}_{9}} \leq n+2$. En effet, si $(x_{1},\ldots,x_{m})$ est une solution de taille $m \geq n+3$ alors $(x_{2},\ldots,x_{m-2})$ contient un sous-uplet $(x_{i},\ldots,x_{i+l-1})$ vérifiant $K_{l}(x_{i},\ldots,x_{i+l-1})=\pm \overline{1}$ (puisque $m-3 \geq n$). En utilisant la proposition \ref{131}, on peut ensuite construire une solution permettant de réduire $(x_{1},\ldots,x_{m})$. Le but concret du programme Maxima tailleirreducneuf est de chercher pour un entier $n$ donné les $n$-uplets d'éléments de $\mathbb{F}_{9}$ ne possédant pas de sous-uplet dont le continuant vaut $\pm \overline{1}$. Concrètement, on part des 1-uplets différents de $\overline{0}, \overline{1}, \overline{-1}$ et on regarde si leur continuant vaut 1. On obtient la liste $M$ des 1-uplets ne possédant pas de sous-uplet dont le continuant vaut $\pm \overline{1}$. Ensuite, on construit à partir de $M$, la liste $N$ des candidats de taille 2. Ce sont tous les 2-uplets $(x_{1},x_{2})$ que l'on peut construire en rajoutant à la fin des 1-uplets un élément de $\mathbb{F}_{9}$ différent de $\overline{0}, \overline{1}, \overline{-1}$. On regarde ensuite si ceux-ci contiennent un sous-uplet de continuant $\pm \overline{1}$ en notant que ces valeurs ne peuvent être atteintes qu'avec $K_{2}(x_{1},x_{2})$ et $K_{1}(x_{2})$. On obtient ainsi la nouvelle liste $M$ des 2-uplets ne possédant pas de sous-uplet dont le continuant vaut $\pm \overline{1}$ à partir de laquelle on peut construire la nouvelle liste $N$ des candidats de taille 3. On continue ainsi de suite jusqu'à ce que les éléments de la liste $M$ soit de taille $n$. C'est cette liste $M$ que l'on récupère à la fin de l’exécution.
\\
\\tailleirreducneuf(n):=block([N,w,t,j,H,q,p,A,M,B],
\\ \phantom{BC} A:[x,1+x,-1+x,-x,1-x,-1-x],
\\ \phantom{BC} N:[[x],[1+x],[-1+x],[-x],[1-x],[-1-x]],
\\ \phantom{BC} p:1,
\\ \phantom{BC} M:[1],
\\ \phantom{BC} while (length(M)>0 and p<n+1) do
\\ \phantom{BCDE} (M:[],
\\ \phantom{BCDE} for h:1 thru length(N) do
\\ \phantom{BCDEFG} (w:1,
\\ \phantom{BCDEFG} j:1,
\\ \phantom{BCDEFG} while (w=1 and j<p+1) do
\\ \phantom{BCDEFGHI} (H:create\textunderscore list(N[h][k],k,j,p),
\\ \phantom{BCDEFGHI} q:gfconti(H),
\\ \phantom{BCDEFGHI} if (q=1 or q=-1 or q=2 or q=-2) then w:0,
\\ \phantom{BCDEFGHI} j:j+1),
\\ \phantom{BCDEFG} if w=1 then M:append(M,[N[h]])),
\\ \phantom{BCDE} if length(M)\# 0 then
\\ \phantom{BCDEFG} (B:[],
\\ \phantom{BCDEFG} for j:1 thru length(M) do
\\ \phantom{BCDEFGHI} (for l:1 thru length(A) do
\\ \phantom{BCDEFGHIAA} (B:append(B,[append(M[j],[A[l]])]))),
\\ \phantom{BCDEFG} N:B),
\\ \phantom{BCDE} p:p+1),
\\return(M));
\\
\\ \noindent Pour $n=12$, on obtient une liste de 124 éléments. Pour $n=13$, la liste est vide, ce qui donne $\ell_{\mathbb{F}_{9}} \leq 15$. Par ailleurs, le 15-uplet $(\overline{1}-x,x,x+\overline{1},x+\overline{1},x-\overline{1},\overline{1}-x,x,x+\overline{1},x+\overline{1},x-\overline{1},\overline{1}-x,x,x+\overline{1},x+\overline{1},x-\overline{1})$ est une $\lambda$-quidité irréductible sur $\mathbb{F}_{9}$. Ainsi, $\ell_{\mathbb{F}_{9}}=15$.

\end{document}